\documentclass[11pt,reqno]{amsart}
\usepackage{float}
\usepackage[pdftex,hyperfootnotes=false,colorlinks=false,hypertexnames=false]{hyperref}

\usepackage[english]{babel}
\usepackage{amsmath,amssymb,anyfontsize,enumerate,stmaryrd,mathtools,thmtools,tikz,
nccmath}
\usepackage[oldstyle]{libertine}%
\usepackage[T1]{fontenc}
\usepackage[final]{microtype}
\usepackage[utf8]{inputenc}
\usepackage{csquotes}
\usepackage{comment}
\makeatletter
\renewcommand*\libertine@figurestyle{LF}
\makeatother
\makeatletter
\renewcommand*\libertine@figurestyle{OsF}
\makeatother

\usepackage{mathabx}
\usepackage{tikz-cd}
\tikzcdset{scale cd/.style={every label/.append style={scale=#1},
    cells={nodes={scale=#1}}}}
\usepackage[noabbrev,capitalise]{cleveref}
\usepackage{autonum}
\usetikzlibrary{%
  matrix,%
  calc,%
  arrows%
}

\usepackage{tabularx}
\usepackage{booktabs}
\usepackage[labelfont=bf,format=plain,justification=raggedright,singlelinecheck=false]{caption}

\usepackage[sanitizesort]{glossaries}
\newglossary[slg]{symbolslist}{syi}{syg}{Index of Notation}
\makenoidxglossaries

\newglossaryentry{TT}{
 name=\ensuremath{\mathbb{T}\mathcal{T}},
 description={The category of real tori with integral structure},
 type=symbolslist
}
\newglossaryentry{TA}{
 name=\ensuremath{\mathbb{T}\mathcal{A}},
 description={The category of tropical abelian varieties},
 type=symbolslist
}
\newglossaryentry{TC}{
 name=\ensuremath{\mathbb{T}\mathcal{C}},
 description={The category of tropical curves},
 type=symbolslist
}
\newglossaryentry{Ab}{
 name=\ensuremath{Ab},
 description={The category of abelian groups},
 type=symbolslist
}

\newglossaryentry{Sigma}{
 name=\ensuremath{\Sigma},
 description={Depending on the context: A real torus with integral structure or pptav},
 type=symbolslist
}
\newglossaryentry{widecheckSigma}{
 name=\ensuremath{\widecheck{\Sigma}},
 description={The dual of Sigma},
 type=symbolslist
}

\newglossaryentry{zeta}{
 name=\ensuremath{\zeta},
 description={Polarization on $\Sigma$},
 type=symbolslist
}
\newglossaryentry{fzeta}{
 name=\ensuremath{f_\zeta},
 description={Isogeny induced by \ensuremath{\zeta}},
 type=symbolslist
}

\newglossaryentry{f}{
 name=\ensuremath{f:\Sigma_1 \rightarrow \Sigma_2},
 description={A morphism in \ensuremath{\mathbb{T}\mathcal{T}}},
 type=symbolslist
}
\newglossaryentry{widecheckf}{
 name=\ensuremath{\widecheck{f}:\widecheck{\Sigma}_2 \rightarrow \widecheck{\Sigma}_1},
 description={The dual of \ensuremath{f}},
 type=symbolslist
}

\newglossaryentry{Ker}{
 name=\ensuremath{Ker(f)_0},
 description={Kernel of \ensuremath{f} in \ensuremath{\mathbb{T}\mathcal{T}} or \ensuremath{\mathbb{T}\mathcal{A}}},
 type=symbolslist
}

\newglossaryentry{ker}{
 name=\ensuremath{\ker(f)},
 description={Kernel of \ensuremath{f} in  \ensuremath{Ab}},
 type=symbolslist
}

\newglossaryentry{zetaGamma}{
 name=\ensuremath{\zeta_\Gamma},
 description={Natural principal polarization on Jac\ensuremath{(\Gamma)}},
 type=symbolslist
}

\newglossaryentry{PhiP0}{
 name=\ensuremath{ \Phi_{P_0}},
 description={The tropical Abel-Jacobi map with reference point \ensuremath{P_0}},
 type=symbolslist
}
\newglossaryentry{Gamma}{
 name=\ensuremath{\Gamma},
 description={Tropical curve},
 type=symbolslist
}

\newglossaryentry{TE}{
 name=\ensuremath{\mathbb{T}E},
 description={Tropical elliptic curve},
 type=symbolslist
}

\newglossaryentry{Jpp}{
name=\ensuremath{J^{pp}},
description={The pptav from Construction \ref{construction_determinesplitting}},
 type=symbolslist
}

\newglossaryentry{zetapp}{
name=\ensuremath{\zeta^{pp}},
description={The principle polarization on \ensuremath{J^{pp}} from Construction \ref{construction_determinesplitting}},
type=symbolslist
}

\newglossaryentry{Qpp}{
name=\ensuremath{Q^{pp}},
 description={The quadractic form associated to \ensuremath{J^{pp}} with pp \ensuremath{\zeta^{pp}}},
 type=symbolslist
}
\theoremstyle{plain}
    \newtheorem{theorem}{Theorem}[section]
    \newtheorem{construction/theorem}[theorem]{Construction/Theorem}
    \newtheorem{corollary}[theorem]{Corollary}
    \newtheorem{lemma}[theorem]{Lemma}
    \newtheorem{proposition}[theorem]{Proposition}
    
\theoremstyle{definition}
    \newtheorem{remark}[theorem]{Remark}
    \newtheorem{remark/reference}[theorem]{Remark/Reference}
    \newtheorem{notation}[theorem]{Notation}
    \newtheorem{example}[theorem]{Example}
    \newtheorem{definition}[theorem]{Definition}
    \newtheorem{construction}[theorem]{Construction}
    \newtheorem{algorithm}[theorem]{Algorithm}
     
     \newtheorem{convention}[theorem]{Convention}
     \newtheorem{reconstruction}[theorem]{Plan}
     \newtheorem{conjecture}[theorem]{Conjecture}
     \newtheorem{q}[theorem]{Question}
     
     \newtheorem{task}[theorem]{Task}

\newenvironment{sketchproof}{\proof[Sketch of Proof.]}
 {\endproof}
\newcounter{diagram}
\crefname{diagram}{Diagram}{Diagrams}
\newcounter{savedequation}
 \newenvironment{diagram}[1][]
 {%
  \setcounter{savedequation}{\value{equation}}%
  \setcounter{equation}{\value{diagram}}%
  \crefalias{equation}{diagram}%
  \begin{equation}\begin{tikzcd}[#1]%
 }
 {%
  \end{tikzcd}\end{equation}%
  \setcounter{diagram}{\value{equation}}%
  \setcounter{equation}{\value{savedequation}}%
  \ignorespacesafterend 
 }

\newcommand\sbullet[1][.5]{\mathbin{\vcenter{\hbox{\scalebox{#1}{$\bullet$}}}}}

\DeclareMathOperator{\Coker}{Coker}
\DeclareMathOperator{\Hom}{Hom}
\DeclareMathOperator{\Prin}{Prin}
\DeclareMathOperator{\Ker}{Ker}
\DeclareMathOperator{\im}{im}
\DeclareMathOperator{\Det}{Det}
\DeclareMathOperator{\lcm}{lcm}
\DeclareMathOperator{\Div}{Div}
\DeclareMathOperator{\rk}{rk}
\DeclareMathOperator{\Jac}{Jac}
\DeclareMathOperator{\Pic}{Pic}
\DeclareMathOperator{\H1}{H_1}
\DeclareMathOperator{\Co}{C_0}
\DeclareMathOperator{\C1}{C_1}
\DeclareMathOperator{\Aut}{Aut}
\DeclareMathOperator{\Stab}{Stab}
\title{Tropical split Jacobians of curves of genus 2 II}
\author[L.~Cobigo]{Lou-Jean Leila Cobigo}
\address{L.~Cobigo: Eberhard Karls Universität Tübingen, Fachbereich Mathematik, Auf der Morgenstelle 10, 72076 Tübingen, Germany}
\email{loco@math.uni-tuebingen.de}

\begin{document}
\begin{abstract}
This paper is the second in a series of two papers which study the phenomenon of tropical split Jacobians. The first paper is a contemplative study, embedded in the broader context of exploring connections between the category of \emph{tropical abelian varieties (tav)}, $\mathbb{T}\mathcal{A}$, and the category of \emph{tropical curves}, $\mathbb{T}\mathcal{C}$. Tropical split Jacobians take on different forms depending on whether we look at them in $\mathbb{T}\mathcal{A}$ or $\mathbb{T}\mathcal{C}$: They appear either as 2 dimensional tavs that decompose into a product of two elliptic curves, or as a pair of optimal coverings. \cite{arXiv:2410.13459} examines both and then focuses on how optimal covers give rise to split Jacobians. This paper takes a different approach. Instead of looking at the phenomenon as a whole, we analyze its building blocks, a pair of elliptic curves together with a finite subgroup of their product, and how to reassemble them into a Jacobian. 
\end{abstract}
\maketitle
\tableofcontents
\glsaddall
\printnoidxglossary[type=symbolslist,sort=def]
\section{Introduction}\label{section_intro}
This paper is the second in a series of two that revolve around the phenomenon of tropical split Jacobians. 

Embedded in the broader context of exploring connections between the category of tropical abelian varieties (tav), $\mathbb{T}\mathcal{A}$, and the category of tropical curves, $\mathbb{T}\mathcal{C}$, \cite{arXiv:2410.13459} is a \emph{phenomenological analysis}: Split Jacobians may appear either in the form of 2 dimensional abelian varieties that decompose into a product of two tropical elliptic curves, or as a pair of tropical optimal coverings. 

\begin{theorem}(\cite{arXiv:2410.13459}, Theorem 61 and Theorem 66)\label{theorem_CovertoSplitting}
Let $\Gamma$ be a tropical curve of genus $2$. For an optimal pair $(\mathbb{T}E,\varphi)$, that is a pair consisting of a tropical elliptic curve $\mathbb{T}E$ and a tropical optimal cover $\varphi: \Gamma \rightarrow \mathbb{T}E$, there exists another optimal pair $(\mathbb{T}E',\varphi')$ such that the tropical Jacobian of $\Gamma$, $\Jac(\Gamma)$, splits. Moreover, the splitting of $\Jac(\Gamma)$ is given by an isogeny
  \begin{align}
      \phi:\mathbb{T}E' \bigoplus \mathbb{T}E \rightarrow \Jac(\Gamma),
  \end{align}
  whose kernel satisfies $\Jac_d(\mathbb{T}E')\cong \ker(\phi) \cong \Jac_d(\mathbb{T}E)$. In this case $\Jac(\Gamma)$ is said to be d-split and $(\mathbb{T}E,\mathbb{T}E',\ker(\phi))$ is referred to as splitting data.
 \end{theorem}
 More precisely, through push-forward and pull-back of divisors $\varphi$ and $\varphi'$ generate two exact sequences

\begin{align}
 &0 \rightarrow \mathbb{T}E' \xrightarrow{\varphi'^{*}} \Jac(\Gamma) \xrightarrow{\varphi_*}  \mathbb{T}E \rightarrow 0  \\ 
   &   0 \rightarrow \mathbb{T}E \xrightarrow{\varphi*} \Jac(\Gamma) \xrightarrow{\varphi_*'}  \mathbb{T}E' \rightarrow 0.
\end{align}
The pair of maps, $(\varphi'^{*},\varphi^{*})$ and $(\varphi'_{*},\varphi_{*})$, in turn give rise to isogenies, $\phi$ and $\tilde{\phi}$, by utilizing the coproduct, respectively the product, property of $\mathbb{T}E' \bigoplus \mathbb{T}E$ (recall that finite products coincide with finite coproducts in $\mathbb{T}\mathcal{A}$) and interact as follows: 
\begin{diagram}\label{diagram_fromctocultimately}
     \mathbb{T}E' \ar[dd,"m_d"] \ar[dr,"\varphi^{'*}",sloped] \ar[r,"\iota_1", hook] & \mathbb{T}E' \bigoplus \mathbb{T}E \ar[d,dashed,"{\phi}" description] &\mathbb{T}E \ar[dl,"\varphi^{*}",sloped] \ar[l,"\iota_2", hook', swap] \ar[dd,"m_d"]\\
        & \ar[dl,"\varphi'_{*}",sloped] \Jac(\Gamma) \ar[d,dashed,"{\tilde{\phi}}" description] \ar[dr,"\varphi_{*}" ] &  \\
        \mathbb{T}E'  & \ar[l, "p_1"] \mathbb{T}E' \bigoplus \mathbb{T}E \ar[r, "p_2",swap]  &\mathbb{T}E  
\end{diagram}
 where $\iota_i$ and $p_i$ are the canonical injections, respectively projections. The diagonals are formed by our exact sequences and a small diagram-chase shows that $\tilde{\phi}\circ \phi$ is the componentwise multiplication-by-$d$ map.

This paper opens with an \emph{atomic perspective}.
\begin{q}\label{question_Intro}
    What are the basic components of split Jacobians and how are they assembled?
\end{q}
Let $\mathbb{T}E$ and $\mathbb{T}E'$ be tropical elliptic curves and $G$ a finite subgroup of their product. Question \ref{question_Intro} translates to:
\begin{task}\label{task_local}
   \emph{Setting:} $\mathbb{T}\mathcal{A}$ and $\mathbb{T}\mathcal{C}$. Determine whether $(\mathbb{T}E,\mathbb{T}E',G)$ is splitting data and, if it is, construct a tropical curve $\Gamma$ of genus $2$ whose Jacobian splits accordingly (see Section \ref{section_reconstruction} for a precise statement).
\end{task}
Task \ref{task_local} paves the way for addressing issues of a more global nature in the moduli space $A^{tr}_2$ of principally polarized tropical abelian varieties and the moduli space $M^{tr}_2$ of tropical curves of genus $2$.
\begin{task}\label{task_global}
   \emph{Setting:} $A^{tr}_2$ and $M^{tr}_2$. Analysis of the structure of the locus of split Jacobians and the locus of curves with split Jacobians. 
\end{task}
\subsection{Local: Task \ref{task_local}}
The expected end product has been described in \cite{arXiv:2410.13459} (see Theorem \ref{theorem_CovertoSplitting}).
In consequence, we adopt the following approach:

\begin{reconstruction} [see Plan \ref{Masterplan}]\label{masterplaninintro}
Given splitting data $(\mathbb{T}E,\mathbb{T}E',G)$, proceed as follows:
\begin{enumerate}
    \item Determine a splitting $\phi: \mathbb{T}E \oplus \mathbb{T}E' \rightarrow J $ and generate a diagram $D_\phi$ modelled on Diagram \ref{diagram_fromctocultimately}, 
    \item Construct $\gls{Gamma}$,
    \item Define covers $\varphi$ and $\varphi'$ of $\mathbb{T}E$ and $\mathbb{T}E'$.
\end{enumerate}
\end{reconstruction}
\subsection{Global: Task \ref{task_global}}
A natural by-product of Plan \ref{masterplaninintro} is a family of split Jacobians and a characterization of their pointwise preimage under the tropical Torelli map $t^{tr}_2$. The moduli spaces $M^{tr}_2$ and $A^{tr}_2$ offer the right setting for a change of perspective. We consider the following subsets: 
\begin{itemize}
   \item $\mathcal{Q}\subset A_2^{tr}$ the \emph{locus of split Jacobians}.
    \item $\mathbb{T}\mathcal{L}_d\subset M_2^{tr}$ the \emph{locus of curves with d-split Jacobians}. 
\end{itemize}
We think of the first as a Schottky type problem for split Jacobians, whereas the second is a tropical analogue of the classical locus of curves with $(d,d)$-split Jacobians, which is known to be an irreducible 2 dimensional subvariety of $\mathcal{M}_2$ (e.g. \cite{shaska2024machinelearningmodulispace} or \cite{zbMATH05564780}).

Concretely, we are interested in
\begin{itemize}
    \item characterizing the intersection of $\mathcal{Q}$ with the boundary of $ A_2^{tr} $.
    \item determining $\mathbb{T}\mathcal{L}_d$ has a "tropical" structure that in some sense reflects the algebraic one.
\end{itemize}
\subsection{Context}
The present work has its roots in classical algebraic geometry, as did the first paper in this series, \cite{arXiv:2410.13459}.
In \cite{MR1085258} Frey and Kani describe how to construct curves of genus 2 covering two elliptic curves $E$ and $E'$ from an isomorphism
\begin{align}
    \alpha: E'[d] \rightarrow E[d]
\end{align}
whose graph $G_\alpha$ is isotropic with respect to the Weil pairing on $(E \times E')[d]$. The isotropy condition for $G_\alpha$ is necessary for one of the main tools in this context, Mumford' Criterion (Proposition 16.8 \cite{zbMATH03975095} or \cite{MR0282985}, p. 231), for which we develop a tropical analogue in Section \ref{subsection_isogeniesandpolarizations}: It guarantees that the candidate for the Jacobian, the quotient $(E' \times E)/G_\alpha$, can be equipped with a principal polarization (pp). Guided by their ideas, we consider an analogue starting point: We take the data contained in $\alpha$ to represent the buildings blocks of a classical split Jacobian and consider analogues of these building blocks for tropical split Jacobians. While the main steps in the construction of $\Gamma$ can be found in the classical and tropical approach (see Plan \ref{masterplaninintro}), their implementation diverges. For more details on the classical side we refer to \cite{MR1085258}.

With two separate but related settings, Tasks \ref{task_local} and \ref{task_global} need to be tackled in different ways:
Working in $\mathbb{T}\mathcal{C}$ and $\mathbb{T}\mathcal{A}$, we build on
techniques developed in \cite{MR3278571}, \cite{MR2661417},\cite{MR3717092} for handling curves and covers
and \cite{MR2457725}, \cite{röhrle2024tropicalngonalconstruction}, \cite{MR4261102}, \cite{MR2275625}, \cite{arXiv:2410.13459} and more for handling tav. Working in $M^{tr}_2$ and $A^{tr}_2$, we rely on results of \cite{MR3375652}, \cite{MR2968636}, \cite{MR2739784} and \cite{MR2641941}.

The philosophy is the one of \cite{arXiv:2410.13459}, we want to use tropical geometry to make abstraction concrete and offer both, an abstract and a constructive approach. The perspective, on the other hand, is a different one: Abstraction appears in the form of moduli spaces, i.e. we leave the "local" point of view, a focus on objects, to notice that split Jacobians arise in families whose structure we tend to grasp.

 \subsection{Results} Here we describe the main contributions of this paper.

\begin{theorem}(Theorem \ref{theorem_reconstructionsumup})\label{theorem_reconstructionsumupinintro}
    Let $(\mathbb{T}E',\mathbb{T}E,G)$ be splitting data, i.e. $G$ is the graph of an isomorphism $\alpha$ between the $d$-torsion points of $\mathbb{T}E'$ and $\mathbb{T}E$. 
    Then there exists:
    \begin{itemize}
        \item either a curve $\Gamma$ or a family of curves $\Gamma_t$ for $t\geqslant 0$, each of genus $2$,
        \item either a pair of optimal covers $(\varphi': \Gamma \rightarrow \mathbb{T}E',\varphi: \Gamma \rightarrow \mathbb{T}E)$ or a family of pairs $(\varphi_t,\varphi'_t)$ for $t\geqslant 0$, each of degree $d$,
    \end{itemize}
    that induce a splitting
    \begin{align}
      \phi:\mathbb{T}E' \bigoplus \mathbb{T}E \rightarrow \Jac(\Gamma)
  \end{align}
  with  $\Jac_d(\mathbb{T}E')\cong \ker(\phi) \cong \Jac_d(\mathbb{T}E)$.

\end{theorem}
Theorem \ref{theorem_reconstructionsumupinintro} shows: Plan \ref{masterplaninintro} has two possible outcomes. 
These correspond to the two maximal combinatorial types of curves of genus 2, the theta and the dumbbell graph, and can be distinguished using Algorithm \ref{algorithm_PreimageofTormap} (see Subsection \ref{subsection_step2} and Section \ref{section_appendix} for its implementation).

This naturally translates into a question posed in terms of moduli spaces, asking for the intersection of the locus $\mathcal{Q}$ of split Jacobians with the boundary of $A_2^{tr}$ (see Subsection \ref{subsection_Schottky-typeProblem}). The results of which we use to obtain the following combinatorial characterization of $\mathbb{T}\mathcal{L}_d$.

\begin{theorem}(Theorem \ref{theorem_locusofcurveswithdsplitJac})\label{theorem_locusofcurveswithdsplitJacinintro}

    The locus of curves with $d$-split Jacobian, $\mathbb{T}\mathcal{L}_d$, decomposes into $\varphi(d)$ subsets
    \begin{align}
        \mathbb{T}\mathcal{L}_d=\bigcup_{k}L_k,
    \end{align}
      where $\varphi(d)$ denotes the Euler phi function. Each $L_k$ corresponds to a 2 dimensional fan obtained by subdividing the positive quadrant: For each maximal cone $\tilde{\sigma}\subset L_k$ there exists a linear map $\phi_{\tilde{\sigma}}:\tilde{\sigma} \rightarrow M^{tr}_2$ that maps the open cone to a 2 dimensional cone in the theta part of $ M^{tr}_2$. Restricting $\phi_{\tilde{\sigma}}$ to a ray $r$ that is a face of two cones yields a linear map $\phi_{\tilde{\sigma} |r}\times id :r \times \mathbb{R}_{\geqslant 0} \rightarrow M^{tr}_2$ whose image in the dumbbell part of $ M^{tr}_2$ is $\phi_{\tilde{\sigma}}(r) \times r_B$, where $r_B$ denotes the ray that corresponds to the bridge edge. Moreover, the linear maps are compatible in the sense that for two neighboring maximal cones $\tilde{\sigma}$ and $\tilde{\sigma}'$ the restriction of $\phi_{\tilde{\sigma}}$, respectively $\phi_{\tilde{\sigma}'}$, to a common face agree.
\end{theorem}
 \subsection{Plan of Paper}
 This paper consists of three parts. In part 1, Sections \ref{section_catoftavspaper2} to \ref{section_bridge}, we mainly cover background material on tavs, tropical curves and their Jacobians.
 
 We start in Section \ref{section_catoftavspaper2} with a description of the category of tavs and review some categorical constructions such as kernels, products and coproducts that will be needed later. The material we discuss is known, except for Subsections \ref{subsection_isogeniesandpolarizations} and \ref{subsection_adjoints}. Here we develop one of the main tools for part 2, a tropical analogue of Mumford's Criterion (Proposition \ref{proposition_tropicalmumford}) and define adjoints in $\mathbb{T}\mathcal{A}$, a notion we will find useful for concise notation. Section \ref{section_catofcurvespaper2} introduces the category of tropical covers, combining \cite{MR3375652} and \cite{MR2772537}, and includes a short recall on tropical optimal covers, which we introduced in \cite{arXiv:2410.13459}.

 In Section \ref{section_bridge} we describe a connection between $\mathbb{T}\mathcal{C}$ and $\mathbb{T}\mathcal{A}$, which on the level of objects is given by $\Gamma \mapsto \Jac(\Gamma)$. Considering the reverse direction is tricky. The transition to moduli spaces (Section \ref{subsection_tropicaltorelli}) offers a constructive approach since the objects parameterizing pptavs and curves are in some sense more concrete. Our exposition is based on \cite{MR2968636} and \cite{MR3752493}. We will be as as brief as possible, only aiming at a visualization of the objects involved. This concludes part 1.

 The reconstruction procedure, which is the main focus of the paper, is discussed in Section \ref{section_reconstruction} with a computational supplement in Section \ref{section_appendix}. Here we prove Theorem \ref{theorem_reconstructionsumupinintro}, which may be retraced computationally using Algorithm \ref{algorithm_PreimageofTormap}. 
 
 In part 3, consisting of Section \ref{section_modulispaceperspective}, we investigate part 2 from the perspective of moduli spaces, study a tropical analogue of the locus of curves with split Jacobians (Theorem \ref{theorem_locusofcurveswithdsplitJacinintro}) and a Schottky type problem in this context. 
 \paragraph{\emph{Acknowledgements}} Thanks are due to Hannah Markwig for countless conversations and Felix Röhrle for thoughtful comments. 
 \section{Category of tropical abelian varieties}\label{section_catoftavspaper2}
 \subsection{Preliminaries}
 A large part of the following work will take place in the category of tropical abelian varieties $\gls{TA}$. Here we give a straight-to-the-point description. For a more complete introduction see \cite{arXiv:2410.13459}, Sections 2 and 3, which is based on \cite{MR4382460}, Section 2.3 and \cite{röhrle2024tropicalngonalconstruction}, Section 4. 
 
Let $\gls{Sigma}$ be a real torus with integral structure, i.e. $\gls{Sigma}=(\Lambda,\Lambda',[\cdot,\cdot])$, where
\begin{itemize}
\item $\Lambda$ and $\Lambda'$ are finitely generated free abelian groups of the same rank,
    \item $[\cdot,\cdot]: \Lambda \times \Lambda' \rightarrow \mathbb{R}$ is a non-degenerate pairing,
\end{itemize}
with topological realization $\Hom(\Lambda,\mathbb{R})/\Lambda'$, where $\Lambda'\subset \Hom(\Lambda,\mathbb{R})$ via 
$[\cdot,\cdot]$, and $\widecheck{\Sigma}:=(\Lambda',\Lambda,[\cdot,\cdot]^t)$ its dual (see \cite{arXiv:2410.13459}, Sections 2.1).
 \begin{definition}\label{definition_tav}
\emph{Objects:} A \emph{tropical abelian variety (tav)} $\gls{Sigma}=(\Lambda,\Lambda',[\cdot,\cdot])$ is a real torus with integral structure together with a \emph{polarization}, i.e. a group homomorphism $\gls{zeta}: \Lambda' \rightarrow \Lambda$  such that the bilinear form $[\zeta(\cdot),\cdot]: \Lambda' \times \Lambda' \rightarrow \mathbb{R}$ is symmetric and positive definite. The \emph{dimension} of $\Sigma$ is the $\mathbb{R}$-vector space dimension of $\Hom(\Lambda,\mathbb{R})$ and equal to $\rk(\Lambda)$ (equivalently  equal to $\rk(\Lambda')$). Moreover, $\zeta$ is called a \emph{principal polarization} and $ \Sigma$ a \emph{principally polarized tropical abelian variety (pptav)}, whenever $\zeta$ is a bijection.

\emph{Morphisms:} For $i=1,2$ let $\Sigma_i:=(\Lambda_i,\Lambda_i',[\cdot,\cdot]_i)$ be tavs. A \emph{morphism of tavs} is a pair of group homomorphisms $f:=(f^\#: \Lambda_2 \rightarrow \Lambda_1,f_\#: \Lambda'_1 \rightarrow \Lambda'_2)$  such that
\begin{diagram}\label{diagram_equationcompatibilitycondition}
    [f^\#(\lambda_2),\lambda'_1]_1=[\lambda_2,f_\#(\lambda'_1)]_2
\end{diagram}
is satisfied for all $\lambda'_1\in \Lambda'_1$ and $\lambda_2\in \Lambda_2$. 

\emph{Dual notions:} The dualization functor $\widecheck{\cdot}: \mathbb{T}\mathcal{A} \rightarrow \mathbb{T}\mathcal{A}$ sends a tav $\Sigma$ to $\widecheck{\Sigma}$ equipped with the so-called \emph{dual polarization} $\widecheck{\zeta}$ (see \cite{arXiv:2410.13459}, Section 3.1, previously defined by Röhrle and Zakharov, see first version of \cite{röhrle2024tropicalngonalconstruction}) and a morphism $f$, accordingly, to the \emph{dual morphism}, $\gls{widecheckf}$ obtained from the pair $(f_\#,f^\#)$.

Since properties of the induced map of quotients are encoded as properties of the pair $(f^\#,f_\#)$ (see \cite{röhrle2024tropicalngonalconstruction}, Definition 4.8), we call $f$
\begin{itemize}
    \item \emph{surjective}, if $f^\#$ is injective.
    \item \emph{finite}, if $f_\#$ is injective (equivalently if $[\Lambda_1 : f^\#(\Lambda_2)] < \infty$).
    \item \emph{injective},  is $f$ is finite and $f_\#(\Lambda'_1)$ is saturated in $\Lambda'_2$.
    \item an \emph{isogeny}, if it is surjective and finite.
\end{itemize}
\end{definition}
The category of tavs is particularly nice, it is abelian (see \cite{MR4261102}), finitely complete and finitely cocomplete (\cite{arXiv:2410.13459}, Lemma 11). This means that we have notions of (co)kernels and (co)products (and the like) and they behave as we would expect from the category of abelian groups. Knowing this will be sufficient for our purposes. Therefore, we consider explicit constructions only to a limited extend and refer the reader to \cite{arXiv:2410.13459} for more details.
 \begin{definition}(\cite{arXiv:2410.13459}, Definition 9) \label{definition_productandcoproducts}
Given tavs $\Sigma_1 $ and $ \Sigma_2$, we define their \emph{product} $\Sigma_1 \bigotimes \Sigma_2$ as follows: 
\begin{itemize}
    \item The underlying real torus with integral structure is $(\Lambda_1 \times \Lambda_2,\Lambda'_1 \times \Lambda'_2,[\cdot,\cdot]_1 + [\cdot,\cdot]_2)$, where $\Lambda_1 \times \Lambda_2$, respectively $\Lambda'_1 \times \Lambda'_2$, denotes the direct product of groups and $[\cdot,\cdot]_1+[\cdot,\cdot]_2$ is given by $((\lambda_1,\lambda_2),(\lambda'_1,\lambda'_2)) \mapsto [\lambda_1,\lambda'_1]_1 + [\lambda_2,\lambda'_2]_2$.
    \item  The polarization  $\zeta_1 \times \zeta_2$ is defined component-wise.
    \item The object  $\Sigma_1 \bigotimes \Sigma_2$ is equipped with a pair of morphisms, $\pi_1:\Sigma_1 \bigotimes \Sigma_2 \rightarrow \Sigma_1 $ and $\pi_2:\Sigma_1 \bigotimes \Sigma_2 \rightarrow \Sigma_2 $, induced by the natural projection and inclusion maps between the lattices. 
\end{itemize}
The resulting object is a tav, that is: 
\begin{itemize}
    \item The pairing $[\cdot,\cdot]_1+[\cdot,\cdot]_2$ is non-degenerate.
    \item The group homomorphism $\zeta_1 \times \zeta_2$ satisfies the condition described in Definition \ref{definition_tav}.
    \item The projection maps $\pi_1$ and $\pi_2$ are morphism of tori.
\end{itemize}
 The \emph{coproduct} of $\Sigma_1 $ and $ \Sigma_2$ is defined analogously and will be denoted by $\Sigma_1 \bigoplus \Sigma_2$.
   
\end{definition}
\begin{lemma}(\cite{arXiv:2410.13459}, Lemma 26)\label{lemma_quotientoftavbyfinitesubgrouppaper2}
Let $\Sigma$ be a tav and suppose $G \subset \Sigma$ is a finite subgroup. Then there exists a tav $\Sigma_G$ and a free isogeny $ \Sigma \overset{q}{\twoheadrightarrow} \Sigma_G$ whose kernel is $G$.
\end{lemma}
We also include the following Lemma proved in the first paper.
\begin{lemma}\label{lemma_dualizingexactsequences}[\cite{arXiv:2410.13459}, Lemma 23]
   Let   
    \begin{align}
        0 \rightarrow \Sigma_1 \xrightarrow{f} \Sigma_2 \xrightarrow{g} \Sigma_3 \rightarrow 0.
    \end{align}
    be a short exact sequence of tropical abelian varieties.
    Then the dual sequence
    \begin{align}
        0 \rightarrow \widecheck{\Sigma}_3 \xrightarrow{\widecheck{g}} \widecheck{\Sigma}_2 \xrightarrow{\widecheck{f}} \widecheck{\Sigma}_1 \rightarrow 0
    \end{align}
    is exact. In other words 
    \begin{align}
        \widecheck{\cdot}:\mathbb{T}\mathcal{A} \rightarrow \mathbb{T}\mathcal{A}
    \end{align}
    is an exact functor.
\end{lemma}
 \subsection{Isogenies and Polarizations}\label{subsection_isogeniesandpolarizations}
 Isogenies form a distinguished class of morphisms of tavs that are in some sense closely related to polarizations and interact well with them.
 \subsubsection{Polarizations induce isogenies}
 A polarization $\zeta$ defines an isogeny $f_\zeta:=(\zeta,\zeta)$ between $\Sigma$ and its dual $ \widecheck{\Sigma}$.
 We can characterize its kernel using the so-called \emph{type} of $\zeta$, which is given by the invariant factors $(\alpha_1,...,\alpha_n)$ (where $n:=\rk(\Lambda)$) of its Smith normal form. We have:
\begin{align}
    \ker(\gls{fzeta})\cong \mathbb{Z}/\alpha_1\mathbb{Z} \times ... \times \mathbb{Z}/\alpha_n\mathbb{Z}
\end{align}
and thus a "canonical" identification of $\Sigma$ and $\widecheck{\Sigma}$, whenever $f_\zeta$ is an isomorphism, or equivalently whenever $\zeta$ has type $(1,...,1)$ (i.e. $\zeta$ is a pp).
\subsubsection{Bidirectional transport of polarizations via isogenies}
Recall that isogenies (\cite{arXiv:2410.13459}, Section 3.1) allow us to transfer polarizations from the target variety to the domain and vice versa. In other words, any integral torus that is isogenous to a tav can be turned into tav itself. 
\begin{definition}\label{definition_induced/pf/pbpolarizationandpolarizedisogeny}
Let $\Sigma_2$ be a tav with polarization $\zeta_2$, $\Sigma_1$ a real torus with integral structure, and $\phi: \Sigma_1 \rightarrow \Sigma_2$ an isogeny. Then $\phi^* \zeta_2:=\phi^\#\circ \zeta_2 \circ \phi_\#$ is a polarization on $\Sigma_1 $ and called the \emph{induced polarization} or alternatively \emph{the pull-back} of $\zeta_2$ by $\phi$. Conversely, suppose $\Sigma_1$ carries a polarization $\zeta_1$. We can define the \emph{push-forward} of $\zeta_1$ by $\phi$ as $\phi_* \zeta_1:=\widecheck{\zeta}$, where $\zeta:= \widecheck{\phi}^* \widecheck{\zeta_1}$. We say that an isogeny $\phi: \Sigma_1 \rightarrow \Sigma_2$ is \emph{polarized} with respect to polarizations $\zeta_1$ on $\Sigma_1$ and $\zeta_2$ on $\Sigma_2$, if $\zeta_1$ is the polarization induced by $\phi$ and $\zeta_2$, in other words, if the diagram 
 \begin{center}
  
\begin{tikzcd}
\Sigma_1 \arrow[r, "\phi"] \arrow[d, "{f_{\zeta_1}}"]
& \Sigma_2 \arrow[d, "{f_{\zeta_2}}"] \\
\widecheck{\Sigma}_1
& \arrow[l, "\widecheck{\phi}"] \widecheck{\Sigma}_2
\end{tikzcd}
 \end{center}
commutes.
    
\end{definition}
\subsubsection{Polarizing isogenies}
Note that pulling back polarizations always turns $\phi$ into a polarized isogeny. The same is not true for the push-forward. Characterizing under which circumstances a polarization on the source is induced by a polarization on the target becomes relevant for Subsection \ref{subsection_step1} and is content of the subsequent Lemma. 
\begin{lemma}\label{lemma_inducingpolarization}
    Let $\phi: \Sigma_1 \rightarrow \Sigma_2$ be an isogeny and $\zeta_1: \Lambda'_1 \rightarrow \Lambda_1$ a polarization on $\Sigma_1$. We denote by $(\gamma_1,...,\gamma_n)$ the invariant factors of $i: \im(\phi_\#) \rightarrow \Lambda'_2$ and by $\widehat{\alpha_1}$ the first invariant factor of $\zeta_1: \Lambda'_1 \rightarrow \im(\phi^\#)$ (meaning the smallest with respect to divisibility). Then, there exists a polarization $\zeta_2$ on $\Sigma_2$ such that $\phi^*\zeta_2=\zeta_1$, if 
    \begin{enumerate}
        \item $\im(\phi^\#) \supset \im(\zeta_1)$.
        \item $\Det(i):=\prod^n_{j=1}\gamma_j$ is a divisor of $\widehat{\alpha_1}$.
    \end{enumerate}
\end{lemma}
While condition (1) in Lemma \ref{lemma_inducingpolarization} is necessary for the existence of the polarization $\zeta_2$, the divisibility condition (condition (2)) is only sufficient: It has been formulated so as to provide the reader with easily verifiable criterion. The sacrifice of necessity is deliberate.
\begin{proof}
  Our goal is to complete 
  \begin{diagram}\label{inproof_lemmainducingpol}
 \Lambda_1  
& \arrow[l, "\phi^\# "] \Lambda_2  \\
\Lambda'_1 \arrow[r,"\phi_\#"] \arrow[u, "\zeta_1 "]  & \arrow[u, dashed, "\zeta_2 "] \Lambda'_2\\
  \end{diagram} to a commutative diagram by connecting the lattices $\Lambda'_2$ and $\Lambda_2$ by means of a polarization $\zeta_2$ on $\Sigma_2$. We proceed in two steps:
\begin{enumerate}[(i)]
    \item Invert \cref{inproof_lemmainducingpol} "as much as possible" and define a polarization $\tilde{\zeta_2}$ on a sublattice of $\Lambda'_2$.
    \item Extend $\tilde{\zeta_2}$ to $\Lambda'_2$ in an appropriate way.
\end{enumerate}
  Consider 
  \begin{diagram}
\im(\phi^\#) \arrow[rr, bend left, "g_2"] &  \Lambda_1 
& \arrow[l, "\phi^\#"] \Lambda_2  \\
& \Lambda'_1 \arrow[r, "\phi_\#"] \arrow[u, "\zeta_1"] \arrow[ul, bend left, "\zeta_1"]
& \Lambda'_2 \\
& &\arrow[ul, bend left, "g_1"] \arrow[u, hook, "i"] \im(\phi_\#)
  \end{diagram}
where $g_1$ and $g_2$ are the inverses of the isomorphisms obtained from $\phi^\#$ and $\phi_\#$ by restricting their respective codomains and define $\tilde{\zeta_2}: \im(\phi_\#) \rightarrow \Lambda_2$ as the composition $g_2\circ \zeta_1 \circ g_1$. By construction $\tilde{\zeta_2}$ satisfies $\phi^\#\circ \tilde{\zeta_2} \circ \phi_\#=\zeta_1$ and for all $\lambda \in \im(\phi_\#)$ (i.e. $\lambda=\phi_\#(\lambda'_1)$ for some $\lambda'_1 \in \Lambda'_1$) we have:
\begin{align}
    [\tilde{\zeta_2}(\lambda),\lambda]_2&=[\tilde{\zeta_2}(\phi_\#(\lambda'_1)),\phi_\#(\lambda'_1)]_2=[\phi^\#\circ \tilde{\zeta_2}\circ \phi_\#(\lambda'_1)),\lambda'_1]_1\\
    &=[\phi^\#\circ g_2\circ \zeta_1 \circ g_1\circ \phi_\#(\lambda'_1)),\lambda'_1]_1=[\zeta_1 (\lambda'_1),\lambda'_1]_1.
\end{align} Hence, $[\tilde{\zeta_2}(\cdot),\cdot]_2: \im(\phi_\#)_\mathbb{R} \times \im(\phi_\#)_\mathbb{R} \rightarrow \mathbb{R}$ is symmetric and positive definite. This finishes step (i). For step (ii), let $n=\dim(\Sigma_1)=\dim(\Sigma_2)$ and recall that both $\im(\phi_\#)$ and $\im(\phi^\#)$ are full-dimensional sublattices since $\phi$ is an isogeny. 
We extend $\tilde{\zeta_2}$ injectively along $i$ by looking at the diagram in Figure \ref{figure_inprooftropicalmumford} (a) in coordinates (see Figure \ref{figure_inprooftropicalmumford} (b)).
\begin{figure}[H]
    \begin{tikzcd}
(a) & \Lambda'_2 \arrow[r, dashed]
&  \Lambda_2  & (b) & \mathbb{Z}^n \arrow[r, dashed]
&  \mathbb{Z}^n  \\
& \im(\phi_\#) \arrow[ur,"\tilde{\zeta_2}", swap] \arrow[u, hook, "i "]  & & & \mathbb{Z}^n \arrow[ur, "M(\tilde{\zeta_2})", swap] \arrow[u, hook, "M(i)"]  & \\
\end{tikzcd}
    \caption{Transition to coordinates in the proof of Lemma \ref{lemma_inducingpolarization}.}
    \label{figure_inprooftropicalmumford}
\end{figure}

On $\im(\phi_\#)$ and $\Lambda_2$ these are chosen such that the transformation matrix of $\tilde{\zeta_2}$, $M(\tilde{\zeta_2})$, is in Smith normal form. We impose no restrictions on the choice for $\Lambda'_2$. Set 
\begin{align}
    M(\zeta_2):= diag(\widehat{\alpha}_1,...,\widehat{\alpha}_n)\cdot M(i)^{-1},
\end{align}
where $(\widehat{\alpha}_1,...,\widehat{\alpha}_n)$ are the invariant factors of $\tilde{\zeta_2}$ and $M(i)^{-1}$ is the inverse of $M(i)$ over $\mathbb{Q}$. Since $M(i)^{-1}$ is the product of its adjoint and the scalar $\frac{1}{\Det(i)}$ and $\Det(i)$ is a divisor of $\widehat{\alpha_1}$, $M(\zeta_2)$ is an integer matrix. By forgetting coordinates, we obtain an extension, $\zeta_2$, of $\tilde{\zeta_2}$ along $i$ as desired. Note that $\zeta_2$ is automatically injective as $\det(M(\zeta_2))$ is non zero.\\
We conclude the proof by verifying that $\zeta_2$ induces a scalar product $\langle \cdot , \cdot \rangle:=[\zeta_2(\cdot),\cdot ]_2$ on $\Lambda'_{2 \thinspace \mathbb{R}}$, equivalently that the matrix associated to $\langle \cdot , \cdot \rangle$ is symmetric and positive definite. In order to provide us with coordinate representations, $C$ and $C'$, of $\langle \cdot , \cdot \rangle$ and its restriction to $\im(\phi_\#)_\mathbb{R} \times \im(\phi_\#)_\mathbb{R}$, we fix lattice basis $S,S'$ and $S''$ of $\Lambda_2, \Lambda'_2$ and $\im(\phi_\#)$. 
Doing so, we observe that $C$ and $C'$ are related by the following equation: 
\begin{align}
    C'= _{S''}{M^t}_{S'}(i) \cdot C \cdot _{S''}{M}_{S'}(i).
\end{align}
Since $C'$ is symmetric, we have
\begin{align}
     _{S''}{M^t}_{S'}(i) \cdot C^t \cdot _{S''}{M}_{S'}(i)=C'^t=C'= _{S''}{M^t}_{S'}(i) \cdot C \cdot _{S''}{M}_{S'}(i)
\end{align}
and conclude that $C^t=C$ using that $_{S''}{M}_{S'}(i)$ is invertible over $\mathbb{R}$. Choosing $S'$ and $S''$ such that $_{S''}{M}_{S'}(i)$ is in Smith normal form with invariant factors $(\beta_1,...,\beta_n)$, we can compute the leading principal minors $\det(C_k)$ of $C$ in terms of the ones of $C'$:
\begin{align}
    \det(C'_k)=\prod^k_{i=1} \beta^2_i \cdot \det(C_k), k=1,...,n.
\end{align} By Sylvester's criterion all of these must be positive, since $C'$ is symmetric and positive definite. Hence, $C$ is as well.
\end{proof}
\begin{remark}
In the setting of Lemma \ref{lemma_inducingpolarization}, note that whenever $\zeta_2$ exists it is unique: The condition that Diagram \ref{inproof_lemmainducingpol} commutes requires any such $\zeta_2$ to be an extension of $\tilde{\zeta_2}$, which is fixed (by $\zeta_1$ and $\phi$) on a sublattice of $\Lambda'_2$ of full rank. In this case we call $\zeta_2$ the \emph{inducing} polarization.\\
The divisibility criterion in Lemma \ref{lemma_inducingpolarization} involves the invariant factors of the polarization $\zeta_1$ 
on $\Sigma_1$, however with restricted target. These relate to the type $(\alpha_1,...,\alpha_n)$ of $\zeta_1$ in the following way (\cite{MR1341069}, Section 4):
\begin{align}
 \lcm( \thinspace \widehat{\alpha}_{i+1} \widehat{\gamma}_{k-i}: \thinspace 0 \leqslant i \leqslant k-1) \thinspace | \thinspace \alpha_{k}\thinspace| \thinspace \gcd(\thinspace \widehat{\alpha}_{k-1+i} \widehat{\gamma}_{n-i+1}: \thinspace 1 \leqslant i \leqslant n-k+1),
\end{align}
where $(\widehat{\gamma}_1,...,\widehat{\gamma}_n)$ denote the invariant factors of the inclusion $\im(\phi^\#) \hookrightarrow \Lambda_1$.
If $n=2$ this yields 
\begin{align}
    \widehat{\alpha}_{1} \widehat{\gamma}_{1} \thinspace | \thinspace \alpha_{1}\thinspace| \thinspace \gcd(\widehat{\alpha}_{1} \widehat{\gamma}_{2}, \widehat{\alpha}_{2} \widehat{\gamma}_{1}) \text{ and } \lcm(\widehat{\alpha}_{1} \widehat{\gamma}_{2},\widehat{\alpha}_{2} \widehat{\gamma}_{1}) \thinspace | \thinspace \alpha_{2}\thinspace| \thinspace \widehat{\alpha}_{2} \widehat{\gamma}_{2}.
\end{align}
\end{remark}
Clearly, the proof of Lemma \ref{lemma_inducingpolarization} has an algorithmic flavor, which we now concretize.
\begin{algorithm}\label{algorithm_constructionofinducingpolarization}
\ \\
Input: An isogeny $\phi$ and a polarization $\zeta_1$.\\
Output: A matrix $M$ that decides whether $\zeta_2$ exists.
\begin{enumerate}
    \item[I.] For $\mathbb{Z}$-basis $T'$ and $S$ of $\Lambda'_1$ and $\Lambda_2$, let $S''$ and $T''$ denote the basis on $\im(\phi_{\#})$ and $\im(\phi^{\#})$ induced by the injective maps $\phi_{\#}$ and $\phi^{\#}$. Set \\
    $A:=_{T'}{M}_{T''}(\zeta_1)$.
    \item[II.] Choose a $\mathbb{Z}$-basis $S'$ of $\Lambda'_2$ and let $B:=_{S''}{M}_{S'}(i)$. 
    \item[III.] Compute $M:=A\cdot B^{-1}$.
    \end{enumerate}
    \end{algorithm}
Algorithm \ref{algorithm_constructionofinducingpolarization} provides us with an if-and-only-if-criterion for the existence of $\zeta_2$, which is accessible for algorithmic verification and turns into a step-by-step construction should $\zeta_2$ exist. The classical counterpart of Proposition \ref{proposition_tropicalmumford} (see \cite{zbMATH05564780}, Proposition 16.8 or \cite{MR0282985}, p.231) is of a different flavor: It places certain requirements on a pairing $e^\lambda$ which is associated to the polarization $\lambda$ and a matching Weil pairing. We consider this difference to be one of the main reasons why Plan \ref{masterplaninintro} can be carried out explicitly in the tropical setting.

\begin{proposition}\label{proposition_tropicalmumford}
    Let $\phi: \Sigma_1 \rightarrow \Sigma_2$ be an isogeny and $\zeta_1: \Lambda'_1 \rightarrow \Lambda_1$ a polarization on $\Sigma_1$. There exists a polarization $\zeta_2$ on $\Sigma_2$ such that $\phi^*\zeta_2=\zeta_1$ if and only if the output $M$ of Algorithm \ref{algorithm_constructionofinducingpolarization} is an integer matrix. In this case, $M$ is a coordinate representation of $\zeta_2$.
\end{proposition}
 \subsection{Adjoints}\label{subsection_adjoints}
 We will find the following to be a useful addition to our toolbox.
\begin{definition}\label{definition_adjoint}
    Let $f: \Sigma_1 \rightarrow \Sigma_2$ and $\tilde{f} : \Sigma_2 \rightarrow \Sigma_1$ be two morphism of integral tori. We say that  $\tilde{f}$ is \emph{adjoint} to $f$ (alternatively that $(f,\tilde{f})$ is an \emph{adjoint pair}), if their exists principal polarizations $\zeta_1$ on $\Sigma_1$ and  $\zeta_2$ on $\Sigma_2$ such that the following diagram commutes:
     \begin{diagram}
 \Sigma_1  \arrow[d, "f_{\zeta_1} "]
&  \arrow[l, swap, "\tilde{f}"] \Sigma_2 \arrow[d, "f_{\zeta_2}."]  \\
\widecheck{\Sigma}_1  & \arrow[l,"\widecheck{f}"] \Check{\Sigma}_2  \\
  \end{diagram}
  
\end{definition}
Note that though the notion of adjoint is not unique in the strict sense of the word, it is "essentially" so (meaning unique up to conjugation by isomorphisms).
\begin{lemma}\label{lemma_propertiesofadjoints}
    (Properties of adjoint morphisms)\\
    Let $(f,\tilde{f})$ be an adjoint pair. For $i=1,2$ denote by $\langle \cdot, \cdot \rangle_i$ the scalar product on $\Lambda'_{i \thinspace \mathbb{R}}$ induced by $\zeta_i$. We have: 
    \begin{enumerate}
        \item $\langle \lambda'_1, \tilde{f}_\#(\lambda'_2) \rangle_1=\langle f_\#(\lambda'_1), \lambda'_2 \rangle_2$ for all $\lambda'_1\in \Lambda'_1$ and $\lambda'_2\in \Lambda'_2$.
        \item $\Ker(\tilde{f})_0\cong \widecheck{\Coker}(f)$.
    \end{enumerate}
    Moreover, if $f$ is an isogeny, then so is $\tilde{f}$ and the composition $\tilde{f} \circ f$ is given by $(f^* \zeta_2 \circ \zeta^{-1}_1, \zeta^{-1}_1  \circ f^* \zeta_2)$.
\end{lemma}
\begin{proof}
    Points $(1)$,$(2)$ and the formula for $\tilde{f} \circ f$ are easily derived by manipulating definitions. Since the dualization functor, $\widecheck{\cdot}$,  takes isogenies to isogenies, we see that $f$ is an isogeny if and only if $\tilde{f}$ is one.
\end{proof}
 \section{Category of tropical Curves}\label{section_catofcurvespaper2}
 \subsection{Preliminaries}
 The second category we work in is the \emph{category of tropical curves}, $\gls{TC}$. Our exposition is based on \cite{MR3375652} and \cite{MR2772537}, relying on the foundational work of Mikhalkin in \cite{MR2275625} and of Mikhalkin and Zharkov in \cite{MR2457739}.
 \begin{definition}\label{definition_tropicalcurve}
 \emph{Objects}: A metric graph $(G,l)$ is a finite graph $G$ (with no legs/ends) and a function $l: E(G) \rightarrow \mathbb{R}_{>0}$. The geometric realization $\Gamma$ of $(G,l)$ is called a \emph{tropical curve} of genus $g(\Gamma):=b_1(G)$ with \emph{model} $(G,l)$ and \emph{combinatorial type} $G$, where $b_1(G)$ denotes the first Betti number of $G$. \\

   \emph{Morphisms:} A continuous and surjective map $\varphi: \Gamma \rightarrow \tilde{\Gamma} $ between two tropical curves $\Gamma$ and $\tilde{\Gamma}$ is called a \emph{tropical cover} (see e.g. \cite{MR2525845} or \cite{MR3278571}, Section 2.1) if there exists models $(G,l)$ and $(\tilde{G},\tilde{l})$ such that
\begin{itemize}
    \item $\varphi(V(G))\subset V(\tilde{G})$ and $\varphi^{-1}(E(\tilde{G}))\subset E(G)$.
    \item $\varphi$ is locally integer affine linear: On each edge $e\in E(G)$, $\varphi$ restricts to an affine function with integer slope $d_e(\varphi)$ (possibly $0$), called the \emph{weight} or \emph{expansion factor} of $\varphi$ at $e$.
    \item $\varphi$ is harmonic/balanced at every $P\in \Gamma$: For any $\tilde{v}\in T_{\varphi(P)} \tilde{\Gamma}$
    \begin{align}
        d_P(\varphi):=\sum_{v\in T_P\Gamma, v \mapsto \tilde{v}} d_{v}(\varphi)
    \end{align}
    is independent of $\tilde{v}$, where $T_P\Gamma$, respectively $T_{\varphi(P)} \tilde{\Gamma}$, is the set of tangent directions emanating from $P$, respectively $\varphi(P)$, and $d_{v}(\varphi)$ is the directional derivative of $\varphi$ in the direction of $v$ (i.e. $d_{v}(\varphi):=d_e(\varphi)$ for the edge $e$ in direction of $v$).
    \end{itemize}
    We say that $\varphi$ is \emph{finite}, if $d_e(\varphi)>0$ for all edges $e$, and \emph{non-finite} else. In any case, $deg(\varphi):=\sum_{P\in \Gamma, P \mapsto \tilde{P}} d_{P}(\varphi)$, where $\tilde{P}\in \tilde{\Gamma}$ is an arbitrary point, provides a well-defined notion of \emph{degree}, a tropical analogue of the corresponding notion in algebraic geometry. 
\end{definition}
 
We will work with tropical curves through choice of a model and by abuse of notation identify $\Gamma$ with $(G,l)$.
\subsection{Setup}
In Sections \ref{section_reconstruction} and \ref{section_modulispaceperspective} we only work in the genus 1 and 2 part of $\mathbb{T}\mathcal{C}$. We recall a class of morphism introduced in \cite{arXiv:2410.13459} that is relevant in this context:

\begin{definition}\label{definition_optimalmappaper2}
Let $\Gamma$ and $\mathbb{T}E$ be tropical curves of genus 2 and 1. We call a cover $\varphi: \Gamma \rightarrow \mathbb{T}E$ \emph{optimal} if it does not factor through a non-trivial cover, i.e. if there exists a curve $\mathbb{T}\tilde{E}$ and maps $\tilde{\varphi}: \Gamma \rightarrow \mathbb{T}\tilde{E}$, $\phi: \mathbb{T}\tilde{E} \rightarrow \mathbb{T}E $ such that 
\begin{center}
       \begin{tikzcd}
\Gamma \arrow[r,"\tilde{\varphi}"] \arrow[dr,"\varphi",swap]
& \mathbb{T}\tilde{E} \arrow[d,"\phi"]\\
& \mathbb{T}E
\end{tikzcd}
\end{center}
 
commutes, then $\phi$ is an isomorphism (i.e. $deg(\phi)=1$).
\end{definition}
 \begin{figure}[hbt!]
    \centering
 
\begin{tikzpicture}[x=0.75pt,y=0.75pt,yscale=-1,xscale=1]

\draw   (108.21,210.52) .. controls (108.23,202.42) and (140.77,195.86) .. (180.88,195.87) .. controls (221,195.88) and (253.5,202.47) .. (253.48,210.57) .. controls (253.46,218.68) and (220.92,225.24) .. (180.8,225.23) .. controls (140.69,225.21) and (108.18,218.63) .. (108.21,210.52) -- cycle ;
\draw    (214.94,123.82) -- (197.63,155.1) ;
\draw [shift={(196.66,156.85)}, rotate = 298.95] [color={rgb, 255:red, 0; green, 0; blue, 0 }  ][line width=0.75]    (10.93,-3.29) .. controls (6.95,-1.4) and (3.31,-0.3) .. (0,0) .. controls (3.31,0.3) and (6.95,1.4) .. (10.93,3.29)   ;
\draw  [draw opacity=0] (185.85,78.52) .. controls (171.5,74.83) and (162.82,70.04) .. (162.83,64.8) .. controls (162.84,53.27) and (205.02,43.97) .. (257.03,44.04) .. controls (309.04,44.1) and (351.2,53.5) .. (351.18,65.03) .. controls (351.18,69.97) and (343.43,74.5) .. (330.5,78.06) -- (257,64.92) -- cycle ; \draw   (185.85,78.52) .. controls (171.5,74.83) and (162.82,70.04) .. (162.83,64.8) .. controls (162.84,53.27) and (205.02,43.97) .. (257.03,44.04) .. controls (309.04,44.1) and (351.2,53.5) .. (351.18,65.03) .. controls (351.18,69.97) and (343.43,74.5) .. (330.5,78.06) ;  
\draw   (184.94,78.28) .. controls (184.96,70.31) and (217.89,63.86) .. (258.48,63.87) .. controls (299.07,63.89) and (331.95,70.36) .. (331.93,78.33) .. controls (331.91,86.3) and (298.98,92.75) .. (258.39,92.73) .. controls (217.8,92.72) and (184.92,86.25) .. (184.94,78.28) -- cycle ;
\draw   (271.75,209.31) .. controls (271.77,201.9) and (302.15,195.91) .. (339.59,195.92) .. controls (377.03,195.93) and (407.36,201.95) .. (407.33,209.35) .. controls (407.31,216.76) and (376.93,222.75) .. (339.49,222.74) .. controls (302.05,222.73) and (271.72,216.71) .. (271.75,209.31) -- cycle ;
\draw    (303.27,124.82) -- (318.56,153.09) ;
\draw [shift={(319.51,154.85)}, rotate = 241.59] [color={rgb, 255:red, 0; green, 0; blue, 0 }  ][line width=0.75]    (10.93,-3.29) .. controls (6.95,-1.4) and (3.31,-0.3) .. (0,0) .. controls (3.31,0.3) and (6.95,1.4) .. (10.93,3.29)   ;
\draw    (184.94,78.28) ;
\draw [shift={(184.94,78.28)}, rotate = 0] [color={rgb, 255:red, 0; green, 0; blue, 0 }  ][fill={rgb, 255:red, 0; green, 0; blue, 0 }  ][line width=0.75]      (0, 0) circle [x radius= 3.35, y radius= 3.35]   ;
\draw    (331.36,77.82) -- (331.93,78.33) ;
\draw [shift={(331.93,78.33)}, rotate = 41.65] [color={rgb, 255:red, 0; green, 0; blue, 0 }  ][fill={rgb, 255:red, 0; green, 0; blue, 0 }  ][line width=0.75]      (0, 0) circle [x radius= 3.35, y radius= 3.35]   ;
\draw    (147.93,222.92) ;
\draw [shift={(147.93,222.92)}, rotate = 0] [color={rgb, 255:red, 0; green, 0; blue, 0 }  ][fill={rgb, 255:red, 0; green, 0; blue, 0 }  ][line width=0.75]      (0, 0) circle [x radius= 3.35, y radius= 3.35]   ;
\draw    (221.03,222.92) ;
\draw [shift={(221.03,222.92)}, rotate = 0] [color={rgb, 255:red, 0; green, 0; blue, 0 }  ][fill={rgb, 255:red, 0; green, 0; blue, 0 }  ][line width=0.75]      (0, 0) circle [x radius= 3.35, y radius= 3.35]   ;
\draw    (303.27,220.92) ;
\draw [shift={(303.27,220.92)}, rotate = 0] [color={rgb, 255:red, 0; green, 0; blue, 0 }  ][fill={rgb, 255:red, 0; green, 0; blue, 0 }  ][line width=0.75]      (0, 0) circle [x radius= 3.35, y radius= 3.35]   ;
\draw    (366.22,220.92) ;
\draw [shift={(366.22,220.92)}, rotate = 0] [color={rgb, 255:red, 0; green, 0; blue, 0 }  ][fill={rgb, 255:red, 0; green, 0; blue, 0 }  ][line width=0.75]      (0, 0) circle [x radius= 3.35, y radius= 3.35]   ;

\draw (133.22,34.69) node [anchor=north west][inner sep=0.75pt]    {$\Gamma $};
\draw (176.57,121.14) node [anchor=north west][inner sep=0.75pt]    {$\tilde{\varphi }$};
\draw (163.1,37.56) node [anchor=north west][inner sep=0.75pt]    {$1$};
\draw (254.13,73.03) node [anchor=north west][inner sep=0.75pt]    {$1$};
\draw (254.13,46.99) node [anchor=north west][inner sep=0.75pt]    {$1$};
\draw (427.86,208.4) node [anchor=north west][inner sep=0.75pt]    {$\mathbb{T} E$};
\draw (336.84,127.09) node [anchor=north west][inner sep=0.75pt]    {$\varphi $};
\draw (125.36,181.62) node [anchor=north west][inner sep=0.75pt]    {$8$};
\draw (388.35,182.62) node [anchor=north west][inner sep=0.75pt]    {$2$};
\draw (66.63,201.1) node [anchor=north west][inner sep=0.75pt]    {$\mathbb{T}\tilde{E}$};
\draw (178.14,208.65) node [anchor=north west][inner sep=0.75pt]    {$4$};
\draw (332.47,201.64) node [anchor=north west][inner sep=0.75pt]    {$1$};

\end{tikzpicture}

    \caption{A curve of genus $2$ covering two elliptic curves. The numbers are edge lengths.}
    \label{figure_optimalandnotoptimalcover}
\end{figure}
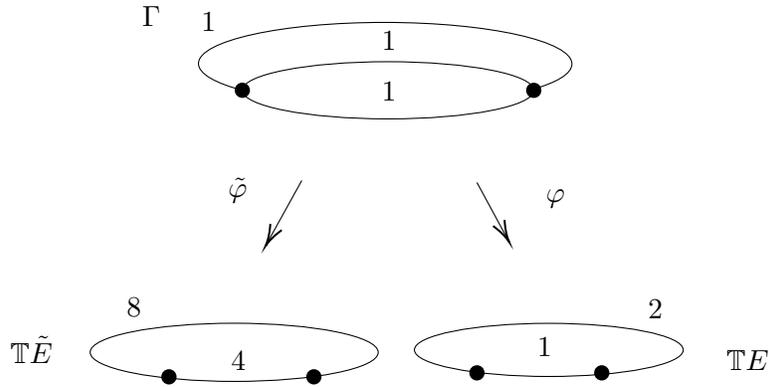
\begin{example}\label{example_optimalandnotoptimalpaper2}
 Figure \ref{figure_optimalandnotoptimalcover} shows a cover of degree $2$ on the right, which must be optimal. The cover on the left, however, is not as it factors, for example, through the first, giving rise to a cover of degree $4$.
\end{example}
 \section{Bridging \ensuremath{\mathbb{T}\mathcal{C}} and \ensuremath{\mathbb{T}\mathcal{A}}}\label{section_bridge}
 \subsection{Preliminaries}
 The categories $\mathbb{T}\mathcal{C}$ and $\mathbb{T}\mathcal{A}$ are connected:
 \begin{itemize}
     \item To each object $\Gamma$ we assign an object $\Jac(\Gamma)$ in $\mathbb{T}\mathcal{A}$,
     \item To each morphism $\varphi: \Gamma_1 \rightarrow \Gamma_2$ we associate a morphism $\varphi_*$ in $\mathbb{T}\mathcal{A}$,
 \end{itemize}
 where $\Jac(\Gamma)$ denotes the Jacobian of $\Gamma$ and $\varphi_*: \Jac(\Gamma_1) \rightarrow \Jac(\Gamma_2)$ the push-forward of $\varphi$ (\cite{arXiv:2410.13459}, Section 5). This connection, especially in the case of curves of genus 2 covering curves of genus 1, was discussed in detail in the first paper. 
 We briefly recall some constructions relevant in this context.
 \begin{construction}\label{construction_theJacobianpaper2}(see \cite{MR2772537} or \cite{MR4382460})
    Let $(G,l)$ be an oriented model of $\Gamma$ and $s,t: E(G) \rightarrow V(G)$ the source and target maps. Then $(G,l)$ comes with two lattices that are related by a non-degenerate pairing:
    \begin{itemize}
        \item \emph{The lattice of harmonic} $1$-\emph{forms}, $\Omega^1_G(\mathbb{Z})$: For each oriented edge $e$ we introduce a formal symbol $de$ called a \emph{basic} $1$-form on $G$ and set $\Omega^1_G(\mathbb{Z})$ to be 
        \begin{align}
            \{ \omega:=\sum_{e} \omega_e de: \omega_e \in \mathbb{Z}, \thinspace \sum_{e:t(e)=V} \omega_e = \sum_{e:s(e)=V} \omega_e \thinspace \forall \thinspace V\in V(G) \}.
        \end{align} 
        It is the subgroup of the free group over $\{ de: e\in E(G)\} $ consisting of harmonic $1-$forms on $G$.
        \item \emph{The lattice of integral} $1$-\emph{cycles}, $\H1(G,\mathbb{Z})$: It is the first simplicial homology group of $G$ given by $\ker(\partial)$, where 
        \begin{align}
            \partial: \C1(G,\mathbb{Z})\rightarrow \Co (G,\mathbb{Z}), e\mapsto t(e)-s(e)
        \end{align} is the boundary operator.
        \item \emph{The integration pairing}, $\int_\cdot \cdot$: We can integrate a basic $1$-form
    \begin{align}
        \int_e de':=\begin{cases}
            l(e), \text{ if } e=e'\\
            0, \text{ else}
        \end{cases}
    \end{align}
    and extend linearly to obtain a perfect pairing
    \begin{align}
       \int_\cdot \cdot: \Omega^1_G(\mathbb{Z}) \times \H1(G,\mathbb{Z}) \rightarrow \mathbb{R}, (\omega, c) \mapsto \int_c \omega.
    \end{align}
    \end{itemize}
    These building blocks are independent of the choice of model (see also \cite{MR2772537}). This means that lattices that arise from different models (that have compatible orientations) are related by isomorphisms, that leave the integration pairing invariant. We will write $\Omega^1_\Gamma(\mathbb{Z})$ and $\H1(\Gamma,\mathbb{Z})$, instead, and complete Construction \ref{construction_theJacobianpaper2} by assigning a pptav to $\Gamma$.
\begin{definition}\label{definition_JacobianandAbelJacobimappaper2}
    The \emph{Jacobian of} $\Gamma$ is the pptav built from $(\Omega^1_\Gamma(\mathbb{Z}),\H1(\Gamma,\mathbb{Z}),\int_{\cdot} \cdot)$ with principal polarization $\gls{zetaGamma}: \H1(\Gamma,\mathbb{Z}) \rightarrow \Omega^1_\Gamma(\mathbb{Z}), \sum a_e e \mapsto \sum a_e de $. It is related to $\Gamma$ by the \emph{tropical Abel-Jacobi map}:
\begin{align}
    \gls{PhiP0}: \Gamma \rightarrow \Jac(\Gamma), P \mapsto \int_{\gamma_{P}}\cdot,
\end{align}
where $P_0\in V(\Gamma)$ is a fixed vertex and $\gamma_P\in \C1(\Gamma,\mathbb{Z})$ is any path connecting $P_0$ to $P$ in $\Gamma$.
\end{definition}
\end{construction}
 As in classical algebraic geometry, we can identify $\Jac(\Gamma)$ with $\Pic^0(\Gamma):=\Div^0(\Gamma)/\Prin(\Gamma)$, where $\Div^0(\Gamma)$ is the group of divisors of degree $0$ and $\Prin(\Gamma)$, the subgroup of principal divisors (\cite{MR2457739}, Section 4.2). Then $\varphi_*$ is the morphism induced by the push-forward of divisors under this identification, whose dual $\varphi^*$ is induced by the pull-back (\cite{arXiv:2410.13459}, Lemma 35).
 \subsection{Setup}
 Specializing to the case of curves of genus 2 covering curves of genus 1, we are interested in a particular class of objects in $\mathbb{T}\mathcal{A}$, so-called \emph{tropical split Jacobians}.
 \begin{definition}\label{definition_splitJacobianspaper2}(\cite{arXiv:2410.13459}, Definition 58)
    Let $\Gamma$ be a tropical curve of genus $2$. We say that $\Jac(\Gamma)$ \emph{splits}, if $\Jac(\Gamma)$ is isogeneous to the coproduct of two elliptic curve $\mathbb{T}E\oplus \mathbb{T}E'$. In this case we call $\mathbb{T}E'$ a \emph{complement} of $\mathbb{T}E$ and vice versa.
\end{definition}
 \begin{convention}\label{convention_thetaanddumbbell}
     We fix a labelling for the combinatorial types we are going to work with (see Figure \ref{figure_conventionforthetaanddumbbell}) 
     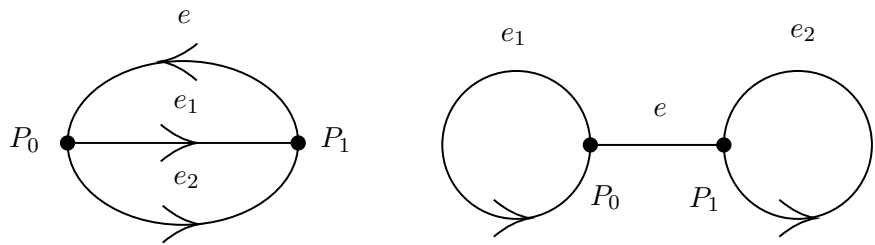
\begin{figure}[H]
         \centering
         \tikzset{every picture/.style={line width=0.75pt}} 

\begin{tikzpicture}[x=0.75pt,y=0.75pt,yscale=-1,xscale=1]

\draw   (62,111.67) .. controls (62,88.84) and (88.04,70.33) .. (120.17,70.33) .. controls (152.29,70.33) and (178.33,88.84) .. (178.33,111.67) .. controls (178.33,134.49) and (152.29,153) .. (120.17,153) .. controls (88.04,153) and (62,134.49) .. (62,111.67) -- cycle ;
\draw   (251,112.67) .. controls (251,92.05) and (267.71,75.33) .. (288.33,75.33) .. controls (308.95,75.33) and (325.67,92.05) .. (325.67,112.67) .. controls (325.67,133.29) and (308.95,150) .. (288.33,150) .. controls (267.71,150) and (251,133.29) .. (251,112.67) -- cycle ;
\draw   (393,112.67) .. controls (393,92.05) and (409.71,75.33) .. (430.33,75.33) .. controls (450.95,75.33) and (467.67,92.05) .. (467.67,112.67) .. controls (467.67,133.29) and (450.95,150) .. (430.33,150) .. controls (409.71,150) and (393,133.29) .. (393,112.67) -- cycle ;
\draw    (62,111.67) -- (178.33,111.67) ;
\draw [shift={(178.33,111.67)}, rotate = 0] [color={rgb, 255:red, 0; green, 0; blue, 0 }  ][fill={rgb, 255:red, 0; green, 0; blue, 0 }  ][line width=0.75]      (0, 0) circle [x radius= 3.35, y radius= 3.35]   ;
\draw [shift={(62,111.67)}, rotate = 0] [color={rgb, 255:red, 0; green, 0; blue, 0 }  ][fill={rgb, 255:red, 0; green, 0; blue, 0 }  ][line width=0.75]      (0, 0) circle [x radius= 3.35, y radius= 3.35]   ;
\draw    (325.67,112.67) -- (393,112.67) ;
\draw [shift={(393,112.67)}, rotate = 0] [color={rgb, 255:red, 0; green, 0; blue, 0 }  ][fill={rgb, 255:red, 0; green, 0; blue, 0 }  ][line width=0.75]      (0, 0) circle [x radius= 3.35, y radius= 3.35]   ;
\draw [shift={(325.67,112.67)}, rotate = 0] [color={rgb, 255:red, 0; green, 0; blue, 0 }  ][fill={rgb, 255:red, 0; green, 0; blue, 0 }  ][line width=0.75]      (0, 0) circle [x radius= 3.35, y radius= 3.35]   ;
\draw   (110,143.33) .. controls (116.11,148.52) and (122.22,151.63) .. (128.33,152.67) .. controls (122.22,153.7) and (116.11,156.81) .. (110,162) ;
\draw   (109,102.33) .. controls (115.11,107.52) and (121.22,110.63) .. (127.33,111.67) .. controls (121.22,112.7) and (115.11,115.81) .. (109,121) ;
\draw   (127.19,80.14) .. controls (121.16,74.86) and (115.1,71.66) .. (109,70.52) .. controls (115.13,69.59) and (121.29,66.57) .. (127.48,61.48) ;
\draw   (277,140.33) .. controls (283.11,145.52) and (289.22,148.63) .. (295.33,149.67) .. controls (289.22,150.7) and (283.11,153.81) .. (277,159) ;
\draw   (421,140.33) .. controls (427.11,145.52) and (433.22,148.63) .. (439.33,149.67) .. controls (433.22,150.7) and (427.11,153.81) .. (421,159) ;

\draw (116,43.4) node [anchor=north west][inner sep=0.75pt]    {$e$};
\draw (114,86.4) node [anchor=north west][inner sep=0.75pt]    {$e_{1}$};
\draw (114,124.4) node [anchor=north west][inner sep=0.75pt]    {$e_{2}$};
\draw (356,90.4) node [anchor=north west][inner sep=0.75pt]    {$e$};
\draw (279,52.4) node [anchor=north west][inner sep=0.75pt]    {$e_{1}$};
\draw (425,50.4) node [anchor=north west][inner sep=0.75pt]    {$e_{2}$};
\draw (31,102.4) node [anchor=north west][inner sep=0.75pt]    {$P_{0}$};
\draw (324,131.4) node [anchor=north west][inner sep=0.75pt]    {$P_{0}$};
\draw (188,102.4) node [anchor=north west][inner sep=0.75pt]    {$P_{1}$};
\draw (374,132.4) node [anchor=north west][inner sep=0.75pt]    {$P_{1}$};

\end{tikzpicture}

         \caption{A curve of type theta on the left and of type dumbbell on the right.}
         \label{figure_conventionforthetaanddumbbell}
     \end{figure}
     and for explicit computations basis
     \begin{itemize}
         \item $(B_1,B_2):=(e +e_2,e_2-e_1)$ of $\H1(\Gamma,\mathbb{Z})$, if $\Gamma$ is of type theta, and $(B_1,B_2)$ $(B_1,B_2):=(e_1,e_2)$, if $\Gamma$ is of type dumbbell. 
         \item in both cases the canonically associated basis $(\omega_1,\omega_2)$ of $\Omega^1_{\Gamma}(\mathbb{Z})$ (see \cite{arXiv:2410.13459}, Section 5).
     \end{itemize}
 \end{convention}
 \subsection{From the perspective of moduli spaces}\label{subsection_tropicaltorelli}
  The connection between $\mathbb{T}\mathcal{A}$ and $\mathbb{T}\mathcal{C}$ focuses on individual objects. Transferred to the level of moduli spaces, it gives rise to the \emph{tropical Torelli map}
 \begin{align}
     t_2^{tr}: M_2^{tr}\rightarrow A_2^{tr}, \Gamma \mapsto \Jac(\Gamma)
 \end{align}
 where $M_2^{tr}$ and $A_2^{tr}$ denote the moduli space of tropical curves of genus $2$ and the moduli space of pptav of dimension $2$, respectively (\cite{MR2641941}, \cite{MR2739784}). 
 For these moduli spaces to be well-behaved, we need slightly more general notions of tropical curves, their Jacobians and correspondingly of pptav. For example, to guarantee that $M_2^{tr}$ is both complete (i.e. closed under specialization) and relates well to its algebraic counterpart $\overline{\mathcal{M}}_2$ (via taking dual graphs) (see \cite{MR2968636}, Remark 2.4), we need to generalize Definition \ref{definition_tropicalcurve} and allow our curves to carry genus at their vertices: A tropical curve $\Gamma$ will be a metric graph $(G,l)$ as in Definition \ref{definition_tropicalcurve} together with a genus function $g: V(G) \rightarrow \mathbb{Z}_{\geqslant 0}$ that satisfies a certain stability condition (see \cite{MR2968636}, Definition 2.1). To include this new piece of data, we define $(G,g)$ to be the \emph{combinatorial type} of $\Gamma$ and the number $g(G) + \sum_{v\in V(G)} g(v)$ to be its genus. The definition of the Jacobian extends to this case and results in a slightly generalized notion of pptav (by relaxing the positive definite-condition to positive semidefinite and renouncing to the non-degeneracy of $\int_\cdot \cdot$ ). 
 We refer to \cite{MR2968636} (Section 4 and 6) or \cite{MR3752493} (Section 4 and 5) for more details as it will be of no further relevance to us.

Now, the tropical Torelli-map $t_2^{tr}$ is a morphism in the category of \emph{stacky fans}, a category that is suitable for the construction of tropical moduli spaces. Naturally, $A_2^{tr}$ and $M_2^{tr}$ are both objects therein.

 A pptav $\Sigma$ is a real integral torus that carries a pp $\zeta$, alternatively a pair $(\frac{\Hom(\Lambda,\mathbb{R})}{\Lambda'}, Q)$, where $Q$ is the quadratic form associated to the symmetric bilinear form $[\zeta(\cdot),\cdot]$. After suitable choice of basis one can find a representative of the isomorphism class of $\Sigma$ of the form $(\mathbb{R}^2/\mathbb{Z}^2,Q)$, where $Q$ is an element of $\tilde{S}_{\geqslant 0}^2$, the space of positive semidefinite $ 2 \times 2$ matrices with rational nullspace. Note that in $A_2^{tr}$ the notion of isomorphism is \emph{stronger} than in $\mathbb{T}\mathcal{A}$: An isomorphism of pptavs (in $A_2^{tr}$) is (in $\mathbb{T}\mathcal{A}$) an isomorphism $f: \Sigma_1 \rightarrow \Sigma_2$ of integral tori such that $f$ is polarized with respect to the pp $\zeta_1$ on $\Sigma_1$ and $\zeta_2$ on $\Sigma_2$ (see Definition \ref{definition_induced/pf/pbpolarizationandpolarizedisogeny} and \cite{röhrle2024tropicalngonalconstruction}, Definition 4.10).

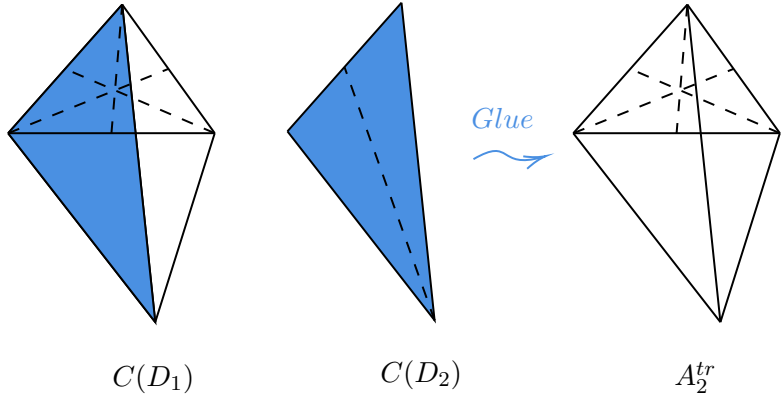
\begin{figure}[hbt!]
     \centering
    \tikzset{every picture/.style={line width=0.75pt}} 

\begin{tikzpicture}[x=0.75pt,y=0.75pt,yscale=-1,xscale=1]

\draw [color={rgb, 255:red, 0; green, 0; blue, 0 }  ,draw opacity=1 ][fill={rgb, 255:red, 245; green, 166; blue, 35 }  ,fill opacity=1 ]   (359.67,82.28) -- (376.34,242.56) ;
\draw [color={rgb, 255:red, 0; green, 0; blue, 0 }  ,draw opacity=1 ][fill={rgb, 255:red, 245; green, 166; blue, 35 }  ,fill opacity=1 ]   (406.33,147.19) -- (359.67,82.28) ;
\draw [color={rgb, 255:red, 0; green, 0; blue, 0 }  ,draw opacity=1 ]   (302.14,147.21) -- (376.34,242.56) ;
\draw [color={rgb, 255:red, 0; green, 0; blue, 0 }  ,draw opacity=1 ]   (406.33,147.19) -- (354.24,147.2) -- (302.14,147.21) ;
\draw    (359.67,82.28) -- (330.9,114.75) -- (302.14,147.21) ;
\draw  [dash pattern={on 4.5pt off 4.5pt}]  (302.14,147.21) -- (383,114.73) ;
\draw  [dash pattern={on 4.5pt off 4.5pt}]  (406.33,147.19) -- (330.9,114.75) ;
\draw  [dash pattern={on 4.5pt off 4.5pt}]  (359.67,82.28) -- (354.24,147.2) ;
\draw    (406.33,147.19) -- (376.34,242.56) ;
\draw [color={rgb, 255:red, 0; green, 0; blue, 0 }  ,draw opacity=1 ][fill={rgb, 255:red, 245; green, 166; blue, 35 }  ,fill opacity=1 ]   (74.64,82.28) -- (91.31,242.56) ;
\draw [color={rgb, 255:red, 0; green, 0; blue, 0 }  ,draw opacity=1 ][fill={rgb, 255:red, 245; green, 166; blue, 35 }  ,fill opacity=1 ]   (121.31,147.19) -- (74.64,82.28) ;
\draw [color={rgb, 255:red, 0; green, 0; blue, 0 }  ,draw opacity=1 ]   (17.11,147.21) -- (91.31,242.56) ;
\draw    (74.64,82.28) -- (45.87,114.75) -- (17.11,147.21) ;
\draw  [dash pattern={on 4.5pt off 4.5pt}]  (74.64,82.28) -- (69.21,147.2) ;
\draw    (121.31,147.19) -- (91.31,242.56) ;
\draw [color={rgb, 255:red, 0; green, 0; blue, 0 }  ,draw opacity=1 ][fill={rgb, 255:red, 74; green, 144; blue, 226 }  ,fill opacity=1 ]   (17.11,147.21) -- (91.31,242.56) -- (77.57,110.5) -- (74.64,82.28) ;
\draw [color={rgb, 255:red, 0; green, 0; blue, 0 }  ,draw opacity=1 ]   (121.31,147.19) -- (69.21,147.2) -- (17.11,147.21) ;
\draw  [dash pattern={on 4.5pt off 4.5pt}]  (121.31,147.19) -- (45.87,114.75) ;
\draw  [dash pattern={on 4.5pt off 4.5pt}]  (17.11,147.21) -- (97.97,114.73) ;
\draw [color={rgb, 255:red, 0; green, 0; blue, 0 }  ,draw opacity=1 ][fill={rgb, 255:red, 74; green, 144; blue, 226 }  ,fill opacity=1 ]   (157.98,146.27) -- (232.18,241.62) -- (222.46,148.15) -- (215.5,81.33) ;
\draw  [dash pattern={on 4.5pt off 4.5pt}]  (232.18,241.62) -- (186.74,113.8) ;
\draw    (74.64,82.28) -- (17.11,147.21) ;
\draw    (215.5,81.33) -- (157.98,146.27) ;
\draw  [dash pattern={on 4.5pt off 4.5pt}]  (69.21,147.2) -- (74.64,82.28) ;
\draw [color={rgb, 255:red, 74; green, 144; blue, 226 }  ,draw opacity=1 ]   (251.8,160.31) .. controls (271.75,148.76) and (263.63,165.74) .. (286.5,159.28) ;
\draw [shift={(288.32,158.73)}, rotate = 162.55] [color={rgb, 255:red, 74; green, 144; blue, 226 }  ,draw opacity=1 ][line width=0.75]    (10.93,-3.29) .. controls (6.95,-1.4) and (3.31,-0.3) .. (0,0) .. controls (3.31,0.3) and (6.95,1.4) .. (10.93,3.29)   ;

\draw (351.74,261.09) node [anchor=north west][inner sep=0.75pt]    {$A_{2}^{tr}$};
\draw (203.35,260.23) node [anchor=north west][inner sep=0.75pt]    {$C( D_{2})$};
\draw (68.25,262.12) node [anchor=north west][inner sep=0.75pt]    {$C( D_{1})$};
\draw (249.39,133.12) node [anchor=north west][inner sep=0.75pt]    {$\textcolor[rgb]{0.29,0.56,0.89}{Glue}$};

\end{tikzpicture}

     \caption{Visualization of $A_2^{tr}$ as in \cite{MR2968636} with dashed lines indicating symmetry.}
     \label{figure_VisualizationofA}
 \end{figure}

 Moreover, $\Sigma_1$ and $\Sigma_2$ are isomorphic if and only if for their selected representatives, $(\mathbb{R}^2/\mathbb{Z}^2,Q_1)$ and $(\mathbb{R}^2/\mathbb{Z}^2,Q_2)$, there exists $X\in Gl(\mathbb{Z})$ with $Q_2=X^tQ_1X$. This suggests the quotient $\tilde{S}_{\geqslant 0}^2/ GL(\mathbb{Z})$, where the action of $GL(\mathbb{Z})$ on $\tilde{S}_{\geqslant 0}^2$ is given by $GL(\mathbb{Z}) \times \tilde{S}_{\geqslant 0}^2 \ni (X,Q) \mapsto X \sbullet Q:= X^TQX$, as a candidate for $A_2^{tr}$. Indeed, there is a point-wise bijection between the two, which leads us to work with this interpretation of $A_2^{tr}$ in Subsection \ref{subsubssection_step2concrete}.
More precisely, $A_2^{tr}$ is a stacky fan with cells 
   \begin{align}
      & C(D_1):=\bar{\sigma}_{D_1}/\Stab({\sigma}_{D_1}), \text{ where } \bar{\sigma}_{D_1}:=\langle \begin{pmatrix}
          1 & -1\\ -1 & 1
      \end{pmatrix},\begin{pmatrix}
          1 & 0\\ 0 & 0
      \end{pmatrix},\begin{pmatrix}
          0 & 0\\ 0 & 1
      \end{pmatrix}\rangle_{\mathbb{R} \geqslant 0} \\
     & C(D_2):=\bar{\sigma}_{D_2}/\Stab({\sigma}_{D_2}), \text{ where } \bar{\sigma}_{D_2}:=\langle \begin{pmatrix}
          1 & 0\\ 0 & 0
      \end{pmatrix},\begin{pmatrix}
          0 & 0\\ 0 & 1
      \end{pmatrix}\rangle_{\mathbb{R} \geqslant 0},\\ 
     &  C(D_3):=\bar{\sigma}_{D_3}/\Stab({\sigma}_{D_3}), \text{ where } \bar{\sigma}_{D_3}:=\langle \begin{pmatrix}
          1 & 0\\ 0 & 0
      \end{pmatrix}\rangle_{\mathbb{R} \geqslant 0}, \\
      & C(D_4):=\bar{\sigma}_{D_4}/\Stab({\sigma}_{D_4}):=\{0\},
   \end{align}
   where $\Stab({\sigma}_{D_i})\leqslant Gl(\mathbb{Z)}$ is the setwise stabilizer of ${\sigma}_{D_i}$. These are
   glued together according to the equivalence relation
   \begin{align}
       Q_1 \sim Q_2 \Leftrightarrow X^tQ_1X= Q_2 \text{ for some } X\in Gl(\mathbb{Z}).
   \end{align}
 As pointed out in \cite{MR2968636} (Section 4), gluing $C(D_i)$ for $i=2,3,4$ does not change $C(D_1)$, so that $A_2^{tr}$ is homeomorphic to $C(D_1)$ (see Figure \ref{figure_VisualizationofA}).
 
\begin{figure}[H]
     \centering
   \tikzset{every picture/.style={line width=0.75pt}} 

\begin{tikzpicture}[x=0.75pt,y=0.75pt,yscale=-1,xscale=1]

\draw [color={rgb, 255:red, 0; green, 0; blue, 0 }  ,draw opacity=1 ][fill={rgb, 255:red, 74; green, 144; blue, 226 }  ,fill opacity=1 ]   (213.38,123.5) -- (304.12,191.37) -- (226.95,35.41) ;
\draw [color={rgb, 255:red, 0; green, 0; blue, 0 }  ,draw opacity=1 ][fill={rgb, 255:red, 245; green, 166; blue, 35 }  ,fill opacity=1 ]   (226.95,35.41) -- (304.12,191.37) ;
\draw [color={rgb, 255:red, 0; green, 0; blue, 0 }  ,draw opacity=1 ][fill={rgb, 255:red, 245; green, 166; blue, 35 }  ,fill opacity=1 ] [dash pattern={on 0.84pt off 2.51pt}]  (213.38,123.5) -- (220.16,79.46) -- (226.95,35.41) ;
\draw [color={rgb, 255:red, 0; green, 0; blue, 0 }  ,draw opacity=1 ]   (330.4,58.44) -- (304.12,191.37) ;
\draw [color={rgb, 255:red, 0; green, 0; blue, 0 }  ,draw opacity=1 ]   (213.38,123.5) -- (330.4,58.44) ;
\draw    (226.95,35.41) -- (330.4,58.44) ;
\draw    (391.33,95) -- (419.33,219) ;
\draw    (391.33,95) -- (421.17,70.5) -- (451,46) ;
\draw [color={rgb, 255:red, 74; green, 144; blue, 226 }  ,draw opacity=1 ]   (451,46) -- (455.17,87.5) -- (459.33,129) ;
\draw  [color={rgb, 255:red, 0; green, 0; blue, 0 }  ,draw opacity=1 ] (282.1,249.4) .. controls (282.1,257.41) and (288.23,263.92) .. (295.78,263.92) .. controls (303.34,263.92) and (309.47,257.41) .. (309.47,249.4) .. controls (309.47,241.38) and (303.34,234.88) .. (295.78,234.88) .. controls (288.23,234.88) and (282.1,241.38) .. (282.1,249.4) -- cycle ;
\draw  [color={rgb, 255:red, 0; green, 0; blue, 0 }  ,draw opacity=1 ] (230.48,248.54) .. controls (230.48,256.56) and (236.61,263.06) .. (244.17,263.06) .. controls (251.73,263.06) and (257.85,256.56) .. (257.85,248.54) .. controls (257.85,240.52) and (251.73,234.02) .. (244.17,234.02) .. controls (236.61,234.02) and (230.48,240.52) .. (230.48,248.54) -- cycle ;
\draw [color={rgb, 255:red, 0; green, 0; blue, 0 }  ,draw opacity=1 ]   (257.85,248.54) -- (282.1,249.4) ;
\draw   (73.72,246.7) .. controls (73.72,238.78) and (86.49,232.36) .. (102.26,232.36) .. controls (118.02,232.36) and (130.8,238.78) .. (130.8,246.7) .. controls (130.8,254.63) and (118.02,261.05) .. (102.26,261.05) .. controls (86.49,261.05) and (73.72,254.63) .. (73.72,246.7) -- cycle ;
\draw    (73.72,246.7) -- (130.8,246.7) ;
\draw [color={rgb, 255:red, 74; green, 144; blue, 226 }  ,draw opacity=1 ]   (340,118.1) .. controls (366.35,105.81) and (356.26,124.03) .. (384.55,116.9) ;
\draw [shift={(386.33,116.43)}, rotate = 164.67] [color={rgb, 255:red, 74; green, 144; blue, 226 }  ,draw opacity=1 ][line width=0.75]    (10.93,-3.29) .. controls (6.95,-1.4) and (3.31,-0.3) .. (0,0) .. controls (3.31,0.3) and (6.95,1.4) .. (10.93,3.29)   ;
\draw [color={rgb, 255:red, 0; green, 0; blue, 0 }  ,draw opacity=1 ][fill={rgb, 255:red, 74; green, 144; blue, 226 }  ,fill opacity=1 ]   (451,46) -- (419.33,219) -- (459.33,129) ;
\draw    (391.33,95) -- (425.33,112) -- (459.33,129) ;
\draw [color={rgb, 255:red, 0; green, 0; blue, 0 }  ,draw opacity=1 ][fill={rgb, 255:red, 74; green, 144; blue, 226 }  ,fill opacity=1 ]   (163.73,123.44) -- (78.53,188.73) -- (153.2,35) ;
\draw [color={rgb, 255:red, 0; green, 0; blue, 0 }  ,draw opacity=1 ][fill={rgb, 255:red, 245; green, 166; blue, 35 }  ,fill opacity=1 ]   (153.2,35) -- (78.53,188.73) ;
\draw [color={rgb, 255:red, 0; green, 0; blue, 0 }  ,draw opacity=1 ][fill={rgb, 255:red, 245; green, 166; blue, 35 }  ,fill opacity=1 ] [dash pattern={on 0.84pt off 2.51pt}]  (163.73,123.44) -- (153.2,35) ;
\draw [color={rgb, 255:red, 0; green, 0; blue, 0 }  ,draw opacity=1 ]   (57.29,55.11) -- (78.53,188.73) ;
\draw [color={rgb, 255:red, 0; green, 0; blue, 0 }  ,draw opacity=1 ]   (163.73,123.44) -- (110.51,89.28) -- (57.29,55.11) ;
\draw    (153.2,35) -- (105.24,45.06) -- (57.29,55.11) ;
\draw [color={rgb, 255:red, 0; green, 0; blue, 0 }  ,draw opacity=1 ]   (532.33,83) -- (419.33,219) ;
\draw [color={rgb, 255:red, 0; green, 0; blue, 0 }  ,draw opacity=1 ]   (459.33,129) -- (532.33,83) ;
\draw    (452.7,46.9) -- (532.33,83) ;
\draw  [dash pattern={on 4.5pt off 4.5pt}]  (57.29,55.11) -- (158.46,79.22) ;
\draw  [dash pattern={on 4.5pt off 4.5pt}]  (163.73,123.44) -- (105.24,45.06) ;
\draw  [dash pattern={on 4.5pt off 4.5pt}]  (153.2,35) -- (110.51,89.28) ;
\draw  [dash pattern={on 4.5pt off 4.5pt}]  (220.16,79.46) -- (330.4,58.44) ;
\draw  [dash pattern={on 4.5pt off 4.5pt}]  (391.33,95) -- (532.33,83) ;
\draw  [dash pattern={on 4.5pt off 4.5pt}]  (459.33,129) -- (421.17,70.5) ;
\draw  [dash pattern={on 4.5pt off 4.5pt}]  (452.7,46.9) -- (425.33,112) ;

\draw (315.57,145.43) node [anchor=north west][inner sep=0.75pt]  [rotate=-359.18]  {$l( e)$};
\draw (232.77,163.81) node [anchor=north west][inner sep=0.75pt]  [rotate=-359.93]  {$l( e_{1})$};
\draw (77.29,108.75) node [anchor=north west][inner sep=0.75pt]    {$l( e_{2})$};
\draw (39.53,117.37) node [anchor=north west][inner sep=0.75pt]    {$l( e)$};
\draw (114.94,159.86) node [anchor=north west][inner sep=0.75pt]    {$l( e_{1})$};
\draw (277.66,114.62) node [anchor=north west][inner sep=0.75pt]  [rotate=-1.41]  {$l( e_{2})$};
\draw (315.4,235) node [anchor=north west][inner sep=0.75pt]    {$e_{2}$};
\draw (212.16,235) node [anchor=north west][inner sep=0.75pt]    {$e_{1}$};
\draw (266.45,235) node [anchor=north west][inner sep=0.75pt]    {$e$};
\draw (95,220) node [anchor=north west][inner sep=0.75pt]    {$e$};
\draw (95,235) node [anchor=north west][inner sep=0.75pt]    {$e_{1}$};
\draw (95,262.88) node [anchor=north west][inner sep=0.75pt]    {$e_{2}$};
\draw (348.72,95) node [anchor=north west][inner sep=0.75pt]    {$\textcolor[rgb]{0.29,0.56,0.89}{Glue}$};
\draw (424,241.4) node [anchor=north west][inner sep=0.75pt]    {$M_{2}^{tr}$};

\end{tikzpicture}

     \caption{Visualization of $M_2^{tr}$ as in \cite{MR2968636} with dashed lines indicating symmetries.}
     \label{figure_VisualizationofM}
 \end{figure}
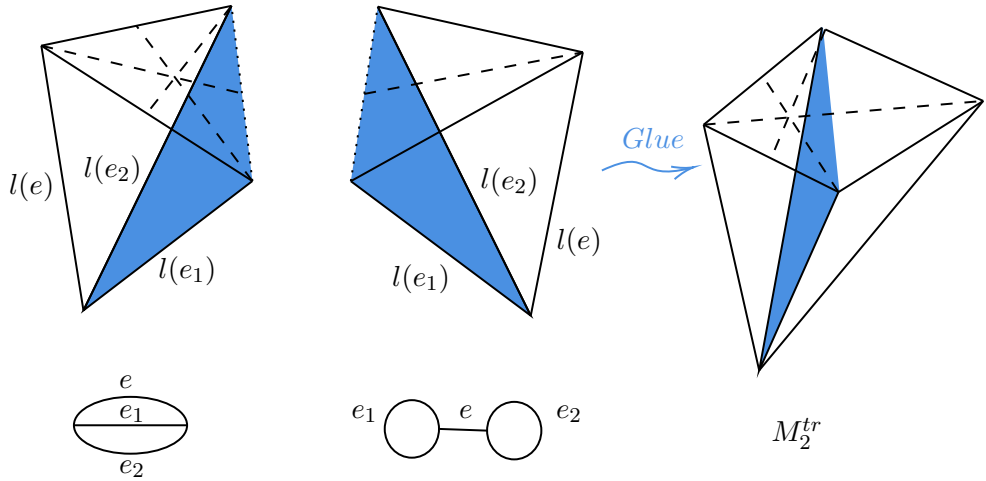

 The  stacky fan $M_2^{tr}$ (see \cite{MR2968636}, Definition and Theorem 3.4) consists of cells 
 \begin{align}
     C(G,g):=\mathbb{R}_{\mathbb{R} \geqslant 0}^{|E(G)|}/\Aut(G,\omega)
 \end{align}
 for each combinatorial type $(G,g)$, where $\Aut(G,g)$ is the automorphism group, i.e. the set of graph automorphisms that preserve $g$, and gluing determining equivalence relation: the relation of specialization (see Figure \ref{figure_VisualizationofM}). In this setting, the tropical Torelli map sends $\Gamma$ to its tropical period matrix (\cite{MR3752493}, Definition 5.1). 
 \section{Constructing Curves of genus 2 covering Curves of genus 1 }\label{section_reconstruction}
 Let $\Gamma$ be a curve of genus $2$. Following \cite{arXiv:2410.13459} we can compute a splitting of $\Jac(\Gamma)$ whenever $\Gamma$ covers an elliptic curve. Such a splitting is given by an isogeny $\phi:\mathbb{T}E' \bigoplus \mathbb{T}E \rightarrow \Jac(\Gamma)$ from which we only want to record the \emph{splitting data}, i.e. $(\mathbb{T}E,\mathbb{T}E',\ker(\phi))$. In this section we reverse the procedure and ask, given $\mathbb{T}E$ and $\mathbb{T}E'$, together with a finite subgroup $G$ of their direct product:
 \begin{itemize}
     \item When is $(\mathbb{T}E,\mathbb{T}E',G)$ actual splitting data?
     \item If it is, is it minimal data, i.e. both necessary and sufficient, for constructing a curve of genus $2$ covering $\mathbb{T}E$ and $\mathbb{T}E'$? 
 \end{itemize}
 
 \subsection{Planning the reconstruction process} 
 Results of \cite{arXiv:2410.13459} (see Section \ref{section_intro}) determine our strategy.

 \begin{reconstruction}\label{Masterplan}
     Given splitting data $(\mathbb{T}E,\mathbb{T}E',G)$, proceed as follows:
\begin{enumerate}
    \item Determine a splitting $\phi: \mathbb{T}E \oplus \mathbb{T}E' \rightarrow J $ and generate a diagram $D_\phi$ modelled on Diagram \ref{diagram_fromctocultimately}.
    \item Construct $\Gamma$.
    \item Define covers $\varphi$ and $\varphi'$.
\end{enumerate}
 \end{reconstruction}
\subsubsection{The properties of splitting data}
We address the feasibility of Plan \ref{Masterplan} first and ask whether splitting data captures the essence of Diagram \ref{diagram_fromctocultimately}, i.e. the data that is both, necessary and sufficient, for its reconstruction.

A closer look reveals: The lower triangle (corresponding to the product property $\mathbb{T}E' \bigoplus \mathbb{T}E$) arises from the upper one (the one corresponding to the coproduct property $\mathbb{T}E' \bigoplus \mathbb{T}E$) as image through the dualization functor, $\widecheck{\cdot}: \mathbb{T}\mathcal{A} \rightarrow \mathbb{T}\mathcal{A}$ (see Definition \ref{lemma_dualizingexactsequences}), and by subsequent identification of the duals with the original spaces. The data required for its reconstruction is therefore: An isogeny $\phi$ and canonical isomorphisms to the duals. We mimic this procedure to construct a diagram $D_\phi$ associated to $\phi$.

\begin{construction}\label{construction_diagraminducedbyisogeny}
     Let $J$ be a pptav and $\phi:\mathbb{T}E' \bigoplus \mathbb{T}E \rightarrow J$ an isogeny  with principal polarizations $\zeta_{\mathbb{T}E' \bigoplus \mathbb{T}E}$ on $\mathbb{T}E' \bigoplus \mathbb{T}E$ (Definition \ref{definition_productandcoproducts}) and $\zeta$ on $J$, both providing isomorphisms to the dual varieties. Start building $D_\phi$ from $\phi$ by drawing arrows that correspond to the canonical injections and projections, $\iota_i$ and $p_i$. These (drawn as solid arrows in Diagram \ref{diagram_inducedbyisogeny}) form its skeleton. Next, complete the upper part to a commutative triangle and set $f_i:=\phi \circ \iota_i$ for $i=1,2$. For the base triangle, let $\tilde{\phi}$ be the adjoint of $\phi$ (with respect to the polarizations $\zeta_{\mathbb{T}E' \bigoplus \mathbb{T}E}$ and $\zeta$ as in Definition \ref{definition_adjoint}) and define the maps $g_i$ by $p_i \circ \tilde{\phi}$ to finish the construction.
 
\begin{diagram}\label{diagram_inducedbyisogeny}
     \mathbb{T}E'  \ar[dr,"f_1",dashed] \ar[r,"\iota_1", hook] & \mathbb{T}E' \bigoplus \mathbb{T}E \ar[d,"{\phi}" ] &\mathbb{T}E \ar[dl,"f_2",dashed,sloped] \ar[l,"\iota_1", hook', swap] \\
        & \ar[dl,"g_1",dashed,sloped] J \ar[d,dashed,"{\tilde{\phi}}"] \ar[dr,"g_2", dashed] &  \\
        \mathbb{T}E'  & \ar[l, "p_1"] \mathbb{T}E' \bigoplus \mathbb{T}E \ar[r, "p_2",swap]  &\mathbb{T}E  
\end{diagram}
\end{construction}
\begin{lemma}\label{lemma_reqfordiagrwithexactdiag}
    Let $\phi:\mathbb{T}E' \bigoplus \mathbb{T}E \rightarrow J$ be an isogeny to a pptav $J$. The diagram $D_\phi$ as in (\ref{diagram_inducedbyisogeny}) has exact diagonals if and only if there exists an integer $d$ such that the following holds:
    \begin{itemize}
        \item $\ker(\phi)\cong \{(\lambda,\alpha(\lambda)): \lambda\in \mathbb{T}E'[d]\}$, where $\alpha: \mathbb{T}E'[d] \rightarrow \mathbb{T}E[d]$ is an isomorphism.
        \item $\phi^*\zeta= (m_d, m_d)$, where $\zeta$ is the principal polarization on $J$.
    \end{itemize}
\end{lemma}
We will use the following characterization of subgroups of the direct product of two groups.
\begin{lemma}(Goursat's Lemma, \cite{zbMATH02692043})\label{lemma_Goursat}
    Let $G_1,G_2$ be groups and $H\leqslant G_1 \oplus G_2$ a subgroup such that the projections $p_1:G_1 \oplus G_2 \rightarrow G_1 $ and $p_2:G_1 \oplus G_2 \rightarrow G_2 $ remain surjective when restricted to $H$. If $N_i$ is the kernel of the restriction of $p_i$ to $H$, then the image of $H$ in the direct product $G_1/N_2 \bigoplus G_2/N_1$ is the graph of an isomorphism
   \begin{align}
       \alpha: G_1/N_2 \rightarrow G_2/N_1.
   \end{align}
\end{lemma}
\begin{proof}[Proof of Lemma \ref{lemma_reqfordiagrwithexactdiag}]
    Let $\phi:\mathbb{T}E' \bigoplus \mathbb{T}E \rightarrow J$ be an isogeny between pptavs and $\tilde{\phi}$ its adjoint with respect to the polarizations $\zeta_{\mathbb{T}E' \bigoplus \mathbb{T}E}$ and $\zeta$ (Definition \ref{definition_adjoint}). For $D_\phi$ to have exact diagonals means that for $i,j\in \{1,2\}$
    \begin{enumerate}
        \item $f_i$ is injective.
        \item $g_i$ is surjective.
        \item $\ker(g_i)=\im(f_j)$ holds whenever $i\neq j$.
    \end{enumerate}
    Recall from Construction \ref{construction_diagraminducedbyisogeny} that for $i\neq j$, $g_j$ is given by $\widecheck{f_i}$ (up to pre- and post-composition with isomorphisms). Then $f_i$ being injective implies that $\widecheck{f_i}$ is surjective and exactness of both diagonals is already equivalent to $(1)$ and $(3)$.
  
    We express these in terms of requirements on $\phi$ that match the statement of Lemma \ref{lemma_reqfordiagrwithexactdiag}. As 
    \begin{align}
        \ker(f_i)=\ker(\phi \circ \iota_i)= \ker(\phi)\cap \im(\iota_i),
    \end{align}
    we have:
    \begin{align}
       (1) \thinspace \Leftrightarrow \thinspace \ker(\phi)\cap \mathbb{T}E' \bigoplus \{0\} = \{0\} \text{ and } \ker(\phi)\cap \{0\} \bigoplus \mathbb{T}E = \{0\}.
    \end{align}
    For $(3)$, we address $"\supset"$ first and observe:
    \begin{align}
       g_2 \circ f_1=0 \text{ and } g_1 \circ f_2=0 &\thinspace \Leftrightarrow \thinspace &&\tilde{\phi}\circ \phi(\mathbb{T}E' \bigoplus \{0\})\subset \mathbb{T}E' \bigoplus \{0\} \text{ and }\\
       &  &&\tilde{\phi}\circ \phi(\{0\} \bigoplus \mathbb{T}E) \subset \{0\} \bigoplus \mathbb{T}E\\
        &\thinspace \Leftrightarrow \thinspace &&\phi^*\zeta=(m_{\alpha_1}, m_{\alpha_2}) \text{ for some } \alpha_1,\alpha_2 \in \mathbb{N},
    \end{align}
   where the last equivalence holds because of the following arguments: From Lemma \ref{lemma_propertiesofadjoints} we know that the composition of $\phi$ with its adjoint, $\tilde{\phi}$, is given by 
   \begin{align}
      \tilde{\phi} \circ \phi = (\phi^*\zeta \circ \zeta^{-1}_{\mathbb{T}E' \bigoplus \mathbb{T}E}, \zeta^{-1}_{\mathbb{T}E' \bigoplus \mathbb{T}E} \circ \phi^*\zeta)
   \end{align}
   So $\tilde{\phi}\circ \phi$ respects the product $\mathbb{T}E' \bigoplus \mathbb{T}E$ (i.e. $\tilde{\phi}\circ \phi$ is a product of morphisms between the individual factors) if and only if its lift $\Hom(\phi^*\zeta \circ \zeta^{-1}_{\mathbb{T}E' \bigoplus \mathbb{T}E})$ does so. This in turn means that $\phi^*\zeta$ itself must be of the form $(h_1,h_2)$ for group homomorphisms $h_1$ and $h_2$, given that $\zeta^{-1}_{\mathbb{T}E' \bigoplus \mathbb{T}E}$ is and that the $\Hom$-functor reflects the property of being a product. Finally, recalling that all lattices involved are $1$-dimensional, yields $h_i=m_{\alpha_i}$, whereby, after suitable choice of basis, we can assume that $\alpha_i\in \mathbb{N}$.\\

   Taken together, $(1)$ and the $"\supset"-$part of $(3)$ determine the shape of $\ker(\phi)$ to some extent: From $\tilde{\phi}\circ \phi=(m_{\alpha_1}, m_{\alpha_2})$ we see $\ker(\tilde{\phi}\circ \phi)=\mathbb{T}E'[\alpha_1] \bigoplus \mathbb{T}E[\alpha_2]$ and $\ker(\phi)$ is a subgroup that projects (via $p_1$ and $p_2$) onto subgroups, say $\mathbb{T}E'[a_1]$, respectively $\mathbb{T}E[a_2]$, of each factor. \\
   Subgroups of direct product that project surjectively onto each factor are characterized by Goursat's lemma (Lemma \ref{lemma_Goursat}): If $N_i$ is the kernel of the restriction of $p_i$ to $\ker(\phi)$, the image of $\ker(\phi)$ in the direct product $\mathbb{T}E'[a_1]/N_2 \bigoplus \mathbb{T}E[a_2]/N_1$ is the graph of an isomorphism
   \begin{align}
       \alpha: \mathbb{T}E'[a_1]/N_2 \rightarrow \mathbb{T}E[a_2]/N_1.
   \end{align}
   We remark that (1) forces both
   \begin{align}
      N_2=(\mathbb{T}E' \bigoplus \{0\})\cap \ker(\phi) \text{ and }
      N_1=(\{0\} \bigoplus \mathbb{T}E) \cap \ker(\phi)
   \end{align}
   to be trivial and thus $d:=a_1=a_2$ as $\alpha$ is an isomorphism.\\
   At this point we should notice that we have almost arrived at the statement of Lemma \ref{lemma_reqfordiagrwithexactdiag}. The only thing left is to take advantage of the $"\subset"-$part of $(3)$ and claim that, under the previous conditions, $"\subset"$ holds if and only if $d=\alpha_1=\alpha_2$ is also true. 
   
   For the "only-if"-direction, suppose (without loss of generality) that $d\neq \alpha_2$. This means that ${p_2}_{|\ker(\phi)}$ does not project surjectively onto $\mathbb{T}E[\alpha_2]$, i.e. that there exists a $\lambda_2\in \mathbb{T}E[\alpha_2] \setminus \mathbb{T}E[d]$. Then $\tilde{\phi}\circ \phi (0,\lambda_2)= (0,0)$ shows $\phi (0,\lambda_2)\in \ker(g_2)$, but $\phi (0,\lambda_2)\not \in \im(f_1)$.
   Otherwise $\phi (0,\lambda_2)=\phi (\lambda_1,0)$ would imply $(-\lambda_1,\lambda_2)\in \ker(\phi)$, a contradiction to $\lambda_2\in \mathbb{T}E[\alpha_2] \setminus \mathbb{T}E[d]$. 
   For the "if"-direction, consider a point $\lambda\in \ker(g_2)$, i.e. $\tilde{\phi}(\lambda)\in \mathbb{T}E' \bigoplus \{0\}$. Since $\phi$ is surjective, we can write $\lambda$ as $\phi(\mu_1,\mu_2)$ for some $(\mu_1,\mu_2)\in \mathbb{T}E' \bigoplus \mathbb{T}E$. Applying $\tilde{\phi}$, yields $\mu_2\in \mathbb{T}E[\alpha_2]$. Then, $\mathbb{T}E[\alpha_2]=\mathbb{T}E[d]$ guarantees that we can find a point $(\tilde{\mu}_1,\mu_2)\in \ker(\phi)$ and obtain
   \begin{align}
       \lambda=\lambda + 0= \phi(\mu_1,\mu_2) - \phi (\tilde{\mu}_1,\mu_2)= \phi(\mu_1-\tilde{\mu}_1,0)\in \im(f_1)
   \end{align}
   as desired. This finishes the proof.
\end{proof}
\begin{notation}
    We now start with Plan \ref{Masterplan}: Take two elliptic curves $\mathbb{T}E'$, $ \mathbb{T}E$ and a finite subgroup $G$ as in Lemma \ref{lemma_reqfordiagrwithexactdiag}, i.e. $G$ is the graph of an isomorphism $\alpha$ between the $d$-torsion points of $\mathbb{T}E'$ and $\mathbb{T}E$, where $d\in \mathbb{N}$.
\end{notation}
 We will continue to use this notation for the remaining sections.
 \subsection{Step 1 of Plan \ref{Masterplan}}\label{subsection_step1}
 A first natural choice for $J$ and $\phi$, suggested by the classical treatment in \cite{MR1085258}, is given by
\begin{align}
  (\mathbb{T}E' \bigoplus \mathbb{T}E)_G \text{ and } q: \mathbb{T}E' \bigoplus \mathbb{T}E \rightarrow (\mathbb{T}E' \bigoplus \mathbb{T}E)_G,
\end{align}
where $(\mathbb{T}E' \bigoplus \mathbb{T}E)_G$ is the quotient of $\mathbb{T}E' \bigoplus \mathbb{T}E$ by $G$, and $q$ the quotient map (Lemma \ref{lemma_quotientoftavbyfinitesubgrouppaper2}). 
\subsubsection{Key observation}
Using the "Cover-to-Splitting" direction discussed in \cite{arXiv:2410.13459} (Theorem \ref{theorem_CovertoSplitting}) to probe the "classical" choice shows: $q$ fails to meet the requirements of Lemma \ref{lemma_reqfordiagrwithexactdiag}. 
This, understood well enough, will eventually show what role $q$ and $(\mathbb{T}E' \bigoplus \mathbb{T}E)_G$ have to play.

In the setting of Theorem \ref{theorem_CovertoSplitting} of \cite{arXiv:2410.13459} we find $q$ to be part of the following commutative triangle 
\begin{diagram}
     \mathbb{T}E' \bigoplus \mathbb{T}E \ar[d,"q"] \ar[r,"\phi"] &  \Jac(\Gamma)  \\
         (\mathbb{T}E' \bigoplus \mathbb{T}E)_G \ar[ur,"{\bar{\phi}}", swap] & 
   \label{diagram_theroleofq}
\end{diagram}
where $G:=\ker(\phi)$ and $\bar{\phi}$ is defined by the universal property of the quotient in $\gls{Ab}$. Breaking down (\ref{diagram_theroleofq}) further, to the level of lattices, shows $\bar{\phi}$ is a morphism of tori:
\begin{diagram}
\Omega^1_{\mathbb{T}E'} (\mathbb{Z}) \bigoplus \Omega^1_{\mathbb{T}E} (\mathbb{Z})  & \arrow[l, "q^\# ",swap]\Omega^1_{\mathbb{T}E'} (\mathbb{Z}) \bigoplus \Omega^1_{\mathbb{T}E} (\mathbb{Z}) & \arrow[l, "\bar{\phi}^\#",swap] \ar[ll,bend left, "\phi^\#"] \Omega^1_{\Gamma} (\mathbb{Z})  \\
\ \\
  \H1(\mathbb{T}E', \mathbb{Z}) \bigoplus \H1(\mathbb{T}E, \mathbb{Z}) \ar[rr,bend right, "\phi_\#",swap] \arrow[r,"q_\#",hook]  &  \pi^{-1}(G) \arrow[r,"\bar{\phi}_\#"]  &  \H1(\Gamma, \mathbb{Z}) \\
\label{diagram_showinglatticesforJandTE'+TE}
  \end{diagram} 
where
\begin{itemize}
    \item $\phi=(\phi^\#,\phi_\#)=((\varphi'_{*},\varphi_{*}),(\varphi^{'*} \oplus \varphi^{*}))$.
    \item  $q=(q^\# ,q_\#)=(Id, \int + \int(\cdot) )$
    \item $\bar{\phi}=(\bar{\phi}^\# ,\bar{\phi}_\#)=(\phi^\#, \int(\cdot)^{-1} \circ \Hom(\phi^\#)_{|\pi^{-1}(G)})$
\end{itemize}

with 
\begin{align}
&\int(\cdot): \H1(\Gamma, \mathbb{Z}) \rightarrow \Hom(\Omega^1_{\Gamma},\mathbb{R}), B \rightarrow \int_{B}\cdot \\ 
&\int + \int(\cdot) : \H1(E', \mathbb{Z}) \bigoplus \H1(E, \mathbb{Z}) \rightarrow \pi^{-1}(G), (B',B) \rightarrow (\int_{B'} + \int_B)(\cdot, \cdot)
\end{align}
the inclusions induced by the natural pairings on $\Jac(\Gamma)$, respectively on $\mathbb{T}E' \bigoplus \mathbb{T}E$, $\pi$ the universal covering of $\mathbb{T}E' \oplus \mathbb{T}E$, and $\int(\cdot)^{-1}$ denotes the isomorphism $\im(\int(\cdot)) \cong \H1(\Gamma,\mathbb{Z})$ given by $\int_{B_i}\cdot \mapsto B_i$
for a lattice Basis $(B_1,B_2)$ of $\H1(\Gamma,\mathbb{Z})$ and the corresponding lattice Basis $(\int_{B_1},\int_{B_2})$ of $\im(\int(\cdot))$. One verifies that $\bar{\phi}_\#$ is a group isomorphism, (\ref{diagram_showinglatticesforJandTE'+TE}) commutes, and $\bar{\phi}$ is compatible with the pairings, using $\pi^{-1}(G)=\Hom(\phi^\#)^{-1}(\H1(\Gamma, \mathbb{Z}))$ (since $\ker(\phi)=G$) and using that $\Hom(\phi^\#)$ is a vector space isomorphism (it is a surjective map between spaces of equal dimension) that agrees with $\phi_\#$ (lifted to a map between the universal covers).

Now both $\mathbb{T}E' \bigoplus \mathbb{T}E$ and $\Jac(\Gamma)$ carry principal polarizations that provide a link between the upper and lower path. We have:
\begin{figure}[H]
    \centering
    \begin{tikzcd}
\Omega^1_{\mathbb{T}E'} (\mathbb{Z}) \bigoplus \Omega^1_{\mathbb{T}E} (\mathbb{Z})  & \arrow[l, "q^\# ",swap]\Omega^1_{\mathbb{T}E'} (\mathbb{Z}) \bigoplus \Omega^1_{\mathbb{T}E} (\mathbb{Z}) & \arrow[l, "\bar{\phi}^\#",swap]  \Omega^1_{\Gamma} (\mathbb{Z})  \\
  \H1(\mathbb{T}E', \mathbb{Z}) \bigoplus \H1 (\mathbb{T}E, \mathbb{Z}) \arrow[u, "m_d \circ \zeta_{\mathbb{T}E' \bigoplus \mathbb{T}E} "] \arrow[r,"q_\#",hook]  &  \pi^{-1}(G) \arrow[r,"\bar{\phi}_\#"] \arrow[u,"\bar{\phi}^*\zeta_{\Gamma}"] &  \H1 (\Gamma, \mathbb{Z}) \arrow[u,"\zeta_{\Gamma}"] . \\
\end{tikzcd}
    \caption{The interaction of the natural pp on $\mathbb{T}E' \bigoplus \mathbb{T}E$ and $\Jac(\Gamma)$, providing a blueprint for $\phi$.}
    \label{diagram_interactionofphiq}
\end{figure}
This brings us to $3$ important observations:
  \begin{enumerate}
      \item By default (meaning with default polarization (see \cite{arXiv:2410.13459}, Lemma 26) $(\mathbb{T}E' \bigoplus \mathbb{T}E)_G$ is \emph{not} a pptav.
      \item The map $\bar{\phi}$ is \emph{not} an isomorphism of tav, though it is, of course, an isomorphism of abelian groups. This factors $\phi$ into a free isogeny $q$ and a dilation $\bar{\phi}$ since $q^\#$ and $\bar{\phi}_\#$ are isomorphisms (compare to \cite{röhrle2024tropicalngonalconstruction}, Lemma 4.9).
      \item The polarization $\bar{\phi}^*\zeta_{\Gamma}$ is \emph{not} principal. 
  \end{enumerate}
    \begin{example}\label{example_analysisofinterplyphiq}
We examine triangle (\ref{diagram_theroleofq}) for the splitting generated by the optimal cover from Figure \ref{figure_optimalandnotoptimalcover} 
with our usual choice of coordinates (for details on the computations see \cite{arXiv:2410.13459}, Section 6). Then
    \begin{diagram}[ampersand replacement=\&]
     \mathbb{R}/\mathbb{Z} \bigoplus \mathbb{R}/3\mathbb{Z} \ar[d,"I_2 "] \ar[r,"{\begin{psmallmatrix}
     1 & 1 \\
     2  & 0
 \end{psmallmatrix}}"] \& \mathbb{R}^2/ \begin{pmatrix}
     2 & 1\\
     1 & 2
 \end{pmatrix}\mathbb{Z}^2  \\
         \mathbb{R}^2/ \langle \begin{pmatrix}
     1 \\
     0 
 \end{pmatrix}, \begin{pmatrix}
    \frac{1}{2}\\
    \frac{3}{2}
 \end{pmatrix} \rangle_\mathbb{Z} \ar[ur,"{\begin{psmallmatrix}
     1 & 1 \\
     2 & 0
 \end{psmallmatrix}}", swap] \&
   \end{diagram} 
breaks down into
\begin{diagram}[ampersand replacement=\&]
     \mathbb{Z}^2    \&[5em] \ar[l,"{M(\phi^\#)=\begin{pmatrix}
     1 & 2 \\
     1 & 0
 \end{pmatrix}}",swap] \mathbb{Z}^2 \ar[dl,"{M(\bar{\phi}^\#)=\begin{pmatrix}
     1 & 2 \\
     1 & 0
 \end{pmatrix}}"] \&  \& \mathbb{Z}^2 \ar[d,"{M(q_\#)=\begin{pmatrix}
     1 & -1 \\
     0 & 2
 \end{pmatrix}}", swap] \ar[r,"{M(\phi_\#)=\begin{pmatrix}
     0 &  2 \\
     1 & -1
 \end{pmatrix}}"] \&[5em] \mathbb{Z}^2  \\ [5ex]
         \mathbb{Z}^2 \ar[u,"M(q^\#)=I_2"]  \&  \& \&  \mathbb{Z}^2 \ar[ur,"{M(\bar{\phi}_\#)=\begin{pmatrix}
     0 & 1 \\
     1 & 0
 \end{pmatrix}}", swap]
\end{diagram}

 This yields coordinate representations for the polarizations $\bar{\phi}^*\zeta_{\Gamma}$ and $\zeta_G$
\begin{align}
   & M(\bar{\phi}^*\zeta_{\Gamma})= M({\phi}^\#)\cdot M(\zeta_{\Gamma}) \cdot M(\bar{\phi}_\#)
    = \begin{pmatrix}
     1 & 2 \\
     1 & 0
 \end{pmatrix}\cdot I_2 \cdot \begin{pmatrix}
     0 & 1 \\
     1 & 0
 \end{pmatrix}= \begin{pmatrix}
     2 & 1 \\
     0 & 1
 \end{pmatrix} \\
 & M(\zeta_G)=M(\check{\zeta}_G)= M({q}_\#)\cdot M(\check{\zeta}_{\mathbb{T}E'\oplus \mathbb{T}E }) \cdot M(q^\#)= \begin{pmatrix}
     1 & -1 \\
     0 & 2
 \end{pmatrix}
\end{align}
that reflect our observations from above:
\begin{enumerate}
    \item $(\mathbb{T}E'\oplus \mathbb{T}E )_G$ with default polarization $\zeta_G$ is \emph{not} a pptav since $\zeta_G$ is of type $(1,2)$.
    \item $\bar{\phi}$ is only a dilation as $\phi^\#$ is not an isomorphism.
    \item $\bar{\phi}^*\zeta_{\Gamma}$ is \emph{not} principal.
\end{enumerate}
 
\end{example}
The previous discussion demonstrates:
  \begin{itemize}
      \item The failure of $q$ to meet the requirements of Lemma \ref{lemma_reqfordiagrwithexactdiag} can be traced back to observations (2) and (3).
      \item The failure of $(\mathbb{T}E'\oplus \mathbb{T}E )_G$ to provide an adequate basis for supporting a principle polarization that interacts "well" with $q$ (see Lemma \ref{lemma_polonquotientJ} for a precise statement) can be traced back to observation (1). It is confirmed by Lemma \ref{lemma_polonquotientJ}.
     
  \end{itemize}
  The previous discussion also demonstrates: These issues can be fixed by using Figure \ref{diagram_interactionofphiq} as a blueprint for $\phi$. We treat $q$ as building block and construct $\phi$ step-wise according to Figure \ref{diagram_interactionofphiq}: The following Lemma takes care of the second vertical arrow.

  \begin{lemma}\label{lemma_polonquotientJ}
  Let $J:=(\mathbb{T}E'\oplus \mathbb{T}E )_G$ and $q: \mathbb{T}E' \bigoplus \mathbb{T}E \rightarrow J$ be the quotient map. Then there exists a unique polarization $\zeta$ on $J$ with $q^*\zeta=m_d \circ \zeta_{\mathbb{T}E' \bigoplus \mathbb{T}E} $. Moreover, $\zeta$ is not principal. 
  \end{lemma}
  \begin{proof}
      We show existence and uniqueness of $\zeta$ using Lemma \ref{lemma_inducingpolarization} and suggest referring to Diagram \ref{diagram_showinglatticesforJandTE'+TE} for an overview of the lattices involved. Since $q^\#$ is an isomorphism, condition (1), that is $\im (q^\#) \supset \im (m_d \circ \zeta_{\mathbb{T}E' \bigoplus \mathbb{T}E})$, is automatically satisfied. For the divisibility condition, condition (2), we determine the invariant factors of 
      \begin{align}
         m_d \circ \zeta_{\mathbb{T}E' \bigoplus \mathbb{T}E}: \H1(\mathbb{T}E', \mathbb{Z}) \bigoplus \H1(\mathbb{T}E, \mathbb{Z}) \rightarrow  \im(q^\#)  \text{ and }  i: \im(q_\#) \hookrightarrow \pi^{-1}(G).
      \end{align}
      The former is just the component-wise multiplication-by-$d$-map whose first and second invariant factor is $d$. For the latter, we first establish a coordinate representation of $\pi^{-1}(G)$ and $\im(q_\#)$ inside $\Hom(\Omega^1_{\mathbb{T}E'}, \mathbb{R}) \oplus \Hom(\Omega^1_{\mathbb{T}E}, \mathbb{R})$ with respect to the dual basis $(\omega,0)^*,(0,\tilde{\omega})^*$. Doing so we identify $\im(q_\#)$ with $\Lambda_1':=\langle \begin{pmatrix}
          l_{\mathbb{T}E'} \\ 0
      \end{pmatrix}, \begin{pmatrix}
         0 \\  l_{\mathbb{T}E}
      \end{pmatrix}\rangle_\mathbb{Z}$, which corresponds to choosing $((B,0),(0,\tilde{B}))$ as basis on $ \H1(\mathbb{T}E', \mathbb{Z}) \bigoplus \H1(\mathbb{T}E, \mathbb{Z})$ and embedding it by means of $\int + \int(\cdot)$. For $\pi^{-1}(G)$, recall that $G$ is the graph of an isomorphism $\alpha$ between cyclic groups of order $d$ (Lemma \ref{lemma_reqfordiagrwithexactdiag}). Under the identification of $\mathbb{T}E$ with $ \mathbb{R}/ (l_{\mathbb{T}E}\mathbb{Z})$, $\frac{l_{\mathbb{T}E}}{d}$ is a generator of $\mathbb{T}E[d]$ whose preimage under $\alpha$, in turn, generates $\mathbb{T}E'[d]\cong \mathbb{R}/ (l_{\mathbb{T}E'}\mathbb{Z})[d] $, i.e. is of the form $\frac{k \cdot l_{\mathbb{T}E'}}{d}$ for some $k\in \mathbb{N}$ with $1\leqslant k \leqslant d-1$ and $\gcd(k,d)=1$. Then, $G=\langle \begin{pmatrix}
          \frac{k \cdot l_{\mathbb{T}E'}}{d} \\ \frac{l_{\mathbb{T}E}}{d}
      \end{pmatrix}\rangle_\mathbb{Z}$ and its lift $\pi^{-1}(G)$ to the universal cover is the $2$-dimensional lattice $G + \Lambda_1'$ generated by $\begin{pmatrix}
          l_{\mathbb{T}E'} \\ 0
      \end{pmatrix}$ and $  \begin{pmatrix}
          \frac{k \cdot l_{\mathbb{T}E'}}{d} \\ \frac{l_{\mathbb{T}E}}{d}
      \end{pmatrix}$. With respect to these we can specify the representation matrix of $i$ since $\begin{pmatrix}
          0 \\ l_{\mathbb{T}E}
      \end{pmatrix}=d\cdot\begin{pmatrix}
          \frac{k \cdot l_{\mathbb{T}E'}}{d} \\ \frac{l_{\mathbb{T}E}}{d}
      \end{pmatrix} - k\cdot \begin{pmatrix}
          l_{\mathbb{T}E'}\\ 0
      \end{pmatrix} $:
      \begin{align}
         M(i)= \begin{pmatrix}
          1 & -k  \\ 0 & d
      \end{pmatrix} 
      \end{align}
      and compute its type $(\gamma_1, \gamma_2)$ by
      \begin{align}
          \gamma_1=\gcd(1,k,d)=1 \text{ and }  \gamma_2=\det(M(i))=d.
      \end{align}
      We see, that $\Det(i)=\gamma_1\cdot \gamma_2=1\cdot d$ is a divisor of the first invariant factor of $m_d\circ \zeta_{\mathbb{T}E' \bigoplus \mathbb{T}E}$. Hence, by Lemma \ref{lemma_inducingpolarization}, $\zeta$ exists and is unique. Moreover, $\zeta$ cannot be principal: Otherwise $q^\#\circ \zeta \circ q_\#$ would have the same type as $q_\#$, $(1,d)$ as previously computed (since $q^\#\circ \zeta$ is an isomorphism), but $q^\#\circ \zeta \circ q_\#=m_d\circ \zeta_{\mathbb{T}E' \bigoplus \mathbb{T}E}$ has type $(d,d)$, a contradiction.
  \end{proof}
  \subsubsection{Determine the splitting $\phi$}
  We now construct $\phi$. For convenience, let $J=(\mathbb{T}E' \bigoplus \mathbb{T}E)_G$, as before, and relabel the target of $\phi$ (see step 1 of Plan \ref{Masterplan}) from $J$ to $J^{pp}$, instead.

  \begin{construction}\label{construction_determinesplitting}
  Let $(\mathbb{T}E',\mathbb{T}E,G)$ be splitting data and $\zeta$ be the unique polarization from Lemma \ref{lemma_polonquotientJ}. For lattice basis $(\omega_1,\omega_2)$ of $\Omega^1_{\mathbb{T}E'} (\mathbb{Z}) \bigoplus \Omega^1_{\mathbb{T}E} (\mathbb{Z})$ and $(B_1,B_2)$ of $\pi^{-1}(G)$  define: 
\begin{itemize}
     \item $\gls{zetapp}:\pi^{-1}(G) \rightarrow  \Omega^1_{\mathbb{T}E'} (\mathbb{Z}) \bigoplus \Omega^1_{\mathbb{T}E} (\mathbb{Z})$ by $B_i\mapsto \omega_i$ for $i=1,2$.
     \item $\bar{\phi}:=(\bar{\phi}^\#,Id)$, where $\bar{\phi}^\#(\omega_i):=\zeta(B_i)$.
    \item  $\gls{Jpp}:=(\Omega^1_{\mathbb{T}E'} (\mathbb{Z}) \bigoplus \Omega^1_{\mathbb{T}E} (\mathbb{Z}),\pi^{-1}(G),[\cdot, \cdot]^{pp})$, where \\ $[\omega, B]^{pp}:= [\bar{\phi}^\#(\omega), B]_G$ for all $B\in \pi^{-1}(G) $ and $\omega \in \Omega^1_{\mathbb{T}E'} (\mathbb{Z}) \bigoplus \Omega^1_{\mathbb{T}E} (\mathbb{Z})$.
\end{itemize}
  \end{construction} 
\begin{proposition}\label{proposition_constructedpairJppandphisatisfiesreq}
    Let  $J^{pp}, \zeta^{pp}$ and $\bar{\phi}$ be defined as in Construction \ref{construction_determinesplitting}. Then, $J^{pp}$ with polarization $\zeta^{pp}$ is a pptav and $\phi:=\bar{\phi} \circ q$ a morphism of tavs satisfying the conditions described in Lemma \ref{lemma_reqfordiagrwithexactdiag}.
\end{proposition}
\begin{proof}
  Proposition \ref{proposition_constructedpairJppandphisatisfiesreq} naturally decomposes into two parts, which we address separately:
  \begin{enumerate}
      \item $J^{pp}$ is a pptav.
      \item $\phi$ satisfies the conditions of Lemma \ref{lemma_reqfordiagrwithexactdiag}.
  \end{enumerate} 
  Because $[\cdot, \cdot]_G$ is a non-degenerate pairing and $\bar{\phi}^\#$ is injective and linear, $[\cdot, \cdot]^{pp}$, which is equal to the composition $[\cdot, \cdot]_G\circ(\bar{\phi}^\#,Id)$, equips the quotient with the structure of an integral torus. To see that $\zeta^{pp}$ even grants $J^{pp}$ the status of a pptav, it is only the "polarized-part" that needs checking (principal it will be, since $\zeta^{pp}$ is a group isomorphism).
 Thus, let us verify that $\zeta^{pp}$ induces a scalar product on $\pi^{-1}(G) \otimes \mathbb{R}$: For $\lambda'_1, \lambda'_2\in \pi^{-1}(G) \otimes \mathbb{R}$ write $\lambda'_1=aB_1+bB_2$ in terms of the basis $(B_1,B_2)$. Unpacking definitions we obtain
  \begin{align}
      [\zeta^{pp}(\lambda'_1), \lambda'_2]^{pp}=& [a\omega_1+b\omega_2, \lambda'_2]^{pp}= [\bar{\phi}^\#(a\omega_1+b\omega_2), \lambda'_2]_G\\
      =&[a\zeta (B_1)+b\zeta(B_1), \lambda'_2]_G=[\zeta (\lambda'_1), \lambda'_2]_G
  \end{align}
   and recognize, in the last step, the bilinear form induced by $\zeta$ on $J$ (with the quotient pairing). Since $\zeta$ is a polarization, this is indeed a scalar product. For part $(2)$, note that by construction
  \begin{align}
      [\bar{\phi}^\#(\lambda_1), \lambda'_2]_G = [\lambda_1,\bar{\phi}_\#(\lambda'_2)]^{pp}
  \end{align}
  holds for all $\lambda'_2\in \pi^{-1}(G)$ and $\lambda_1\in \Omega^1_{\mathbb{T}E'} (\mathbb{Z}) \bigoplus \Omega^1_{\mathbb{T}E} (\mathbb{Z})$. Hence, $\bar{\phi}$ is a morphism, even an isogeny, and $\phi$ as composition of isogenies, is one as well. Target and domain varieties are principally polarized by part $(1)$. Its kernel agrees with $\ker(q)$, since $\bar{\phi}$ is a group isomorphism, and is, therefore, the graph of an isomorphis $\alpha$ as required by Lemma \ref{lemma_reqfordiagrwithexactdiag} (condition 1). Condition 2, finally, holds by construction:
  \begin{align}
      \phi^*\zeta^{pp}=(\bar{\phi}\circ q)^*\zeta^{pp}=q^*(\bar{\phi}^*\zeta^{pp})=q^*(\zeta)=m_d \circ \zeta_{\mathbb{T}E' \bigoplus \mathbb{T}E}.
  \end{align}
\end{proof}
\begin{corollary}
  Let $\phi$ be the map from Proposition \ref{proposition_constructedpairJppandphisatisfiesreq} and $D_\phi$ the diagram associated to $\phi$ as in Construction \ref{construction_diagraminducedbyisogeny}. Then $D_\phi$ has exact diagonals.
\end{corollary}
\begin{proof}
    This follows from Lemma \ref{lemma_reqfordiagrwithexactdiag}.
\end{proof}
\begin{example}\label{example_runningexamplewayhome}
   Let $\mathbb{T}E$ and $\mathbb{T}E'$ be two elliptic curves with $l_{\mathbb{T}E}=3$ and $l_{\mathbb{T}E'}=1$ and $\alpha: \mathbb{T}E'[2] \rightarrow \mathbb{T}E[2] $ the unique isomorphism. We choose $\omega_1:=(\omega,0),\omega_2:=(0,\tilde{\omega})$ as basis of $\Omega^1_{\mathbb{T}E'} (\mathbb{Z}) \bigoplus \Omega^1_{\mathbb{T}E} (\mathbb{Z})$ and let the two generators (of $\pi^{-1}(G)$) play the role of $B_1$ and $B_2$:
   \begin{align}
      \pi^{-1}(G) =  \langle(\int^{Q'}_{P'},\int^{Q}_{P}), (\omega,0)^* \rangle_\mathbb{Z} \subset \Hom(\Omega^1_{\mathbb{T}E'},\mathbb{R})\oplus \Hom(\Omega^1_{\mathbb{T}E},\mathbb{R}),
   \end{align}
   where $Q'\in \mathbb{T}E'$ and $Q\in \mathbb{T}E$ are the only non-trivial points of order two and $P'\in \mathbb{T}E'$ and $P\in \mathbb{T}E$ are fixed reference points (that yield isomorphisms to the respective Jacobians). Under the identification $\mathbb{T}E' \cong \mathbb{R}/\mathbb{Z}$ and $\mathbb{T}E\cong \mathbb{R}/3\mathbb{Z}$ we have
   \begin{align}
       G=\{\begin{pmatrix}
           0\\
           0
       \end{pmatrix}, \begin{pmatrix}
           \frac{1}{2}\\
           \frac{3}{2}
       \end{pmatrix}\} \subset \mathbb{R}/\mathbb{Z} \oplus \mathbb{R}/3\mathbb{Z} \text{ and } \pi^{-1}(G) =  \langle\begin{pmatrix}
           1\\
           0
       \end{pmatrix}, \begin{pmatrix}
           \frac{1}{2}\\
           \frac{3}{2}
       \end{pmatrix} \rangle_\mathbb{Z}.
   \end{align}
   We continue working in coordinates and compute $\zeta$ (see Lemma \ref{lemma_polonquotientJ}), more precisely $_{(B_1,B_2)}M_{(\omega_1,\omega_2)} (\zeta)$, via Algorithm \ref{algorithm_constructionofinducingpolarization}:
   \begin{align}
      _{(B_1,B_2)}M_{(\omega_1,\omega_2)} (\zeta)= \begin{pmatrix}
           2 & 0\\
           0 & 2
       \end{pmatrix}\cdot \begin{pmatrix}
           1 & -1 \\
         0 & 2
       \end{pmatrix}^{-1}= \begin{pmatrix}
           2 & 1 \\
         0 & 1
       \end{pmatrix}.
   \end{align}
   Setting $\zeta^{pp}, \bar{\phi}$ and $J^{pp}$ as described above Proposition \ref{proposition_constructedpairJppandphisatisfiesreq} yields 
   \begin{align}
      _{(B_1,B_2)}M_{(B_1,B_2)} (\bar{\phi}_\#)= I_2 \text{ and }  _{(\omega_1,\omega_2)}M_{(\omega_1,\omega_2)} (\bar{\phi}^\#) = \begin{pmatrix}
           2 & 1 \\
         0 & 1
       \end{pmatrix},
   \end{align}
    such that all maps interact as expected (see Diagram \ref{diagram_interactionofphiq}):
    \begin{diagram}[ampersand replacement=\&]
\mathbb{Z}^2  \& [3em] \arrow[l, "I_2 ",swap]\mathbb{Z}^2 \& [3em] \arrow[l, dashed, "{\begin{pmatrix}
           2 & 1 \\
         0 & 1
       \end{pmatrix}}",swap]  \mathbb{Z}^2  \\ [5ex]
  \mathbb{Z}^2 \arrow[u, "{\begin{pmatrix}
           2 & 0 \\
         0 & 2
       \end{pmatrix}} "] \arrow[r,"{\begin{pmatrix}
           1 & -1 \\
         0 & 2
       \end{pmatrix}}", swap]  \&  \mathbb{Z}^2 \arrow[r,"I_2",swap] \arrow[u, dashed, "{\begin{pmatrix}
           2 & 1 \\
         0 & 1
       \end{pmatrix}}"] \&  \mathbb{Z}^2 \arrow[u,"I_2"] . \\
    \end{diagram}

With coordinate representations at hand, we check that $\zeta^{pp}$ is a polarization simply by computing its Gram matrix
\begin{align}
    M ([\zeta^{pp}(\cdot),\cdot ]^{pp}) = M (\zeta)^tM ([\cdot,\cdot ]_G)=\begin{pmatrix}
           2 & 1 \\
         1 & 2
       \end{pmatrix}
\end{align}
and obtain 
\begin{align}
    \phi: \mathbb{R}/\mathbb{Z} \oplus \mathbb{R}/3\mathbb{Z} \rightarrow \mathbb{R}^2/\begin{pmatrix}
           2 & 1 \\
         1 & 2
       \end{pmatrix}\mathbb{Z}^2, x \mapsto \begin{pmatrix}
           2 & 0 \\
         1 & 1
       \end{pmatrix}x,
\end{align}
an isogeny of pptav whose associated diagram $D_\phi$ has exact diagonals:
 \begin{figure}[H]
   \centering
    \begin{tikzcd}[ampersand replacement = \&]
     \mathbb{R}/\mathbb{Z} \ar[dddd,"\cdot 2",swap] \ar[ddr,"{\begin{psmallmatrix}
     2 \\
     1 \\
    \end{psmallmatrix}}",sloped] \ar[r,"{\begin{psmallmatrix}
     0\\
     1\\
     \end{psmallmatrix}}", hook] \& \mathbb{R}/\mathbb{Z} \bigoplus \mathbb{R}/3\mathbb{Z} \ar[dd,"{\begin{psmallmatrix}
     2 &0 \\
     1 & 1\\
      \end{psmallmatrix}}" description] \& \mathbb{R}/3\mathbb{Z} \ar[ddl,"{\begin{psmallmatrix}
     1  \\
      0\\
     \end{psmallmatrix}}",sloped] \ar[l,"{\begin{psmallmatrix}
     0 \\
     1\\
     \end{psmallmatrix}}", hook', swap] \ar[dddd,"\cdot 2"]\\
      \& \& \& \\
        \& \ar[ddl,"{\begin{psmallmatrix}
     1 & 0 \\
      \end{psmallmatrix}}",sloped] \mathbb{R}^2/ \begin{pmatrix}
     2 & 1\\
     1 & 2\\
      \end{pmatrix}\mathbb{Z}^2 \ar[dd,"{\begin{psmallmatrix}
     1 & 0 \\
     -1 & 2\\
      \end{psmallmatrix}}"description ] \ar[ddr,"{\begin{psmallmatrix}
     -1 & 2\\
      \end{psmallmatrix}}" ,sloped] \&  \\
       \& \& \& \\
        \mathbb{R}/\mathbb{Z}  \& \ar[l, "{\begin{psmallmatrix}
     1 & 0\\
      \end{psmallmatrix}}"] \mathbb{R}/\mathbb{Z} \bigoplus \mathbb{R}/3\mathbb{Z} \ar[r, "{\begin{psmallmatrix}
     0 &1 \\
     \end{psmallmatrix}}", swap]  \& \mathbb{R}/3\mathbb{Z}  \\ 
   \end{tikzcd}
\caption{The diagram associated to the isogeny $\phi$ from Example \ref{example_runningexamplewayhome}.}\label{figure_diagramrunningexamplewayhome}
 \end{figure}
\end{example}
Note that in Example \ref{example_runningexamplewayhome} we used the output data from Example 62 (\cite{arXiv:2410.13459}) 
to generate $D_\phi$, but ended up in a slightly different place. Remark \ref{remark_choicesforbarphi} clarifies this.
 
\begin{remark}\label{remark_choicesforbarphi}

Construction \ref{construction_determinesplitting} should produce a pptav $J^{pp}$ and a dilation $\bar{\phi}$ that extend
\begin{diagram}
\Omega^1_{\mathbb{T}E'} (\mathbb{Z}) \bigoplus \Omega^1_{\mathbb{T}E} (\mathbb{Z})  & \arrow[l, "q^\# ",swap]\Omega^1_{\mathbb{T}E'} (\mathbb{Z}) \bigoplus \Omega^1_{\mathbb{T}E} (\mathbb{Z})  \\
  \H1(\mathbb{T}E', \mathbb{Z}) \bigoplus \H1(\mathbb{T}E, \mathbb{Z}) \arrow[u, "m_d \circ \zeta_{\mathbb{T}E' \bigoplus \mathbb{T}E} "] \arrow[r,"q_\#",hook]  &  \pi^{-1}(G)  \arrow[u,"\zeta"]  \\
\end{diagram}
to a commutative Diagram as in Figure \ref{diagram_interactionofphiq}. 
The specific pp $\zeta^{pp}$ and the isomorphism $\bar{\phi}_\#$ (see Construction \ref{construction_determinesplitting}) are just a matter of choice. 
We argue that these are natural in anticipation of Step 2 of Plan \ref{Masterplan} (see also Definition \ref{definition_JacobianandAbelJacobimappaper2}):
Given an arbitrary $\bar{\phi}_\#\in \Aut(\pi^{-1}(G))$, Construction \ref{construction_determinesplitting} has to be adapted as follows. Set
\begin{align}
\zeta^{\tilde{pp}}( \bar{\phi}_\#(B_i)):= \omega_i \text{ and }
    [\omega, B]^{\tilde{pp}}:= [\bar{\phi}^\#(\omega), \bar{\phi}^{-1}_\#(B)]_G, 
\end{align}
where $ \omega \in \Omega^1_{\mathbb{T}E'} (\mathbb{Z}) \bigoplus \Omega^1_{\mathbb{T}E} (\mathbb{Z})$ and $B\in \pi^{-1}(G)$.
Then $J^{\tilde{pp}}$ and $J^{pp}$ differ by an isomorphism $f:=(id,\bar{\phi}_\#): J^{pp} \rightarrow J^{\tilde{pp}}$ of pptavs that respects the pp, i.e. $f$ satisfies $f^*\zeta^{\tilde{pp}}=\zeta^{pp}$.
\end{remark}
 \subsection{Step 2 of Plan \ref{Masterplan} }\label{subsection_step2}
 Specifically, we ask whether $\gls{Jpp}$ from Construction \ref{construction_determinesplitting} is the Jacobian of a curve of genus $2$. Widening our perspective, we recognize here the incarnation of a classical problem, \emph{the tropical Schottky problem}. It is useful to embed step $2$ in this more general framework and to rephrase it as follows:
\begin{itemize}
    \item Does $J^{pp}$ lie in the \emph{tropical Schottky locus}, i.e. in the image of the \emph{tropical Torelli map} $t_2^{tr}: M^{tr}_2 \rightarrow A^{tr}_2$? 
    \item If yes, can we compute the preimage(s) of $J^{pp}$ under $t_2^{tr}$ explicitly?
\end{itemize}
When $g\leqslant 3$ the tropical Torelli map is surjective, hence feasibility (the first question) will not be an issue. It is the second point that is of real interest, especially since $t_2^{tr}$ is, in contrast to the classical case, \emph{not} injective (see \cite{MR3752493}, p. 24,25).
 \subsubsection{A concrete approach.}\label{subsubssection_step2concrete}
 Here, we take a practical approach to step 2 and ignore (as far as possible) the structural framework it is embedded in. We also provide a computational tool for constructing $\Gamma$ and use it in Section \ref{section_modulispaceperspective} for gaining experience of how $\Gamma$ depends on the choice of splitting data (see Lemma \ref{lemma_bdrycharacterizationforfamilies} and Lemma \ref{lemma_bdrycharacterizationforfamilies_d-1}).

 \paragraph{\emph{Interim Setup.}} Recall from Subsection \ref{subsection_tropicaltorelli}: The points of $A^{tr}_2$ are in bijection with the points of the quotient $\tilde{S}_{\geqslant 0}^2/ GL(\mathbb{Z})$, where $\tilde{S}_{\geqslant 0}^2$ denotes the space of positive semidefinite matrices with rational nullspace and the action of $GL(\mathbb{Z})$ on $\tilde{S}_{\geqslant 0}^2$ is given by $GL(\mathbb{Z}) \times \tilde{S}_{\geqslant 0}^2 \ni (X,Q) \mapsto X \sbullet Q:= X^TQX$ (\cite{MR2968636}, Section 4.3).

 An essential ingredient in this context is a result by Selling. Paraphrased for our purposes it reads:
\begin{proposition}[Selling]\label{propostion_Selling}
    Any positive definite $2 \times 2$ matrix $Q$ is $GL(\mathbb{Z})$-equivalent to a matrix in
    \begin{align}
        \langle \begin{pmatrix}
          1 & 0\\ 0 & 0
      \end{pmatrix},\begin{pmatrix}
          0 &0\\ 0 & 1
      \end{pmatrix},\begin{pmatrix}
          1 &-1\\ -1 & 1
      \end{pmatrix}\rangle_{\mathbb{R}\geqslant 0},
    \end{align}
    which can be computed explicitly using \emph{Selling's Reduction Algorithm} (see \cite{thesisVallentin}, Section 2.3.3, and \cite{zbMATH02716590}). If $Q$ satisfies $q_{12}<0$, then transformations of type (1) $\begin{pmatrix} 1 & 0\\
    1 & 1
        
    \end{pmatrix}$ and (2) $\begin{pmatrix} 1 & 1\\
    0 & 1 
    \end{pmatrix}$ are already sufficient.
\end{proposition}
 \begin{proposition}\label{proposition_fibreofTorelliMap}
 Let $J^{pp}$ be the pptav from Construction \ref{construction_determinesplitting}. Then one of two cases can occur:

\begin{enumerate}
    \item $J^{pp}$ is the Jacobian of a unique tropical curve $C$ whose combinatorial type is the theta-graph.
    \item There exists a family of tropical curves $\mathcal{C}$ over $\mathbb{R}_{\geqslant 0}$ such that $J^{pp}$ is the Jacobian of each fiber.
\end{enumerate}
Using Algorithm \ref{algorithm_PreimageofTormap} we can identify which one does and give an explicit description of the tropical curve(s) for each specific case.
\end{proposition}

\begin{proof}
In order to access the methods developed in the context of the tropical Schottky problem (\cite{MR3752493}), we first have to reformulate our input ($J^{pp}$ with principal polarization $\gls{zetapp}$) in adequate language:  
  The equivalence class of $(J^{pp},\zeta^{pp})$ corresponds to $\bar{Q}^{pp}\in \tilde{S}_{\geqslant 0}^2/ GL(\mathbb{Z})$, where $Q^{pp}$ is a coordinate representation of the scalar product $[\zeta^{pp}(\cdot), \cdot]^{pp}$. We have to compute $Q^{pp}$ next:
   We know from the proof of Proposition \ref{proposition_constructedpairJppandphisatisfiesreq} that $[\zeta^{pp}(\cdot), \cdot]^{pp}$ agrees with $[\zeta (\cdot), \cdot]_G$. Thus, upon fixing $\mathbb{Z}$-basis
   \begin{itemize}
       \item $S':=(\begin{pmatrix}
          l_{\mathbb{T}E'} \\ 0
      \end{pmatrix},  \begin{pmatrix}
          \frac{k \cdot l_{\mathbb{T}E'}}{d} \\ \frac{l_{\mathbb{T}E}}{d}
      \end{pmatrix})$ of $\pi^{-1}(G)$,
       \item $S:=((\omega,0),(0,\tilde{\omega}))$ on $\Omega^1_{\mathbb{T}E'} \bigoplus \Omega^1_{\mathbb{T}E} $ with dual basis $S^*$,
   \end{itemize}
   we see that $Q^{pp}:=M_{(S',S')}([\zeta^{pp}(\cdot), \cdot]^{pp})$ factors as $_{S'}M_S(\zeta)^t\cdot  M_{(S,S')}([\zeta(\cdot), \cdot]_G)$.
Since the quotient pairing (\cite{arXiv:2410.13459}, Lemma 26) 
is defined by
    \begin{align}
       [\cdot,\cdot]_G:  \Omega^1_{\mathbb{T}E'} \bigoplus \Omega^1_{\mathbb{T}E} \times \pi^{-1}(G) \rightarrow \mathbb{R}, (\lambda,\lambda') \rightarrow j(\lambda')(\lambda),
    \end{align}
    where $j$ denotes the embedding $\pi^{-1}(G) \hookrightarrow \Hom(\Omega^1_{\mathbb{T}E'},\mathbb{R}) \bigoplus \Hom(\Omega^1_{\mathbb{T}E},\mathbb{R})$, we have:
    \begin{align}
       M_{(S,S')}([\zeta(\cdot), \cdot]_G)= _{S'}M_{S^*}(j)=\begin{pmatrix}
          l_{\mathbb{T}E'} & \frac{k \cdot l_{\mathbb{T}E'}}{d}\\ 0 & \frac{l_{\mathbb{T}E}}{d}
      \end{pmatrix}.
    \end{align}
     To obtain a coordinate representation for $\zeta$ (the polarization inducing $m_d \circ \zeta_{\mathbb{T}E' \bigoplus \mathbb{T}E} $) next, we use Algorithm \ref{algorithm_constructionofinducingpolarization}, just to note that most of the work has already been carried out in the proof of Lemma \ref{lemma_polonquotientJ}: With the notation therein, fix $T':=((B,0),(0,\tilde{B}))$ as basis for $ \H1(\mathbb{T}E', \mathbb{Z}) \bigoplus \H1(\mathbb{T}E, \mathbb{Z})$ and take the basis $S$ and $S'$ as above. Then $_{S'}M_S(\zeta)$ is given by the matrix product
    \begin{align}
        _{T'}M_S(m_d \circ \zeta_{\mathbb{T}E' \bigoplus \mathbb{T}E} ) _{S''}M_{S'}^{-1}(i)= \begin{pmatrix}
          d & 0\\ 0 & d
      \end{pmatrix} \begin{pmatrix}
          1 & -k\\ 0 & d
      \end{pmatrix}^{-1}=\begin{pmatrix}
          d & k\\ 0 & 1
      \end{pmatrix}
      \end{align}
      and $Q^{pp}$ finally by
       \begin{align}
       M_{(S',S')}([\zeta^{pp}(\cdot), \cdot]^{pp})=
      \begin{pmatrix}
          d\cdot l_{\mathbb{T}E'} & k \cdot l_{\mathbb{T}E'}\\  k \cdot l_{\mathbb{T}E'} & \frac{k^2l_{\mathbb{T}E'}+l_{\mathbb{T}E} }{d}
      \end{pmatrix}.
      \end{align}
   The idea is to recognize in $Q^{pp}$ the period matrix of genus $2$ graph $\Gamma$ and recover $\Gamma$ by just reading off edge lengths from its matrix entries. However, the non-injectivity of $t_2^{tr}$ requires a more systematic approach: In the spirit of \cite{MR3752493}, we therefore compute a representative $Q_\Gamma$ for each relevant combinatorial type $\Gamma$ and check whether $Q^{pp}$ is $GL(\mathbb{Z})$-equivalent to one of these. From Section \ref{section_bridge}, with the corresponding base choice, we get
    \begin{align}
        Q^1_\Gamma=\begin{pmatrix}
          l(e)+l(e_2) &l(e_2)\\ l(e_2) & l(e_1)+l(e_2)
      \end{pmatrix} \text{ and } Q^2_\Gamma=\begin{pmatrix}
          l(e_1) &0\\ 0 & l(e_2)
      \end{pmatrix}
    \end{align} 
    for the theta-graph and the dumbbell-graph with corresponding closed cones 
    \begin{align}
        \sigma^1_\Gamma=\langle \begin{pmatrix}
          1 & 0\\ 0 & 0
      \end{pmatrix},\begin{pmatrix}
          0 &0\\ 0 & 1
      \end{pmatrix},\begin{pmatrix}
          1 &1\\ 1 & 1
      \end{pmatrix}\rangle_{\mathbb{R}\geqslant 0} \text{ and }\sigma^2_\Gamma=\langle \begin{pmatrix}
          1 & 0\\ 0 & 0
      \end{pmatrix},\begin{pmatrix}
          0 &0\\ 0 & 1
      \end{pmatrix}\rangle_{\mathbb{R} \geqslant 0}.
    \end{align} 
    Since $\sigma^2_\Gamma$ is just a face of $\sigma^1_\Gamma$, it suffices to consider $\sigma^1_\Gamma$, which we will, for the sake of convenience, replace by
    \begin{align}
        \sigma:=\begin{pmatrix}1 &0\\ 0 & -1
      \end{pmatrix} \sbullet \sigma^1_\Gamma =\langle \begin{pmatrix}
          1 & 0\\ 0 & 0
      \end{pmatrix},\begin{pmatrix}
          0 &0\\ 0 & 1
      \end{pmatrix},\begin{pmatrix}
          1 &-1\\ -1 & 1
      \end{pmatrix}\rangle_{\mathbb{R}\geqslant 0},
    \end{align} 
    accordingly $Q^{pp}$ by $\begin{pmatrix}1 &0\\ 0 & -1 \end{pmatrix} \sbullet Q^{pp}$, but retain our old notation $\gls{Qpp}$ in this case.
    Next, we perform Selling's Reduction Algorithm (Proposition \ref{propostion_Selling}) to reduce $Q^{pp}$ to a form lying in $\sigma$. Non-uniqueness poses a problem for determining edge lengths in a consistent way (see Remark \ref{remark_consistentedgelengthdetermination} and \ref{remark_consistentlengthdetermin2}). We therefore extend Selling's Reduction Algorithm (as described in \cite{thesisVallentin}, Section 2.3.3) by applying transformations $X\in \Stab(\sigma)$ to get a unique representative $\tilde{Q}^{pp}$ in 
    \begin{align}
        F:=\langle \begin{pmatrix}
          0 & 0\\ 0 & 1
      \end{pmatrix},\begin{pmatrix}
          1 &0\\ 0 & 1
      \end{pmatrix},\begin{pmatrix}
          2 &-1\\ -1 & 2
      \end{pmatrix}\rangle_{\mathbb{R}\geqslant 0}.
    \end{align} 
    Since $F\subset \sigma$ we can write $\tilde{Q}^{pp}$ as linear combination of the extreme rays of $\sigma$:
    \begin{align}
       \tilde{Q}^{pp}= l_1 \begin{pmatrix}
          1 & 0\\ 0 & 0
      \end{pmatrix}+ l_2\begin{pmatrix}
          0 &0\\ 0 & 1
      \end{pmatrix}+ l_3\begin{pmatrix}
          1 &-1\\ -1 & 1
      \end{pmatrix}.
    \end{align}

    If $\tilde{Q}^{pp}$ lies in the interior of $\sigma$, we conclude that $Q^{pp}$ is the period matrix of the theta-graph whose metric is fully determined by the coefficients $l_i$ for $i=1,2,3$, i.e. we have 
    \begin{align}
        l(e)=l_1,l(e_1)=l_2,l(e_2)=l_3.
    \end{align} 
    Suppose, however, that $\tilde{Q}^{pp}$ lies in a face of $\sigma$, then it can only lie in the interior of a $2$-dimensional face since $\det(Q^{pp})\neq 0$. Thus $Q^{pp}$ is equivalent to a form in $\sigma^2_\Gamma$. In other words, it is the period matrix of a family of curves $\mathcal{C} \rightarrow \mathbb{R}_{\geqslant 0}$ with
    \begin{itemize}
        \item general fibre $C_t$ whose combinatorial type is the dumbbell-graph with $l(e_1),l(e_2)$ determined by the two non-zero coefficients and bridge edge of length $t$. 
        \item special fibre $C_0$ whose combinatorial type is the common specialization of the theta- and the dumbbell-graph obtained by contracting one edge.
    \end{itemize}
 
\end{proof}
An algorithmic perspective on Proposition \ref{proposition_fibreofTorelliMap}:
\begin{algorithm}\label{algorithm_PreimageofTormap}
Input: A pptav $J^{pp}$ as in Proposition \ref{proposition_fibreofTorelliMap}. \\
Output: A tropical curve $C$ with $\Jac(C)=J^{pp}$.\\
    \begin{enumerate}
    \item Perform Selling's reduction algorithm to obtain a representative $Q$ of $Q^{pp}$ in $\sigma$. (Output option: the Selling Parameters of $Q$ and a list of integers that records the number of transformations of type (1) $\begin{pmatrix} 1 & 0\\
    1 & 1
        
    \end{pmatrix}$ and (2) $\begin{pmatrix} 1 & 1\\
    0 & 1 
    \end{pmatrix}$ yielding $Q$ (in reverse order of their application, ending with type (1))).
    \item For all $X\in \Stab(\sigma)$ compute $X\sbullet Q$. If $X\sbullet Q \in F$, set $Q:=X\sbullet Q $.
    \item Let $L$ be a list. Write $Q$ as linear combination of the extremal rays of $\sigma$. Add the coefficients to $L$. 
    \item Return $L$.
    \end{enumerate}
    \end{algorithm}
    We implemented this algorithm in SINGULAR and provide instructions for use that are of interest to us here, as they clarify the question of length assignement in the case that $Q$ lives in the boundary. For more details on Algorithm \ref{algorithm_PreimageofTormap} see Section \ref{section_appendix}.
    \begin{remark}\label{remark_consistentlengthdetermin2}
       Length output interpretation for Algorithm \ref{algorithm_PreimageofTormap}:
        \begin{enumerate}
            \item If all entries of $L$ are non-zero, then $Q^{pp}$ is the period matrix of the theta-graph and the entries of $L$  correspond to the edge lengths $l(e),l(e_1),l(e_2)$ (in this order).
        \item Else, $Q^{pp}$ is the period matrix of a family of curves whose combinatorial is the dumbbell graph. Using the labeling from Figure \ref{figure_conventionforthetaanddumbbell} 
        we interpret the entries of $L$ as follows:
     \begin{itemize}
         \item If $l_3=0$ (i.e. $Q\in \sigma^2_\Gamma$), set $l(e_i)=l_i$.
         \item If $l_2=0$ (i.e. $Q \not \in \sigma^2_\Gamma$), consider $\begin{pmatrix}
          -1 &0\\ -1 & 1
      \end{pmatrix}\sbullet Q$ and set $l(e_1)=l_1$ and $l(e_2)=l_3$.
         \item If $l_1=0$ (i.e. $Q\not \in \sigma^2_\Gamma$), consider $\begin{pmatrix}
          1 &-1\\ 0 & -1
      \end{pmatrix}\sbullet Q$ and set $l(e_1)=l_3$ and $l(e_2)=l_2$.
     \end{itemize}
     \end{enumerate}
    \end{remark}

 \subsubsection{From the perspective of moduli spaces.}\label{subsubsection_step2modulispaceversion}

As we move away from a sole focus on the algorithmic part of the problem, we bring the structural framework back into the picture: From the perspective of moduli spaces (Section \ref{subsection_tropicaltorelli}) Proposition \ref{proposition_fibreofTorelliMap} asks for a special fibre, the preimage of $Q^{pp}$ under $t_2^{tr}$. Since $\{\bar{\sigma}_{D_1},\bar{\sigma}_{D_2},\bar{\sigma}_{D_3},\bar{\sigma}_{D_4}\}$ is a complete set of representatives for $Gl(\mathbb{Z})$-equivalence classes of secondary cones (see \cite{MR2968636}, Example 4.10),
we (only) know that $Q^{pp}\in X\cdot \bar{\sigma}_{D_1}$ for an $X\in Gl(\mathbb{Z})$. We first find \emph{an} $X$ with $X^tQ^{pp}X\in \bar{\sigma}_{D_1}$ using Selling's Reduction Algorithm. The cone $\bar{\sigma}_{D_1}$, however, still has symmetries that arise from the action of its stabilzer $\Stab({\sigma}_{D_1})$. So we pick a fundamental domain, the cone $F$ in the proof of Proposition \ref{proposition_fibreofTorelliMap} (see blue-shaded area in Figure \ref{figure_Visualizationofreconstructionofgamma}), and compute the \emph{unique} representative $\tilde{Q}^{pp}$ of $Q^{pp}$ in $F$. If $l_i\neq 0$ holds for $i=1,2,3$ (see Proposition \ref{proposition_fibreofTorelliMap} and its proof), in other words, if $Q^{pp}$ is a \emph{generic} point (see \cite{MR2739784}, Definition 2.1.1) of $A_2^{tr}$, we have
 \begin{align}
   (t_2^{tr})^ {-1} (Q^{pp})=\{C\}.
 \end{align}
 This reflects the fact that $t_2^{tr}$ is of tropical degree $1$. Else
 \begin{align}
   (t_2^{tr})^ {-1} (Q^{pp})=\{ \mathcal{C} \},
 \end{align}
where $\mathcal{C}$ is the family of tropical curves from Proposition \ref{proposition_fibreofTorelliMap}.
Figure \ref{figure_Visualizationofreconstructionofgamma} offers a visualization of both situations.
 \begin{figure}[H]
     \centering
 \tikzset{every picture/.style={line width=0.75pt}} 

\begin{tikzpicture}[x=0.75pt,y=0.75pt,yscale=-1,xscale=1]

\draw [color={rgb, 255:red, 0; green, 0; blue, 0 }  ,draw opacity=1 ][fill={rgb, 255:red, 245; green, 166; blue, 35 }  ,fill opacity=1 ]   (332.95,46.88) -- (351.92,207.94) ;
\draw [color={rgb, 255:red, 0; green, 0; blue, 0 }  ,draw opacity=1 ][fill={rgb, 255:red, 245; green, 166; blue, 35 }  ,fill opacity=1 ]   (386.04,112.1) -- (332.95,46.88) ;
\draw [color={rgb, 255:red, 0; green, 0; blue, 0 }  ,draw opacity=1 ]   (267.51,112.13) -- (351.92,207.94) ;
\draw    (332.95,46.88) -- (300.23,79.51) -- (267.51,112.13) ;
\draw  [dash pattern={on 4.5pt off 4.5pt}]  (267.51,112.13) -- (359.5,79.49) ;
\draw  [dash pattern={on 4.5pt off 4.5pt}]  (386.04,112.1) -- (300.23,79.51) ;
\draw  [dash pattern={on 4.5pt off 4.5pt}]  (332.95,46.88) -- (326.77,112.12) ;
\draw    (386.04,112.1) -- (351.92,207.94) ;
\draw [color={rgb, 255:red, 0; green, 0; blue, 0 }  ,draw opacity=1 ][fill={rgb, 255:red, 245; green, 166; blue, 35 }  ,fill opacity=1 ]   (119.73,53.31) -- (138.7,214.37) ;
\draw [color={rgb, 255:red, 0; green, 0; blue, 0 }  ,draw opacity=1 ][fill={rgb, 255:red, 245; green, 166; blue, 35 }  ,fill opacity=1 ] [dash pattern={on 0.84pt off 2.51pt}]  (177.58,89.41) -- (119.73,53.31) ;
\draw [color={rgb, 255:red, 0; green, 0; blue, 0 }  ,draw opacity=1 ]   (76.33,89.41) -- (138.7,214.37) ;
\draw  [dash pattern={on 4.5pt off 4.5pt}]  (76.33,89.41) -- (220.97,53.31) ;
\draw  [dash pattern={on 4.5pt off 4.5pt}]  (119.73,53.31) -- (126.96,89.41) ;
\draw    (177.58,89.41) -- (138.7,214.37) ;
\draw    (138.7,214.37) -- (220.97,53.31) ;
\draw [color={rgb, 255:red, 0; green, 0; blue, 0 }  ,draw opacity=1 ][fill={rgb, 255:red, 74; green, 144; blue, 226 }  ,fill opacity=1 ]   (326.77,112.12) -- (351.92,207.94) -- (386.04,112.1) -- (327.73,88.83) ;
\draw [color={rgb, 255:red, 0; green, 0; blue, 0 }  ,draw opacity=1 ]   (386.04,112.1) -- (326.77,112.12) -- (267.51,112.13) ;
\draw   (119.73,53.31) -- (220.97,53.31) -- (177.58,89.41) -- (76.33,89.41) -- cycle ;
\draw    (76.33,89.41) -- (126.96,89.41) -- (177.58,89.41) ;
\draw    (76.33,89.41) -- (98.03,71.36) -- (119.73,53.31) ;
\draw  [dash pattern={on 4.5pt off 4.5pt}]  (98.03,71.36) -- (177.58,89.41) ;
\draw   (348.62,165.19) -- (355.58,172.31)(355.58,164.96) -- (348.62,172.54) ;
\draw   (120.6,161.51) -- (127.56,168.63)(127.56,161.28) -- (120.6,168.86) ;
\draw    (105.8,179.77) .. controls (4.48,269.34) and (455.47,276.15) .. (382.82,185.73) ;
\draw [shift={(381.69,184.36)}, rotate = 49.55] [color={rgb, 255:red, 0; green, 0; blue, 0 }  ][line width=0.75]    (10.93,-3.29) .. controls (6.95,-1.4) and (3.31,-0.3) .. (0,0) .. controls (3.31,0.3) and (6.95,1.4) .. (10.93,3.29)   ;
\draw [color={rgb, 255:red, 208; green, 2; blue, 27 }  ,draw opacity=1 ]   (190.22,87.91) .. controls (228.01,42.7) and (272.29,70.2) .. (286.29,83.74) ;
\draw [shift={(287.7,85.15)}, rotate = 226.55] [color={rgb, 255:red, 208; green, 2; blue, 27 }  ,draw opacity=1 ][line width=0.75]    (10.93,-3.29) .. controls (6.95,-1.4) and (3.31,-0.3) .. (0,0) .. controls (3.31,0.3) and (6.95,1.4) .. (10.93,3.29)   ;
\draw  [color={rgb, 255:red, 208; green, 2; blue, 27 }  ,draw opacity=1 ] (284.95,88.74) -- (293.92,92.59)(291.32,85.78) -- (287.55,95.55) ;
\draw  [color={rgb, 255:red, 208; green, 2; blue, 27 }  ,draw opacity=1 ] (139.75,131.2) -- (146.7,138.32)(146.7,130.97) -- (139.75,138.55) ;
\draw [color={rgb, 255:red, 208; green, 2; blue, 27 }  ,draw opacity=1 ]   (195.36,55.91) .. controls (206.67,67.85) and (174.47,89.89) .. (181.43,104.59) .. controls (188.4,119.29) and (172.73,133.07) .. (143.14,133.99) ;
\draw  [dash pattern={on 0.84pt off 2.51pt}]  (128.14,116.08) -- (158.14,151.89) ;

\draw (359.6,197.22) node [anchor=north west][inner sep=0.75pt]    {$A_{2}^{tr}$};
\draw (100.05,194.46) node [anchor=north west][inner sep=0.75pt]    {$M_{2}^{tr}$};
\draw (227.06,223.1) node [anchor=north west][inner sep=0.75pt]    {$t_{2}^{tr}$};
\draw (231.41,39.38) node [anchor=north west][inner sep=0.75pt]    {$t_{2}^{tr}$};

\end{tikzpicture}

     \caption{Propostion \ref{proposition_fibreofTorelliMap} for case (1) in black and case (2) in red.}
     \label{figure_Visualizationofreconstructionofgamma}
 \end{figure}
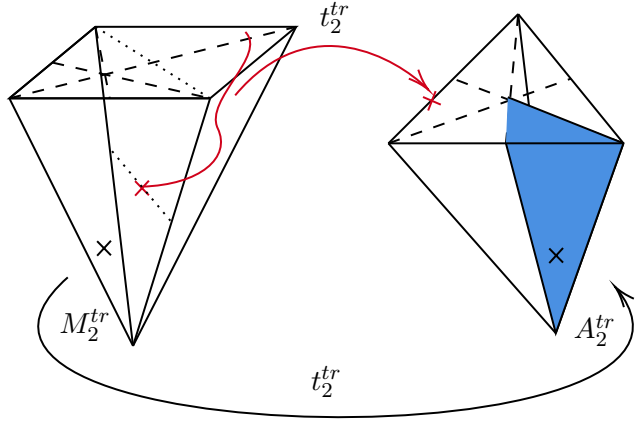
 \begin{remark/reference}\label{remark_consistentedgelengthdetermination}
 
    The reconstruction is based on the approach of \cite{MR3752493}, where they construct an \emph{arbitrary} curve with a given tropical Jacobian. Not only arbitrary in the sense, that a positive definite quadratic form $Q$ can be the period matrix associated to curves with different combinatorial types (the natural non-injectivity of $t^{tr}_2$). But also arbitrary when specifying the metric for a fixed combinatorial type $(G,g)$. There is not necessarily a unique $Gl(\mathbb{Z})$ representative of $Q$ in the corresponding cone $\sigma_{(G,g)}$, since we can still operate on $Q$ with $X\in \Stab(\sigma_{(G,g)})$. Depending on our choice of representative, we get different metrics on $(G,g)$. The reason for this is as follows: The action of $\Stab(\sigma_{(G,g)})$ on elements of $\sigma_{(G,g)}$ really corresponds to permuting edges of $G$. To determine edge lengths, however, we fix a cycle base and an order on $E(G)$, write $Q$ as linear combination of the corresponding extreme rays and match the resulting coefficients to the edges according to our chosen order.
    In order to do justice to the geometric meaning behind the $\Stab(\sigma_{(G,g)})$-action, we remove the element of arbitrariness and only use the representative of $Q$ in $F$ to determine the metric.
\end{remark/reference}
 \subsection{Step 3 of Plan \ref{Masterplan}}\label{subsection_step3} 
 Fix $\Gamma:=C$ or $\Gamma:=C_0, \Gamma_t:=C_t$ for $t\in \mathbb{R}_{> 0}$, $P_0\in \Gamma$ as in Figure \ref{figure_conventionforthetaanddumbbell}, depending on whether case (1) or (2) occurs in Proposition \ref{proposition_fibreofTorelliMap}, and define maps
\begin{align}
    \varphi:= g_2 \circ \Phi_{P_0} \text{ and } \varphi':= g_1 \circ \Phi_{P_0}
\end{align}
adding a subscript $t$ $(\varphi_t,\varphi'_t)$ for case $(2)$. At this point, all that remains is to assemble the building blocks from steps 1 and 2 of Plan \ref{Masterplan}.
\begin{theorem}\label{theorem_reconstructionsumup}
    Let $(\mathbb{T}E',\mathbb{T}E,G)$ be splitting data, i.e. $G$ is the graph of an isomorphism $\alpha$ between the $d$-torsion points of $\mathbb{T}E'$ and $\mathbb{T}E$. 
    Then there exists
    \begin{itemize}
        \item a (family of) curve(s) of genus $2$ $\Gamma$ ($\Gamma_t$ for $t\geqslant 0$),
        \item a (family of) pair(s) of optimal covers $(\varphi': \Gamma \rightarrow \mathbb{T}E',\varphi: \Gamma \rightarrow \mathbb{T}E)$ ($(\varphi_t,\varphi'_t)$ for $t\geqslant 0$) of degree $d$,
    \end{itemize}
    that induce a splitting
    \begin{align}
      \phi:\mathbb{T}E' \bigoplus \mathbb{T}E \rightarrow \Jac(\Gamma)
  \end{align}
  with  $\Jac_d(\mathbb{T}E')\cong \ker(\phi) \cong \Jac_d(\mathbb{T}E)$.
\end{theorem}
\begin{proof}
   Let $\phi, \Gamma,\varphi'$ and $\varphi$ be as in steps 1-3, omitting the subscript $t$ for clarity. All that remains is to show that the maps $\varphi$ and $\varphi'$ are optimal covers of degree $d$ with $g_1=\varphi'_*$, $g_2=\varphi_*$, $f_1=\varphi^*$ and $f_2=\varphi^{'*}$. 
    Compare the definition of $\varphi$ and $\varphi'$ with the complementary cover (\cite{arXiv:2410.13459}, Section 6.2) and note that they are constructed according to the same pattern. We therefore refrain from proving the covering-property and refer to Lemma 64 in \cite{arXiv:2410.13459}, instead. 
    With two covers at our disposal we embark/set out to reinterpret Diagram \ref{diagram_inducedbyisogeny} in terms of $\varphi$ and $\varphi'$: We specify embeddings
\begin{align}
    & j_1: \mathbb{T}E' \rightarrow \Jac(\mathbb{T}E'), P' \mapsto P' - 0_{\Jac(\mathbb{T}E')}\\
     & j_2: \mathbb{T}E \rightarrow \Jac(\mathbb{T}E), P \mapsto P - 0_{{\Jac(\mathbb{T}E)}}
\end{align}
that are compatible with $\Phi_{P_0}$. Under these identifications (i.e. $\mathbb{T}E' \cong \Jac(\mathbb{T}E')$ and $\mathbb{T}E' \cong \Jac(\mathbb{T}E')$ via $j_1$ and $j_2$) we see that
\begin{itemize}
    \item $g_1$ ($g_2$) corresponds to the push-forward $\varphi'_*$ ($\varphi_*$) since
     \begin{diagram}
 \Gamma \arrow[r, "{\Phi}_{P_0} "] \arrow[d, "\varphi' "] 
&  \Jac(\Gamma) \arrow[d, "g_1"] & \Gamma \arrow[r, "{\Phi}_{P_0} "] \arrow[d, "\varphi "] 
&  \Jac(\Gamma) \arrow[d, "g_2"] \\
\mathbb{T}E' \arrow[r, hook ,"j_1"] & \Jac(\mathbb{T}E') & \mathbb{T}E \arrow[r, hook ,"j_2"] & \Jac(\mathbb{T}E),   \\
  \end{diagram}
  commute.
    \item $f_1$ ($f_2$) corresponds to the pull-back $\varphi'^*$ ($\varphi^*$) (\cite{arXiv:2410.13459}, Lemma 35) 
\end{itemize}
which, (when combined/together) with the exactness of the diagonals in \ref{diagram_inducedbyisogeny}, implies optimality of $\varphi$ and $\varphi'$. We can now compute $deg(\varphi)$ and $deg(\varphi')$ as in the proof of Theorem 66 (\cite{arXiv:2410.13459}) 
For $\varphi$ we have
\begin{align}
   deg(\varphi)id & = \varphi_* \circ \varphi^*= (g_2 \circ j_2) \circ (j^{-1}_2 \circ f_2)\\
   & = p_2 \circ \tilde{\phi} \circ \phi \circ \iota_2, 
    =  p_2 \circ (m_d,m_d) \circ  \iota_2= m_d
\end{align}
and thus $deg(\varphi)=d$. The analogous computation for $\varphi'$ shows $deg(\varphi')=d$ as well. 
\end{proof}
We let Example \ref{example_runningexamplewayhome} rest for reconstruction step 2, as it was already conceivable from the "Cover-to-Splitting" direction (see Figure \ref{figure_optimalandnotoptimalcover} and Example 67 (\cite{arXiv:2410.13459}), 
that $J^{pp}$ is the Jacobian of the theta-graph whose edges have all length $1$. We pick it up again to illustrate step 3 of Plan \ref{Masterplan}.
\begin{example}
We define 
  \begin{diagram}[ampersand replacement = \&]
\varphi': \Gamma \arrow[r, hook, "\Phi_{P_0} "] \& J^{pp}=\Jac(\Gamma) \arrow[r,"g_1"] \arrow[d, "\cong "] \& \Jac(\mathbb{T}E') \arrow[r,"j_1"] \arrow[d, "\cong "]  \& \mathbb{T}E' \\
\& \mathbb{R}^2/ \begin{pmatrix}
     2 & 1\\
    1 & 2
 \end{pmatrix}\mathbb{Z}^2  \arrow[r,"{\begin{psmallmatrix}
    1 & 0
 \end{psmallmatrix}}"] \&  \mathbb{R}/\mathbb{Z} 
 \end{diagram}
 and
\begin{diagram}[ampersand replacement = \&]
 \varphi:\Gamma \arrow[r, hook, "\Phi_{P_0} "] \& J^{pp}=\Jac(\Gamma) \arrow[r,"g_2"] \arrow[d, "\cong "] \& \Jac(\mathbb{T}E) \arrow[r,"j_2"] \arrow[d, "\cong "]  \& \mathbb{T}E \\
\& \mathbb{R}^2/ \begin{pmatrix} 2 & 1\\    1 & 2 \end{pmatrix}\mathbb{Z}^2  \arrow[r,"{\begin{psmallmatrix}     -1 & 2
 \end{psmallmatrix}}"] \&  \mathbb{R}/3\mathbb{Z} 
  \end{diagram}
  where coordinate representations are taken from Figure \ref{figure_diagramrunningexamplewayhome}. In order to understand $\varphi$ and $\varphi'$ as morphisms of graphs we work locally as \cite{arXiv:2410.13459} (see Example 67 for details). 
  We can easily read off parametrizations of $e,e_1$ and $e_2$ in $J^{pp}$ from Figure 10 (\cite{arXiv:2410.13459}). 
  By acting on these with $M(g_i)$ we finally obtain
 \begin{align}
     & \varphi'_{|e} (t) = t + 0_{\Jac(\mathbb{T}E')} , \varphi'_{|e_1}(t) =  0_{\Jac(\mathbb{T}E')}, \varphi'_{|e_2}(t) =  t + 0_{\Jac(\mathbb{T}E')}, \\
     & \varphi_{|e}(t) =  -t + 0_{\Jac(\mathbb{T}E)} , \varphi_{|e_1}(t) =  -2t + 0_{\Jac(\mathbb{T}E)}, \varphi_{|e_2}(t) = t + 0_{\Jac(\mathbb{T}E)},
 \end{align}
 where $t\in [0,1]$.
 \begin{figure}[H]
     \centering
   \begin{tikzpicture}[x=0.75pt,y=0.75pt,yscale=-1,xscale=1]

\draw  [color={rgb, 255:red, 0; green, 0; blue, 0 }  ,draw opacity=1 ] (57.44,240.74) .. controls (57.46,232.59) and (88.86,226) .. (127.57,226.02) .. controls (166.28,226.03) and (197.64,232.64) .. (197.62,240.79) .. controls (197.6,248.93) and (166.2,255.53) .. (127.49,255.51) .. controls (88.78,255.5) and (57.41,248.88) .. (57.44,240.74) -- cycle ;
\draw [color={rgb, 255:red, 0; green, 0; blue, 0 }  ,draw opacity=1 ]   (57.44,240.74) ;
\draw [shift={(57.44,240.74)}, rotate = 0] [color={rgb, 255:red, 0; green, 0; blue, 0 }  ,draw opacity=1 ][fill={rgb, 255:red, 0; green, 0; blue, 0 }  ,fill opacity=1 ][line width=0.75]      (0, 0) circle [x radius= 3.35, y radius= 3.35]   ;
\draw  [color={rgb, 255:red, 0; green, 0; blue, 0 }  ,draw opacity=1 ] (119.05,251.9) .. controls (126.39,253.91) and (133.73,255.11) .. (141.07,255.52) .. controls (133.73,255.92) and (126.39,257.13) .. (119.05,259.14) ;
\draw   (111.19,85.55) .. controls (117.73,85.05) and (124.16,85.21) .. (130.47,86.01) .. controls (124.28,84.57) and (118.21,82.49) .. (112.25,79.76) ;
\draw    (54.65,50.64) -- (54.65,120.37) ;
\draw    (54.65,50.64) .. controls (81.44,30.64) and (203.76,50.94) .. (197.51,87.71) ;
\draw    (54.65,120.37) .. controls (72.51,149.49) and (189.48,130.37) .. (197.51,87.71) ;
\draw    (54.65,50.64) .. controls (52.87,88.59) and (212.69,91.24) .. (170.73,93.01) ;
\draw    (54.65,120.37) .. controls (68.94,93.89) and (151.98,93.01) .. (170.73,93.01) ;
\draw   (112.69,130.37) .. controls (119.57,131.14) and (126.3,131.14) .. (132.89,130.34) .. controls (126.45,131.93) and (120.16,134.3) .. (114.02,137.46) ;
\draw    (123.08,154.26) -- (123.65,197.27) ;
\draw [shift={(123.68,199.27)}, rotate = 269.24] [color={rgb, 255:red, 0; green, 0; blue, 0 }  ][line width=0.75]    (10.93,-3.29) .. controls (6.95,-1.4) and (3.31,-0.3) .. (0,0) .. controls (3.31,0.3) and (6.95,1.4) .. (10.93,3.29)   ;
\draw    (54.65,120.37) ;
\draw [shift={(54.65,120.37)}, rotate = 0] [color={rgb, 255:red, 0; green, 0; blue, 0 }  ][fill={rgb, 255:red, 0; green, 0; blue, 0 }  ][line width=0.75]      (0, 0) circle [x radius= 3.35, y radius= 3.35]   ;
\draw    (54.65,50.64) ;
\draw [shift={(54.65,50.64)}, rotate = 0] [color={rgb, 255:red, 0; green, 0; blue, 0 }  ][fill={rgb, 255:red, 0; green, 0; blue, 0 }  ][line width=0.75]      (0, 0) circle [x radius= 3.35, y radius= 3.35]   ;
\draw   (51.08,94.18) .. controls (53.04,88.85) and (54.23,83.52) .. (54.66,78.18) .. controls (55,83.53) and (56.1,88.88) .. (57.96,94.24) ;
\draw  [draw opacity=0] (276.87,108.92) .. controls (263.43,105.1) and (255.35,100.22) .. (255.36,94.9) .. controls (255.38,82.75) and (297.59,72.96) .. (349.65,73.02) .. controls (401.71,73.08) and (443.9,82.98) .. (443.88,95.13) .. controls (443.87,100.41) and (435.91,105.24) .. (422.64,109.01) -- (349.62,95.02) -- cycle ; \draw   (276.87,108.92) .. controls (263.43,105.1) and (255.35,100.22) .. (255.36,94.9) .. controls (255.38,82.75) and (297.59,72.96) .. (349.65,73.02) .. controls (401.71,73.08) and (443.9,82.98) .. (443.88,95.13) .. controls (443.87,100.41) and (435.91,105.24) .. (422.64,109.01) ;  
\draw   (276.02,110.2) .. controls (276.04,101.8) and (308.99,95) .. (349.62,95.02) .. controls (390.25,95.03) and (423.17,101.85) .. (423.14,110.25) .. controls (423.12,118.65) and (390.16,125.45) .. (349.54,125.43) .. controls (308.91,125.42) and (275.99,118.6) .. (276.02,110.2) -- cycle ;
\draw  [color={rgb, 255:red, 0; green, 0; blue, 0 }  ,draw opacity=1 ] (273.44,241.74) .. controls (273.46,233.59) and (304.86,227) .. (343.57,227.02) .. controls (382.28,227.03) and (413.64,233.64) .. (413.62,241.79) .. controls (413.6,249.93) and (382.2,256.53) .. (343.49,256.51) .. controls (304.78,256.5) and (273.41,249.88) .. (273.44,241.74) -- cycle ;
\draw [color={rgb, 255:red, 0; green, 0; blue, 0 }  ,draw opacity=1 ]   (273.44,241.74) ;
\draw [shift={(273.44,241.74)}, rotate = 0] [color={rgb, 255:red, 0; green, 0; blue, 0 }  ,draw opacity=1 ][fill={rgb, 255:red, 0; green, 0; blue, 0 }  ,fill opacity=1 ][line width=0.75]      (0, 0) circle [x radius= 3.35, y radius= 3.35]   ;
\draw  [color={rgb, 255:red, 0; green, 0; blue, 0 }  ,draw opacity=1 ] (335.05,252.9) .. controls (342.39,254.91) and (349.73,256.11) .. (357.07,256.52) .. controls (349.73,256.92) and (342.39,258.13) .. (335.05,260.14) ;
\draw    (342.08,156.26) -- (342.65,199.27) ;
\draw [shift={(342.68,201.27)}, rotate = 269.24] [color={rgb, 255:red, 0; green, 0; blue, 0 }  ][line width=0.75]    (10.93,-3.29) .. controls (6.95,-1.4) and (3.31,-0.3) .. (0,0) .. controls (3.31,0.3) and (6.95,1.4) .. (10.93,3.29)   ;
\draw    (276.87,108.92) -- (276.02,110.2) ;
\draw [shift={(276.44,109.56)}, rotate = 123.78] [color={rgb, 255:red, 0; green, 0; blue, 0 }  ][fill={rgb, 255:red, 0; green, 0; blue, 0 }  ][line width=0.75]      (0, 0) circle [x radius= 3.35, y radius= 3.35]   ;
\draw    (422.64,109.01) -- (424.62,108.43) ;
\draw [shift={(423.63,108.72)}, rotate = 343.6] [color={rgb, 255:red, 0; green, 0; blue, 0 }  ][fill={rgb, 255:red, 0; green, 0; blue, 0 }  ][line width=0.75]      (0, 0) circle [x radius= 3.35, y radius= 3.35]   ;
\draw [color={rgb, 255:red, 0; green, 0; blue, 0 }  ,draw opacity=1 ]   (413.62,241.79) ;
\draw [shift={(413.62,241.79)}, rotate = 0] [color={rgb, 255:red, 0; green, 0; blue, 0 }  ,draw opacity=1 ][fill={rgb, 255:red, 0; green, 0; blue, 0 }  ,fill opacity=1 ][line width=0.75]      (0, 0) circle [x radius= 3.35, y radius= 3.35]   ;
\draw   (355.1,128.64) .. controls (348.45,126.7) and (341.82,125.55) .. (335.19,125.22) .. controls (341.81,124.75) and (348.41,123.49) .. (355,121.43) ;
\draw   (335.33,69.12) .. controls (341.87,71.39) and (348.44,72.86) .. (355.04,73.53) .. controls (348.41,73.66) and (341.75,74.59) .. (335.07,76.33) ;
\draw   (336.17,91.53) .. controls (342.8,93.52) and (349.43,94.73) .. (356.05,95.12) .. controls (349.43,95.53) and (342.82,96.73) .. (336.21,98.74) ;

\draw (207.38,237.72) node [anchor=north west][inner sep=0.75pt]  [color={rgb, 255:red, 0; green, 0; blue, 0 }  ,opacity=1 ]  {$\mathbb{T} E'$};
\draw (127.43,27.81) node [anchor=north west][inner sep=0.75pt]    {$e$};
\draw (32.44,71.57) node [anchor=north west][inner sep=0.75pt]    {$e_{1}$};
\draw (78.55,56.34) node [anchor=north west][inner sep=0.75pt]    {$e_{2}$};
\draw (22.98,29.1) node [anchor=north west][inner sep=0.75pt]    {$P_{0}$};
\draw (25.66,115.59) node [anchor=north west][inner sep=0.75pt]    {$P_{1}$};
\draw (26.17,249.33) node [anchor=north west][inner sep=0.75pt]  [color={rgb, 255:red, 0; green, 0; blue, 0 }  ,opacity=1 ]  {$0_{\Jac(\mathbb{T} E')}$};
\draw (143.25,163.43) node [anchor=north west][inner sep=0.75pt]    {$\varphi '$};
\draw (213.63,55.3) node [anchor=north west][inner sep=0.75pt]    {$\Gamma $};
\draw (423.38,238.72) node [anchor=north west][inner sep=0.75pt]  [color={rgb, 255:red, 0; green, 0; blue, 0 }  ,opacity=1 ]  {$\mathbb{T} E$};
\draw (242.17,250.33) node [anchor=north west][inner sep=0.75pt]  [color={rgb, 255:red, 0; green, 0; blue, 0 }  ,opacity=1 ]  {$0_{\Jac(\mathbb{T} E)}$};
\draw (362.25,165.43) node [anchor=north west][inner sep=0.75pt]    {$\varphi $};
\draw (254.98,114.73) node [anchor=north west][inner sep=0.75pt]    {$P_{0}$};
\draw (428.66,118.59) node [anchor=north west][inner sep=0.75pt]    {$P_{1}$};
\draw (317.44,49.57) node [anchor=north west][inner sep=0.75pt]    {$e_{1}$};
\draw (296.55,81.34) node [anchor=north west][inner sep=0.75pt]    {$e_{2}$};
\draw (314.43,125.73) node [anchor=north west][inner sep=0.75pt]    {$e$};
\draw (271,61.73) node [anchor=north west][inner sep=0.75pt]    {$2$};
\draw (450.63,55.3) node [anchor=north west][inner sep=0.75pt]    {$\Gamma $};

\end{tikzpicture}

     \caption{$\varphi'$ and  $\varphi$ as morphisms of graphs.}
     \label{fig:enter-label}
 \end{figure}
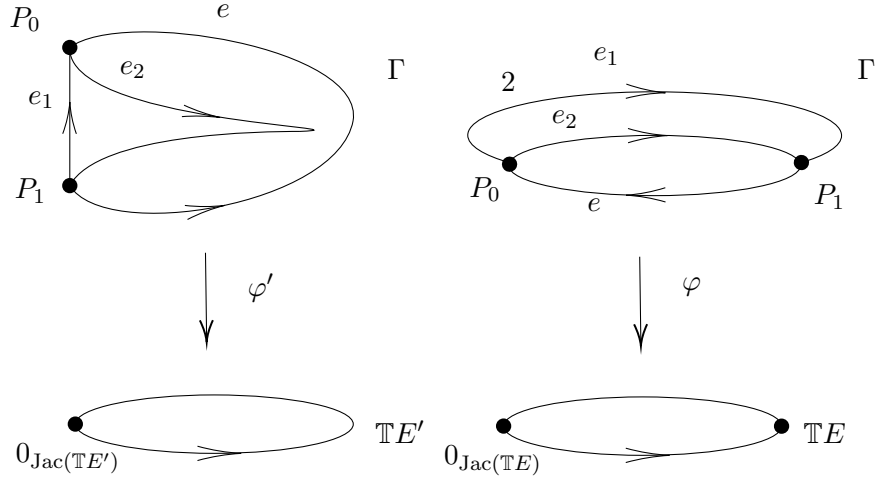
\end{example}
 \begin{remark}
     We can now answer the two questions at the beginning Section \ref{section_reconstruction}: $(\mathbb{T}E',\mathbb{T}E,G)$ is splitting data if and only if it satisfies the conditions described in Lemma \ref{lemma_reqfordiagrwithexactdiag}. Moreover, Theorem \ref{theorem_reconstructionsumup} shows, it is both necessary and sufficient.
 \end{remark}
 \section{Moduli space perspective}\label{section_modulispaceperspective}
 We investigate a Schottky-type problem and study the fibres of $t_2^{tr}$ for the case of split Jacobians. In this context, we are interested in the following subsets:
\begin{enumerate}
    \item $\mathcal{Q}\subset A_2^{tr}$ the \emph{locus of split Jacobians}.
    \item $\mathbb{T}\mathcal{L}_d\subset M_2^{tr}$ the \emph{locus of curves with d-split Jacobians}. 
\end{enumerate}
The second is a tropical analogue of the locus of genus 2 curves with $(d,d)$-split Jacobian studied for example in \cite{shaska2024machinelearningmodulispace}, \cite{zbMATH05564780}) or \cite{zbMATH02074360}. The notation from the previous Section is maintained.
\subsection{A Schottky-type Problem}\label{subsection_Schottky-typeProblem}
As a result of our reconstruction procedure in Section \ref{section_reconstruction}, we have a representative for each isomorphism-class of pptavs that gives rise to a split Jacobian, i.e. we have a description for $\mathcal{Q}$:
\begin{align}
    \mathcal{Q}=\{Q^{pp}(d,k,l_{\mathbb{T}E'},l_{\mathbb{T}E}) \thinspace | \thinspace d\in \mathbb{N}_{>1}, k\in\mathbb{Z}^*_d, l_{\mathbb{T}E'},l_{\mathbb{T}E}\in \mathbb{R}_{>0}  \}.
\end{align}
With this family of pptavs coming from geometry at hand, it is natural to ask for $\mathcal{Q} \cap \partial A_2^{tr} $. For each specific case, Proposition \ref{proposition_fibreofTorelliMap} tells us whether $\gls{Qpp}$ is a generic point, but it does so only \emph{a posteriori} meaning after running through Algorithm \ref{algorithm_PreimageofTormap}. Much more satisfying would be an \emph{a priori} characterization, one that is inspired by the following observation
    \begin{align}
        Q^{pp}\in \sigma \text{ is not a generic point } \Leftrightarrow \thinspace & l_{\mathbb{T}E}=kl_{\mathbb{T}E'}(d-k)\text{ for } k,d\in \mathbb{N} ,\\
        & k<d, \gcd(k,d)=1. 
    \end{align}
    This can be seen immediately by writing $Q^{pp}$ as linear combination of the extreme rays of $\sigma$ and noting that only one of the coefficients can vanish. More generally:
\begin{remark}\label{remark_QliesinbdrywhenSellingParametervanishes}
   We have used the description of $\sigma$ in terms of generators. As a system of inequalities it is given by
   \begin{align}
     \{ \begin{pmatrix}
         q_{11}\\
         q_{12} \\
         q_{22}
     \end{pmatrix} \in \mathbb{R}^3 \thinspace |\thinspace q_{12}\leqslant 0, q_{11}+q_{12} \geqslant 0, q_{22}+q_{12} \geqslant 0\},  
   \end{align}
   where we identify a $2 \times 2$ symmetric matrix with entries $q_{ij}$ with a vector in $\mathbb{R}^3$. Then, $Q\in \sigma$ lies in a $2$-dimensional face if and only if one of the inequalities is not strict.
\end{remark}
We exploit Remark \ref{remark_QliesinbdrywhenSellingParametervanishes} and Algorithm \ref{algorithm_PreimageofTormap} to give a concrete description of the intersection with the boundary for two subfamilies of $\mathcal{Q}$.
\begin{lemma}\label{lemma_bdrycharacterizationforfamilies}
Consider the same setting as in Proposition \ref{proposition_fibreofTorelliMap} and let $k=1$. Then case $(2)$ holds if and only if there exists $\alpha\in \mathbb{N}_{>0}$ with
\begin{align}
    \alpha\cdot l_{\mathbb{T}E}= (d-\alpha)\cdot l_{\mathbb{T}E'}.
\end{align}
\end{lemma}
\begin{proof}
    The proof is based on two simple observations:
    \begin{enumerate}
        \item A reduced representative $Q$ of $Q^{pp}$ inside $\sigma$ can be obtained only using the transformation $\begin{pmatrix}
            1 & 0 \\
            1 & 1 
        \end{pmatrix}$.
        \item $Q$ lies in a $2$-dimensional face if and only if $q_{22}+q_{12}=0$ holds.
    \end{enumerate}
    We start with $(1)$: Consider the sequence 
    \begin{align}
       a_n:=-q^{pp}_{12}-(n+1)q^{pp}_{22}, \thinspace  n\in \mathbb{N}
    \end{align}
    and note that $(a_n)_{n \in \mathbb{N}}$ is strictly decreasing. Let
    $\tilde{\alpha}:=\min \{ n \in \mathbb{N}| a_n \leqslant 0  \}$.
    We claim that $Q:=\begin{pmatrix}
            1 & 0 \\
            \tilde{\alpha} & 1 
        \end{pmatrix}\sbullet Q^{pp} \in \sigma$, i.e. that $Q$ satisfies the inequalities from Remark \ref{remark_QliesinbdrywhenSellingParametervanishes}. A computation shows
        \begin{align}
           & q_{12}= q^{pp}_{12}+\tilde{\alpha} q^{pp}_{22},\\
           & q_{11}+q_{12}=q^{pp}_{11}+(2\tilde{\alpha}+1)q^{pp}_{12}+\tilde{\alpha}(\tilde{\alpha}+1)q^{pp}_{22},\\
           & q_{22}+q_{12} =q^{pp}_{12}+(\tilde{\alpha}+1)q^{pp}_{22},
        \end{align}
        i.e.
    
        \begin{align}
           & q_{12}= -a_{\tilde{\alpha}-1} < 0\\
           & q_{11}+q_{12}= q^{pp}_{11}+(\tilde{\alpha}+1)q^{pp}_{12} - \tilde{\alpha}a_{\tilde{\alpha}}\geqslant q^{pp}_{11}+(\tilde{\alpha}+1)q^{pp}_{12}\\
           & q_{22}+q_{12} =-a_{\tilde{\alpha}} \geqslant 0
        \end{align}
        as an immediate consequence of the definition of $\tilde{\alpha}$ and the fact that $(a_n)_{n \in \mathbb{N}}$ is strictly decreasing.
        Recalling that $q^{pp}_{11}=dl_{\mathbb{T}E'}$ and $q^{pp}_{12}=-kl_{\mathbb{T}E'}$, we obtain
        \begin{align}
            q^{pp}_{11}+(\tilde{\alpha}+1)q^{pp}_{12} \geqslant 0
        \end{align}
        since $d \geqslant \tilde{\alpha} +1$ (again by minimality of $\tilde{\alpha}$) and 
        conclude (1). Note that point (2) is just a special case of Remark \ref{remark_QliesinbdrywhenSellingParametervanishes} when combined with the following observation:
        \begin{align}
            q_{11}+q_{12}=0 &\Leftrightarrow q^{pp}_{11}+(\tilde{\alpha}+1)q^{pp}_{12}=0  \text{ and }  \tilde{\alpha}a_{\tilde{\alpha}}=0\\
           &\Leftrightarrow d=\tilde{\alpha}+1  \text{ and }  \tilde{\alpha}(l_{\mathbb{T}E'} - (\tilde{\alpha} +1)\frac{l_{\mathbb{T}E'}+l_{\mathbb{T}E}}{d})=0\\
           &\Leftrightarrow d=\tilde{\alpha}+1  \text{ and }  (l_{\mathbb{T}E'} - (\tilde{\alpha} +1)\frac{l_{\mathbb{T}E'}+l_{\mathbb{T}E}}{d})=0,
           \end{align}
       where the last equivalence holds since $d>1$ (and thus $\tilde{\alpha}\neq 0)$, so we have a contradiction.\\
       As intermediate result (from $(1)$ and $(2)$) we get: Case $(2)$ in Proposition \ref{proposition_fibreofTorelliMap}) holds if and only if $q_{22}+q_{12}$, i.e. $a_{\tilde{\alpha}}=l_{\mathbb{T}E'} - (\tilde{\alpha} +1)\frac{l_{\mathbb{T}E'}+l_{\mathbb{T}E}}{d}$ vanishes. 
       This shows the "only if"-direction (by setting $\alpha:=\tilde{\alpha}+1$  of the statement of Lemma \ref{lemma_bdrycharacterizationforfamilies}. For the "if"-direction suppose there exists such an $\alpha\in \mathbb{N}_{>0}$, then $a_{\alpha -1}=0$ and $\alpha-1=\tilde{\alpha}$ since $(a_n)_{n\in \mathbb{N}}$ is strictly decreasing. This means $Q\in \sigma$ by point $(1)$, in particular $Q$ is in a 2-dimensional face by $(2)$.
\end{proof}
A symmetrical result is content of the following Lemma.
\begin{lemma}\label{lemma_bdrycharacterizationforfamilies_d-1}
Consider the same setting as in Proposition \ref{proposition_fibreofTorelliMap} and let $d\geqslant 3$, $k=d-1$. Then case $(2)$ holds if and only if there exists $\beta\in \mathbb{N}_{>0}$ with
\begin{align}
    \beta\cdot l_{\mathbb{T}E}= (d-\beta)\cdot l_{\mathbb{T}E'}.
\end{align}
\end{lemma}
The proof is based on the same methods, but their implementation is tedious and more complicated. We leave it at a rough sketch. 
    \begin{sketchproof} 
    If $Q^{pp}\in \sigma$, then case $(2)$ holds if and only if $l_{\mathbb{T}E}= (d-1)\cdot l_{\mathbb{T}E'}$ (see Observation above Remark \ref{remark_QliesinbdrywhenSellingParametervanishes}. 
    If not, replace the two observations from the proof of Lemma \ref{lemma_bdrycharacterizationforfamilies} by:
    \begin{enumerate}
        \item We can find a reduced representative $Q$ of $Q^{pp}$ inside $\sigma$ of the form 
        \begin{align}
            \begin{pmatrix}
            1 & \tilde{\beta} \\
            0 & 1 
        \end{pmatrix}\begin{pmatrix}
            1 & 0 \\
            1 & 1 
        \end{pmatrix}\sbullet Q^{pp}.
        \end{align}
        \item $Q$ lies in a $2$-dimensional face if and only if $q_{11}+q_{12}=0$ holds.
    \end{enumerate}
    Consider a sequence $(b_n)_{n\in \mathbb{N}}$ analogous to $(a_n)_{n\in \mathbb{N}}$ and set 
    \begin{align}
        \tilde{\beta}:=\min \{ n \in \mathbb{N}| b_n \leqslant 0  \}. 
    \end{align}
    Then check that the entries of $Q$ satisfy the inequalities from Remark \ref{remark_QliesinbdrywhenSellingParametervanishes} and argue that equality can hold at most once, namely for $q_{11}+q_{12}$. A final substitution provides us with the claimed equation.
    \end{sketchproof}
\begin{example}\label{example_hittingbdry}
Consider two elliptic curves $\mathbb{T}E$ and $\mathbb{T}E'$ with $l_{\mathbb{T}E}=5$ and $l_{\mathbb{T}E'}=3$. We are looking for subgroups $G$ of the direct product $\mathbb{T}E' \oplus \mathbb{T}E$ that arise from isomorphisms between their $d$-torsions. These should be as required by Lemma \ref{lemma_bdrycharacterizationforfamilies} and \ref{lemma_bdrycharacterizationforfamilies_d-1} of the form (here in coordinates, i.e. as map $\mathbb{R}/3\mathbb{Z}[d]\rightarrow \mathbb{R}/5\mathbb{Z}[d] $):
\begin{align}
    \frac{3}{d} \mapsto \frac{5}{d} \thinspace \text{ (as in Lemma \ref{lemma_bdrycharacterizationforfamilies}) and } 
    \frac{(d-1)3}{d} \mapsto \frac{5}{d} \thinspace \text{( as in Lemma \ref{lemma_bdrycharacterizationforfamilies_d-1}) }
\end{align}
and yield for each $d$ a pptav $J^{pp}$ (Constrution \ref{construction_determinesplitting}) that is the Jacobian of a family of curves whose combinatorial type is the dumbbell-graph. According to the previous discussion consider 
 \begin{align}
    a\cdot 5= (d-a)\cdot 3,
\end{align}   
 which is satisfied, for example, by the tuple $(d=16,a=6)$ or $(d=24,a=9)$. Running through Algorithm \ref{algorithm_PreimageofTormap} we find that the special fiber of the corresponding family of genus $2$ curves is defined by $(l(e_1)=\frac{1}{2}$,$l(e_2)=30)$, for the tuple $(d=16,a=6)$, and $(l(e_1)=\frac{1}{3}$,$l(e_2)=45)$, for the tuple $(d=24,a=9)$. The general fibers are as in the proof of Proposition \ref{proposition_fibreofTorelliMap}. 
\end{example}
Example \ref{example_hittingbdry} illustrates: Geometrically, the subfamilies of pptavs from Lemma \ref{lemma_bdrycharacterizationforfamilies} and \ref{lemma_bdrycharacterizationforfamilies_d-1} correspond to fixing a special class of isomorphisms between the $d$-torsion points of $\mathbb{T}E'$ and $\mathbb{T}E$. 
The reason for explicit results in these cases is simply that the (type of) reduced representative in $\sigma$ that is produced by Algorithm \ref{algorithm_PreimageofTormap} is fixed. For a general $Q^{pp} \in \mathcal{Q}$ we only know that there exists a sequence $(\alpha_1,\beta_1,...,\alpha_n,\beta_n)$ of positive integers such that
\begin{align}
   \begin{pmatrix}
            1 & \beta_n \\
            0 & 1 
        \end{pmatrix}\begin{pmatrix}
            1 & 0 \\
            \alpha_n & 1 
        \end{pmatrix}\cdots \begin{pmatrix} 1 & \beta_1 \\
            0 & 1 
        \end{pmatrix}\begin{pmatrix}
            1 & 0 \\
            \alpha_1 & 1 
        \end{pmatrix} \sbullet Q^{pp} \in \sigma .
\end{align}
If we want to extend our previous approach, the question arises whether there exist other subfamilies of $\mathcal{Q}$ that restrict the possibilities for $(\alpha_1,\beta_1,...,\alpha_n,\beta_n)$ (e.g. that fix $n$). In other words whether $(\alpha_1,\beta_1,...,\alpha_n,\beta_n)$ has geometric meaning.
\begin{example}
In the vein of previous results, we aim at classes of isomorphisms whose graphs are the subgroups $G$ (Lemma \ref{lemma_reqfordiagrwithexactdiag}). But more conservatively, work with fixed $d=5$ and $k=2$ first. The result is sobering: 
For $l_{\mathbb{T}E}=2$ and $l_{\mathbb{T}E'}=3$, for example, we have type $(2,0)$, and for $l_{\mathbb{T}E}=10$ and $l_{\mathbb{T}E'}=100$, type $(3,1)$. 
\end{example}
Hence, from now on we will fix $(\alpha_1,\beta_1,...,\alpha_n,\beta_n)$ as additional, purely algebraic datum and split the initial question about the intersection $\mathcal{Q} \cap \partial A_2^{tr} $:
\begin{enumerate}
    \item Given $G$ (i.e. given $k,d$), can we find $\mathbb{T}E$ and $\mathbb{T}E'$ such that $\gls{Jpp}$ as constructed in \ref{construction_determinesplitting} is equal to $\Jac(\mathcal{C})$ (see Proposition \ref{proposition_fibreofTorelliMap}, case (2))?
    \item Given $\mathbb{T}E$ and  $\mathbb{T}E'$, can we find $G$ such that $\gls{Jpp}=\Jac(C_t)$ for each member $C_t$ of $\mathcal{C}$?
\end{enumerate}
\paragraph{Question (1).} The following Lemma can be viewed as a generalization of Lemma \ref{lemma_bdrycharacterizationforfamilies} and \ref{lemma_bdrycharacterizationforfamilies_d-1}, assuming that the representative is fixed (!).
\begin{lemma}\label{lemma_question1inmodulispaceperspective}
Let $1\leqslant k\leqslant d$ be coprime and $(\alpha_1,\beta_1,...,\alpha_n,\beta_n)\in \mathbb{N}\setminus \{0\}^{2n-1} \cup \mathbb{N}$.
    Then there exists elliptic curves $\mathbb{T}E$ and $\mathbb{T}E'$ with $\gls{Jpp}=\Jac(\mathcal{C})$ as in case (2) of Proposition \ref{proposition_fibreofTorelliMap}, if one of the following two linear equations has a solution $(l_{\mathbb{T}E'},l_{\mathbb{T}E})\in \mathbb{R}_{>0}$:
    \begin{align}
    c^j_1(\alpha_i,\beta_i,k,d)l_{\mathbb{T}E'} + c^j_2(\alpha_i,\beta_i)l_{\mathbb{T}E} =0, \thinspace j=1,2,
\end{align} 
where the coefficients $c^j_1$ and $c^j_2$ for $j=1,2$ are fixed by $G$ (i.e. by $k,d$) and $(\alpha_1,\beta_1,...,\alpha_n,\beta_n)$.
\end{lemma}
\begin{proof}
    Let $Q:=\begin{pmatrix}
            1 & \beta_n \\
            0 & 1 
        \end{pmatrix}\begin{pmatrix}
            1 & 0 \\
            \alpha_n & 1 
        \end{pmatrix}\cdots \begin{pmatrix} 1 & \beta_1 \\
            0 & 1 
        \end{pmatrix}\begin{pmatrix}
            1 & 0 \\
            \alpha_1 & 1 
        \end{pmatrix} \sbullet Q^{pp}.$ We claim that $Q\in \partial \sigma$, whenever one of the following holds:
        \begin{diagram}
        \label{diagram_insection62equations}
            q_{11} + q_{12} = 0 \text{ or } q_{22} + q_{12} = 0.
        \end{diagram}
        Suppose $q_{11} + q_{12} = 0$. We have to show that the entries of $Q$ satisfy the remaining inequalities from Remark \ref{remark_QliesinbdrywhenSellingParametervanishes}. Note that $q_{11}>0$ since $q^{pp}_{11}$ is and, given that transformations of type $\begin{pmatrix} 
            1 & \beta \\
            0 & 1 
        \end{pmatrix}$ and $\begin{pmatrix}
            1 & 0 \\
            \alpha & 1 
        \end{pmatrix}$ keep one of the diagonal entries of $Q^{pp}$ fixed, the other does not change sign as well.\\
        Exemplary for 
        $ \tilde{Q}:=\begin{pmatrix}
            1 & 0 \\
            \alpha & 1 
        \end{pmatrix} \sbullet Q^{pp}$:
        We have $\tilde{q}_{11}=q^{pp}_{11} + \alpha q^{pp}_{12}$ and $\tilde{q}_{22}=q^{pp}_{22}$. Then $0 < \tilde{q}_{11}\tilde{q}_{22}-(\tilde{q}_{12})^2$ ($\tilde{Q}$ is positive definite) implies that $\tilde{q}_{11}$ and $\tilde{q}_{22}$ have the same sign.

        This shows $q_{12} = -q_{11} < 0$ and $q_{22} + q_{12}= q_{22} - q_{11} >0$ follows from  $0 < q_{11}q_{22}-(q_{12})^2= q_{11}q_{22}-(q_{11})^2$.
        Proceeding analogously for the case $q_{22} + q_{12}=0$ proves the claim.

        To obtain equations of the form given in Lemma \ref{lemma_question1inmodulispaceperspective}, consider, equivalently to Equations (\ref{diagram_insection62equations}), their $d$-multiple. It then follows from the structure of $Q^{pp}$ that these are linear in $l_{\mathbb{T}E'}$ and $l_{\mathbb{T}E}$, i.e. with no constant term: 
        By acting with $\begin{pmatrix} 
            1 & \beta \\
            0 & 1 
        \end{pmatrix}$ and $\begin{pmatrix}
            1 & 0 \\
            \alpha & 1 
        \end{pmatrix}$ on $Q^{pp}$ we simply swap the entries of $Q^{pp}$, which are linear functions in $l_{\mathbb{T}E'}$ and $l_{\mathbb{T}E}$, for linear combinations of them.
\end{proof}
\begin{remark}\label{remark_insection6existenceofsolutiondependsonsigncoeff}
\ \\
\begin{itemize}
    \item The existence of a solution to the equations in Lemma \ref{lemma_question1inmodulispaceperspective} is equivalent to $c^j_1$ and $c^j_2$ having different sign. A closer look at both reveals that $c^j_2>0$. The sign of $c^j_1$, however, depends on $G$ and $(\alpha_1,\beta_1,...,\alpha_n,\beta_n)$.
    \item The proof shows that Equations (\ref{diagram_insection62equations}) cannot be satisfied simultaneously.
\end{itemize}
\end{remark}
\paragraph{Question (2).}
Unlike Question (1), answering Question (2) (amounts to/involves) finding a special solution to a quadratic equation. By considering
\begin{align}
     q_{11} + q_{12} = 0 \text{ or } q_{22} + q_{12} = 0
\end{align}
as equations in $d$ and $k$ we get:
\begin{align}
    c^j_1(\alpha_i,\beta_i)l_{\mathbb{T}E'}\cdot k^2 + c^j_2(\alpha_i,\beta_i)l_{\mathbb{T}E'}\cdot kd + c^j_3(\alpha_i,\beta_i)l_{\mathbb{T}E'}\cdot d^2 + c^j_4(\alpha_i,\beta_i)l_{\mathbb{T}E} =0, 
\end{align} 
where $c^j_i\in \mathbb{Z}$ with $c^j_1,c^j_3,c^j_4>0$ and $c^j_2 < 0$ for $j=1,2$ are fixed by $(\alpha_1,\beta_1,...,\alpha_n,\beta_n)$. Then $G$ exists, if one of the above has a solution $(d,k)$ with $d\in \mathbb{N}$ and $k\in \mathbb{Z}^*_d$. 
\subsection{Locus of curves with split Jacobian}\label{subsection_LocusofcurveswithsplitJacobian}
Subsections \ref{subsubssection_step2concrete} and \ref{subsection_Schottky-typeProblem} give us a way to understand and organize the locus of curves with $d$-split Jacobians
that follows naturally from previous considerations.
\begin{construction}\label{construction_fanDeltak}
The locus $\mathbb{T}\mathcal{L}_d:=(t_2^{tr})^{-1}(\mathcal{Q}_d)$ decomposes into 
\begin{align}
   \mathbb{T}\mathcal{L}_d=\bigcup_{k \in \mathbb{Z}^*_d} L_k,
\end{align}
where $L_k:=(t_2^{tr})^{-1}(\mathcal{Q}_{d,k})$ and $\mathcal{Q}_d\subset \mathcal{Q}$ and $\mathcal{Q}_{d,k}\subset \mathcal{Q}$ are the subsets obtained by fixing $d$, respectively $d$ and $k$. Now, $\mathcal{Q}_{d,k}$ is the disjoint union of $\mathcal{Q}_{d,k} \cap \partial \mathcal{A}^{tr}_2$ and $\mathcal{Q}_{d,k}\setminus (\mathcal{Q}_{d,k} \cap \partial \mathcal{A}^{tr}_2)$.

By Lemma \ref{lemma_question1inmodulispaceperspective} we know that $\mathcal{Q}_{d,k} \cap \partial \mathcal{A}^{tr}_2$ is in 1:1 correspondence with a union of half rays:
\begin{align}
    \bigcup_{(\alpha_i,\beta_i)}\{ (l_{\mathbb{T}E'},l_{\mathbb{T}E}) \in \mathbb{R}^2_{>0} \thinspace| \thinspace c_1(\alpha_i,\beta_i,k,d)l_{\mathbb{T}E'} + c_2(\alpha_i,\beta_i,k,d)l_{\mathbb{T}E} =0\},
\end{align}

where $(\alpha_i,\beta_i)$ are so that $c_1(\alpha_i,\beta_i,k,d)>0$ (Remark \ref{remark_insection6existenceofsolutiondependsonsigncoeff}). In this case we say $(\alpha_i,\beta_i)$ is \emph{feasible} (for $d$ and $k$).

We define the fan $\Delta_k$ associated to $L_k$ as the collection of cones obtained by subdividing $\mathbb{R}^2_{>0}$ as above and taking the Euclidean closure.
\end{construction}
\begin{theorem}\label{theorem_locusofcurveswithdsplitJac}
   The locus of curves with $d$-split Jacobian decomposes into $\varphi(d)$ subsets
    \begin{align}
        \mathbb{T}\mathcal{L}_d=\bigcup_{k\in \mathbb{Z}^*_d}L_k,
    \end{align}
  where $\varphi(d)$ denotes the Euler phi function. 
  The fan $\Delta_k$ (see Construction \ref{construction_fanDeltak}) has maximal cones $\sigma_{(\alpha_i,\beta_i)}$ indexed by feasible types and relates to $L_k$ via a family of linear maps $\{ \phi_{\sigma_{(\alpha_i,\beta_i)}} : \sigma_{(\alpha_i,\beta_i)} \rightarrow M^{tr}_2 \}$ such that
  \begin{itemize}
      \item  $\phi_{\tilde{\sigma}}$  
      maps the interior of $\tilde{\sigma}:=\sigma_{(\alpha_i,\beta_i)}$ to a 2 dimensional cone in the theta part of $ M^{tr}_2$.
      \item restricting $\phi_{\tilde{\sigma}}$ to a ray $r$ that is a face of two maximal cones yields a linear map $\phi_{\tilde{\sigma} |r}\times id :r \times \mathbb{R}_{\geqslant 0} \rightarrow M^{tr}_2$ whose image in the dumbbell part of $ M^{tr}_2$ is $\phi_{\tilde{\sigma}}(r) \times r_B$, where $r_B$ denotes the ray that corresponds to the bridge edge.
      \item the maps in $\{\phi_{\sigma_{(\alpha_i,\beta_i)}}\}$ are compatible in the sense that for two neighboring maximal cones $\tilde{\sigma}$ and $\tilde{\sigma}'$ the restriction of $\phi_{\tilde{\sigma}}$, respectively $\phi_{\tilde{\sigma}'}$, to a common face agree.
  \end{itemize}
\end{theorem}
\begin{proof}
   
We address the structure of $\Delta_k$ first and then describe its realization in $M^{tr}_2$.

\emph{1. Maximal cones are indexed by feasible types:}
 
A point $(l_{\mathbb{T}E'},l_{\mathbb{T}E})$ in the interior of a maximal cone $\tilde{\sigma}\in \Delta_k$ corresponds to $Q^{pp}(l_{\mathbb{T}E'},l_{\mathbb{T}E})\in (\mathcal{Q}_{d,k}\setminus (\mathcal{Q}_{d,k} \cap \partial \mathcal{A}^{tr}_2))$, i.e. we have 
\begin{diagram}\label{diagram_equationinproofTLd}
    q_{12}<0, - q_{12} - q_{11}<0, - q_{12} - q_{22}<0,
\end{diagram}
where  $Q:=(q_{ij}):=(\alpha_i,\beta_i)\sbullet Q^{pp}(l_{\mathbb{T}E'},l_{\mathbb{T}E})\in \sigma$ denotes the representative (in $\sigma$) obtained by Selling's reduction algorithm (see Subsection \ref{subsubssection_step2concrete}). Since inequalities (\ref{diagram_equationinproofTLd}) depend continuously on $(l_{\mathbb{T}E'},l_{\mathbb{T}E})$, points in a small neighborhood of $(l_{\mathbb{T}E'},l_{\mathbb{T}E})$ yield pp of type $(\alpha_i,\beta_i)$ as well. This is true as long as inequalities (\ref{diagram_equationinproofTLd}) are preserved.

If they are not, then $Q(l_{\mathbb{T}E'},l_{\mathbb{T}E})\in \partial A^{tr}_2$ and hence $Q^{pp}(l_{\mathbb{T}E'},l_{\mathbb{T}E})\in \partial A^{tr}_2$. 
More precisely: If $- q_{12} - q_{11}=0$ or  $- q_{12} - q_{22}=0$ holds, then $Q^{pp}(l_{\mathbb{T}E'},l_{\mathbb{T}E})$ has type $(\alpha_i,\beta_i)$ as well and $(l_{\mathbb{T}E'},l_{\mathbb{T}E})$ lies on the ray that corresponds to $(\alpha_i,\beta_i)$. If $q_{12}=0$, then $(l_{\mathbb{T}E'},l_{\mathbb{T}E})$ is a point on the ray defined by $(\tilde{\alpha}_i,\tilde{\beta}_i)$, where $(\tilde{\alpha}_i,\tilde{\beta}_i)$ differs from $(\alpha_i,\beta_i)$ by an $\alpha$ or by a $\beta$ transformation: $q_{12}=0$ does not occur when performing Selling's reduction algorithm, but $Q\in \sigma$. We compare inequalities (\ref{diagram_equationinproofTLd}) before and after an $\alpha$, respectively $\beta$ transformation:

\begin{align}
    & q^{after}_{12}=-q^{before}_{23}, q^{after}_{13}=q^{before}_{13} + 2 q^{before}_{23}, q^{after}_{23}=2q^{before}_{23} + q^{before}_{12} \\
    & q^{after}_{12}=-q^{before}_{13}, q^{after}_{13}=q^{before}_{12} + 2 q^{before}_{13}, q^{after}_{23}=q^{before}_{23} + 2 q^{before}_{13}
\end{align}
We see $q^{after}_{12}=0$ means, either $q^{before}_{23}=0$, or $q^{before}_{13}=0$. Selling's Algorithm has to terminate earlier.

\emph{2. Construction of } $\phi_{\tilde{\sigma}}$:
For each maximal cone $\tilde{\sigma}$ with corresponding type $(\alpha_i,\beta_i)$ write $Q(l_{\mathbb{T}E'},l_{\mathbb{T}E}):=(\alpha_i,\beta_i)\sbullet Q^{pp}(l_{\mathbb{T}E'},l_{\mathbb{T}E})$ as linear combination of the extreme rays of $\sigma$
    \begin{align}
      l_1(l_{\mathbb{T}E'},l_{\mathbb{T}E}) \begin{pmatrix}
          1 & 0\\ 0 & 0
      \end{pmatrix}+ l_2(l_{\mathbb{T}E'},l_{\mathbb{T}E})\begin{pmatrix}
          0 &0\\ 0 & 1
      \end{pmatrix}+ l_3(l_{\mathbb{T}E'},l_{\mathbb{T}E})\begin{pmatrix}
          1 &-1\\ -1 & 1
      \end{pmatrix},
    \end{align}
    and set $\phi_{\tilde{\sigma}}(l_{\mathbb{T}E'},l_{\mathbb{T}E}):=(l_1,l_2,l_3)$. Note that $\phi_{\tilde{\sigma}}: \tilde{\sigma} \rightarrow M^{tr}_2$ is linear since
   
    \begin{align}
     l_1= q_{12} + q_{11}, l_2= q_{12} + q_{22}, l_3=- q_{12}.
    \end{align}   
\emph{3. Compatibility:} Combining Points 1 and 2, we see 
\begin{align}
    \phi_{\tilde{\sigma}}=\begin{pmatrix}
        q^{after}_{12} + q^{after}_{11}\\
        q^{after}_{12} + q^{after}_{22}\\
        -q^{after}_{12}
    \end{pmatrix} \text{ and } \phi_{\tilde{\sigma}'}=\begin{pmatrix}
        q^{before}_{12} + q^{before}_{11}\\
        q^{before}_{12} + q^{before}_{22}\\
        -q^{before}_{12}
    \end{pmatrix}
\end{align}
for cones $\tilde{\sigma}$ and $\tilde{\sigma}'$ having a common face $r$. As $r$ is defined by $q^{before}_{12} + q^{before}_{11}=0$ or $q^{before}_{12} + q^{before}_{22}=0$, $\phi_{\tilde{\sigma}|r}$ and $\phi_{\tilde{\sigma}'|r}$ agree up to permutation. 

\emph{4. The interior of }$\tilde{\sigma}$: $\phi_{\tilde{\sigma}}(\tilde{\sigma})$ is a cone that maps the interior of $\tilde{\sigma}$ to the theta part of $M^{tr}_2$. This is an immediate consequence of Point 2.
We show that $\phi_{\tilde{\sigma}}(\tilde{\sigma})$ is 2 dimensional: First, suppose $\mathbb{R}_{\geqslant 0}\cdot l_{\mathbb{T}E'}$ is not a face of $\tilde{\sigma}$
and let $\gamma$ be a path between the two extreme rays of $\tilde{\sigma}$ with constant $l_{\mathbb{T}E}$ coordinate (see Figure \ref{figure_inproofTLd}). Given $(l_{\mathbb{T}E'},l_{\mathbb{T}E}),(\tilde{l}_{\mathbb{T}E'},\tilde{l}_{\mathbb{T}E})\in \tilde{\sigma}$, observe that $Q^{pp}(l_{\mathbb{T}E'},l_{\mathbb{T}E})\sim Q^{pp}(\tilde{l}_{\mathbb{T}E'},\tilde{l}_{\mathbb{T}E})$ in $A^{tr}_2$ implies $\tilde{l}_{\mathbb{T}E'}\cdot\tilde{l}_{\mathbb{T}E}=l_{\mathbb{T}E'}\cdot l_{\mathbb{T}E}$ (by taking determinants). Then $\phi_{\tilde{\sigma}|\gamma}$ is injective and for a shift $\gamma_\epsilon:=\gamma+(0,\epsilon)$, where $\epsilon>0$ is big enough, we have $\phi_{\tilde{\sigma}}(\im(\gamma))\cap \phi_{\tilde{\sigma}}(\im(\gamma_\epsilon))=\emptyset$ by the same argument. This proves that $\phi_{\tilde{\sigma}}(\tilde{\sigma})$ is 2 dimensional. Now, suppose $\mathbb{R}_{\geqslant 0} \cdot l_{\mathbb{T}E'}$ is a face of $\tilde{\sigma}$ and consider a path with constant $l_{\mathbb{T}E'}$ coordinate instead.
\begin{figure}[H]
    \centering
    \tikzset{every picture/.style={line width=0.75pt}} 

\begin{tikzpicture}[x=0.75pt,y=0.75pt,yscale=-1,xscale=1]

\draw  (148.1,235.1) -- (345.43,235.1)(167.83,57.8) -- (167.83,254.8) (338.43,230.1) -- (345.43,235.1) -- (338.43,240.1) (162.83,64.8) -- (167.83,57.8) -- (172.83,64.8)  ;
\draw [color={rgb, 255:red, 74; green, 144; blue, 226 }  ,draw opacity=1 ]   (167.83,235.1) -- (250.33,75) ;
\draw [color={rgb, 255:red, 74; green, 74; blue, 74 }  ,draw opacity=1 ][fill={rgb, 255:red, 155; green, 155; blue, 155 }  ,fill opacity=0.64 ]   (250.33,75) -- (167.83,235.1) -- (320.33,159) ;
\draw [color={rgb, 255:red, 0; green, 0; blue, 0 }  ,draw opacity=1 ][fill={rgb, 255:red, 74; green, 74; blue, 74 }  ,fill opacity=0.26 ]   (167.83,235.1) -- (320.33,159) ;
\draw    (224.33,207.67) -- (224.33,123.67) ;
\draw [shift={(224.33,123.67)}, rotate = 270] [color={rgb, 255:red, 0; green, 0; blue, 0 }  ][fill={rgb, 255:red, 0; green, 0; blue, 0 }  ][line width=0.75]      (0, 0) circle [x radius= 3.35, y radius= 3.35]   ;
\draw [shift={(224.33,159.67)}, rotate = 90] [color={rgb, 255:red, 0; green, 0; blue, 0 }  ][line width=0.75]    (10.93,-4.9) .. controls (6.95,-2.3) and (3.31,-0.67) .. (0,0) .. controls (3.31,0.67) and (6.95,2.3) .. (10.93,4.9)   ;
\draw [shift={(224.33,207.67)}, rotate = 270] [color={rgb, 255:red, 0; green, 0; blue, 0 }  ][fill={rgb, 255:red, 0; green, 0; blue, 0 }  ][line width=0.75]      (0, 0) circle [x radius= 3.35, y radius= 3.35]   ;
\draw    (244.08,197.05) -- (243.08,87.05) ;
\draw [shift={(243.08,87.05)}, rotate = 269.48] [color={rgb, 255:red, 0; green, 0; blue, 0 }  ][fill={rgb, 255:red, 0; green, 0; blue, 0 }  ][line width=0.75]      (0, 0) circle [x radius= 3.35, y radius= 3.35]   ;
\draw [shift={(243.53,136.05)}, rotate = 89.48] [color={rgb, 255:red, 0; green, 0; blue, 0 }  ][line width=0.75]    (10.93,-4.9) .. controls (6.95,-2.3) and (3.31,-0.67) .. (0,0) .. controls (3.31,0.67) and (6.95,2.3) .. (10.93,4.9)   ;
\draw [shift={(244.08,197.05)}, rotate = 269.48] [color={rgb, 255:red, 0; green, 0; blue, 0 }  ][fill={rgb, 255:red, 0; green, 0; blue, 0 }  ][line width=0.75]      (0, 0) circle [x radius= 3.35, y radius= 3.35]   ;

\draw (322,245.4) node [anchor=north west][inner sep=0.75pt]    {$l_{\mathbb{T} E}$};
\draw (112,56.4) node [anchor=north west][inner sep=0.75pt]    {$l_{\mathbb{T} E'}$};
\draw (322.33,162.4) node [anchor=north west][inner sep=0.75pt]  [color={rgb, 255:red, 0; green, 0; blue, 0 }  ,opacity=1 ]  {$r_{1}$};
\draw (243,49.4) node [anchor=north west][inner sep=0.75pt]  [color={rgb, 255:red, 0; green, 0; blue, 0 }  ,opacity=1 ]  {$r_{2}$};
\draw (205,171.4) node [anchor=north west][inner sep=0.75pt]    {$\gamma $};
\draw (257,140.4) node [anchor=north west][inner sep=0.75pt]    {$\gamma _{\epsilon }$};

\end{tikzpicture}

    \caption{Paths $\gamma$ and $\gamma_\epsilon$ in the proof of Theorem \ref{theorem_locusofcurveswithdsplitJac}}
    \label{figure_inproofTLd}
\end{figure}
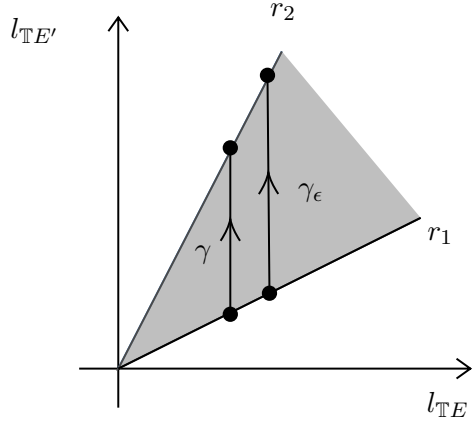
\emph{5. The boundary of } $\tilde{\sigma}$: If $(l_{\mathbb{T}E'},l_{\mathbb{T}E})$ lies on a common face $r$ of two maximal cones, then $l_{\mathbb{T}E'}$ and $l_{\mathbb{T}E}$ are both non-zero and give rise to a pp that lives in the boundary of $\mathcal{A}^{tr}_2$. Then, by Proposition \ref{proposition_fibreofTorelliMap}, its preimage under $t^{tr}_2$ consists of a family of tropical curves of type dumbbell whose cycle lengths are fixed by the non-zero entries of $\phi_{\tilde{\sigma}}(l_{\mathbb{T}E'},l_{\mathbb{T}E})$. This induces a map
\begin{align}
    \phi_{\tilde{\sigma} |r}\times id :r \times \mathbb{R}_{\geqslant 0} \rightarrow M^{tr}_2,
\end{align}
whose image in the dumbbell part of $ M^{tr}_2$ is $\phi_{\tilde{\sigma}}(r) \times r_B$, where $r_B$ denotes the ray that corresponds to the bridge edge.

\end{proof}
\begin{remark}\ \\
\begin{itemize}
    \item The maps $\phi_{\sigma_{(\alpha_i,\beta_i)}}$ can be viewed as local inverses of the Torelli map.
    \item By taking the Euclidean closure in $\mathbb{R}^2$ of the cones induced by subdividing $\mathbb{R}^2_{>0}$ according to feasible types, we get points that do not fit into the setting of Section \ref{section_reconstruction}, namely points whose $l_{\mathbb{T}E}$- or $l_{\mathbb{T}E'}$-coordinate is $0$. Nevertheless, these give insight into the behavior of the reconstruction procedure in the limit $l({\mathbb{T}E}) \rightarrow 0$ or $l(\mathbb{T}E')\rightarrow 0$ for the inputs $(\mathbb{T}E,\mathbb{T}E',G)$. The corresponding genus 2 curve tends towards a curve whose combinatorial type is a cycle with a vertex of genus 1.
\end{itemize}
\end{remark}
\begin{remark}\label{remark_explicitcharacterizationL1andLd-1}
    With Lemma \ref{lemma_bdrycharacterizationforfamilies} and \ref{lemma_bdrycharacterizationforfamilies_d-1} we have an explicit description of the fans $\Delta_1$ and $\Delta_{d-1}$ associated to $L_1$ and $L_{d-1}$: For $\Delta_1$ we have a collection of $d$ cones obtained by subdividing $\mathbb{R}^2_{\geqslant 0}$ as follows:
    \begin{align}
        \bigcup^{d-1}_{\alpha=1}\{ (l_{\mathbb{T}E'},l_{\mathbb{T}E}) \in \mathbb{R}^2_{>0} \thinspace| \thinspace \alpha l_{\mathbb{T}E} + (\alpha-d)l_{\mathbb{T}E'} =0\}.
    \end{align}
    For $\Delta_{d-1}$ the subdivision is induced by
    \begin{align}
        \bigcup^{d-1}_{\beta=1}\{ (l_{\mathbb{T}E'},l_{\mathbb{T}E}) \in \mathbb{R}^2_{>0} \thinspace| \thinspace \beta l_{\mathbb{T}E} + (\beta-d)l_{\mathbb{T}E'} =0\}
    \end{align}
    and generates $d$ cones.
\end{remark}
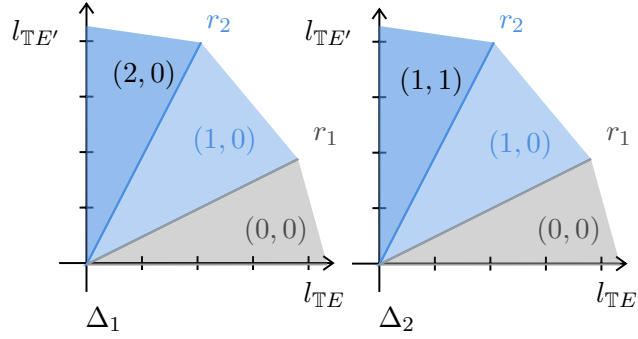
\begin{figure}
    \centering
    \tikzset{every picture/.style={line width=0.75pt}} 

\begin{tikzpicture}[x=0.75pt,y=0.75pt,yscale=-0.7,xscale=0.7]
\tikzset{every picture/.style={line width=0.75pt}} 

\draw  (339.1,217.42) -- (536.43,217.42)(358.83,34) -- (358.83,237.8) (529.43,212.42) -- (536.43,217.42) -- (529.43,222.42) (353.83,41) -- (358.83,34) -- (363.83,41) (398.83,212.42) -- (398.83,222.42)(438.83,212.42) -- (438.83,222.42)(478.83,212.42) -- (478.83,222.42)(518.83,212.42) -- (518.83,222.42)(353.83,177.42) -- (363.83,177.42)(353.83,137.42) -- (363.83,137.42)(353.83,97.42) -- (363.83,97.42)(353.83,57.42) -- (363.83,57.42) ;
\draw   ;
\draw [color={rgb, 255:red, 74; green, 144; blue, 226 }  ,draw opacity=1 ]   (358.83,218.1) -- (441.33,58) ;
\draw [color={rgb, 255:red, 155; green, 155; blue, 155 }  ,draw opacity=1 ][fill={rgb, 255:red, 74; green, 74; blue, 74 }  ,fill opacity=0.24 ]   (511.33,142) -- (358.83,218.1) -- (532.33,218) ;
\draw [color={rgb, 255:red, 74; green, 144; blue, 226 }  ,draw opacity=1 ][fill={rgb, 255:red, 74; green, 144; blue, 226 }  ,fill opacity=0.39 ]   (441.33,58) -- (358.83,218.1) -- (511.33,142) ;
\draw [color={rgb, 255:red, 155; green, 155; blue, 155 }  ,draw opacity=1 ][fill={rgb, 255:red, 74; green, 74; blue, 74 }  ,fill opacity=0.26 ]   (358.83,218.1) -- (511.33,142) ;
\draw  (128.1,217.02) -- (325.43,217.02)(147.83,30) -- (147.83,237.8) (318.43,212.02) -- (325.43,217.02) -- (318.43,222.02) (142.83,37) -- (147.83,30) -- (152.83,37) (187.83,212.02) -- (187.83,222.02)(227.83,212.02) -- (227.83,222.02)(267.83,212.02) -- (267.83,222.02)(307.83,212.02) -- (307.83,222.02)(142.83,177.02) -- (152.83,177.02)(142.83,137.02) -- (152.83,137.02)(142.83,97.02) -- (152.83,97.02)(142.83,57.02) -- (152.83,57.02) ;
\draw   ;
\draw [color={rgb, 255:red, 74; green, 144; blue, 226 }  ,draw opacity=1 ]   (147.83,218.1) -- (230.33,58) ;
\draw [color={rgb, 255:red, 155; green, 155; blue, 155 }  ,draw opacity=1 ][fill={rgb, 255:red, 74; green, 74; blue, 74 }  ,fill opacity=0.24 ]   (300.33,142) -- (147.83,218.1) -- (321.33,218) ;
\draw [color={rgb, 255:red, 74; green, 144; blue, 226 }  ,draw opacity=1 ][fill={rgb, 255:red, 74; green, 144; blue, 226 }  ,fill opacity=0.39 ]   (230.33,58) -- (147.83,218.1) -- (300.33,142) ;
\draw [color={rgb, 255:red, 155; green, 155; blue, 155 }  ,draw opacity=1 ][fill={rgb, 255:red, 74; green, 74; blue, 74 }  ,fill opacity=0.26 ]   (147.83,218.1) -- (300.33,142) ;
\draw [draw opacity=0][fill={rgb, 255:red, 74; green, 144; blue, 226 }  ,fill opacity=0.59 ]   (147.33,46) -- (147.83,218.1) -- (230.33,58) ;
\draw [draw opacity=0][fill={rgb, 255:red, 74; green, 144; blue, 226 }  ,fill opacity=0.59 ]   (358.33,46) -- (358.83,218.1) -- (441.33,58) ;

\draw (356,246.4) node [anchor=north west][inner sep=0.75pt]    {$\Delta_{2}$};
\draw (513,228.4) node [anchor=north west][inner sep=0.75pt]    {$l_{\mathbb{T} E}$};
\draw (303,39.4) node [anchor=north west][inner sep=0.75pt]    {$l_{\mathbb{T} E'}$};
\draw (470,181.4) node [anchor=north west][inner sep=0.75pt]  [color={rgb, 255:red, 74; green, 74; blue, 74 }  ,opacity=1 ]  {$( 0,0)$};
\draw (435,119.4) node [anchor=north west][inner sep=0.75pt]  [color={rgb, 255:red, 74; green, 144; blue, 226 }  ,opacity=1 ]  {$( 1,0)$};
\draw (371,73.4) node [anchor=north west][inner sep=0.75pt]    {$( 1,1)$};
\draw (520,117.4) node [anchor=north west][inner sep=0.75pt]  [color={rgb, 255:red, 74; green, 74; blue, 74 }  ,opacity=1 ]  {$r_{1}$};
\draw (443,36.4) node [anchor=north west][inner sep=0.75pt]  [color={rgb, 255:red, 74; green, 144; blue, 226 }  ,opacity=1 ]  {$r_{2}$};
\draw (145,246.4) node [anchor=north west][inner sep=0.75pt]    {$\Delta_{1}$};
\draw (302,228.4) node [anchor=north west][inner sep=0.75pt]    {$l_{\mathbb{T} E}$};
\draw (92,39.4) node [anchor=north west][inner sep=0.75pt]    {$l_{\mathbb{T} E'}$};
\draw (259,181.4) node [anchor=north west][inner sep=0.75pt]  [color={rgb, 255:red, 74; green, 74; blue, 74 }  ,opacity=1 ]  {$( 0,0)$};
\draw (222,116.4) node [anchor=north west][inner sep=0.75pt]  [color={rgb, 255:red, 74; green, 144; blue, 226 }  ,opacity=1 ]  {$( 1,0)$};
\draw (164,68.4) node [anchor=north west][inner sep=0.75pt]    {$( 2,0)$};
\draw (309,117.4) node [anchor=north west][inner sep=0.75pt]  [color={rgb, 255:red, 74; green, 74; blue, 74 }  ,opacity=1 ]  {$r_{1}$};
\draw (232,36.4) node [anchor=north west][inner sep=0.75pt]  [color={rgb, 255:red, 74; green, 144; blue, 226 }  ,opacity=1 ]  {$r_{2}$};

\end{tikzpicture}

    \caption{The fan $\Delta_1$ on the left and the fan $\Delta_2$ on the right with cones labeled by their type $(\alpha,\beta)$. }
    \label{figure_exampleforfanassociatedtoL1Ld-1}
\end{figure}
\begin{example}\label{example_L1L2agree}
Let $d=3$, then $\mathbb{T}\mathcal{L}_3=L_1 \cup L_2$. By Remark \ref{remark_explicitcharacterizationL1andLd-1} we have an explicit description for both: $\Delta_1$ is induced by
\begin{align}
    \{ (l_{\mathbb{T}E'},l_{\mathbb{T}E}) \in \mathbb{R}^2_{>0} \thinspace| \thinspace l_{\mathbb{T}E} -2 l_{\mathbb{T}E'} =0\} \cup
    \{ (l_{\mathbb{T}E'},l_{\mathbb{T}E}) \in \mathbb{R}^2_{>0} \thinspace| \thinspace 2 l_{\mathbb{T}E} - l_{\mathbb{T}E'} =0\}
\end{align}
and the fan $\Delta_{3-1}$ by:
\begin{align}
    \{ (l_{\mathbb{T}E'},l_{\mathbb{T}E}) \in \mathbb{R}^2_{>0} \thinspace| \thinspace  l_{\mathbb{T}E} - 2 l_{\mathbb{T}E'} =0\}\cup
    \{ (l_{\mathbb{T}E'},l_{\mathbb{T}E}) \in \mathbb{R}^2_{>0} \thinspace| \thinspace 2 l_{\mathbb{T}E} - l_{\mathbb{T}E'} =0\}.
\end{align}
See Figure \ref{figure_exampleforfanassociatedtoL1Ld-1}. 
Following the proof of Theorem \ref{theorem_locusofcurveswithdsplitJac} we compute the maps associated to the maximal cones of $\Delta_1$:

\begin{diagram}\label{diagram_inexampleTL}
    \phi_{(0,0)}=\begin{pmatrix}
        2l_{\mathbb{T}E'} \\ \frac{l_{\mathbb{T}E}}{3}-\frac{2l_{\mathbb{T}E'}}{3}\\ l_{\mathbb{T}E'}
    \end{pmatrix} \thinspace
    \phi_{(1,0)}=\begin{pmatrix}
        \frac{2l_{\mathbb{T}E}}{3}+\frac{2l_{\mathbb{T}E'}}{3} \\ \frac{-l_{\mathbb{T}E'}}{3}+\frac{2l_{\mathbb{T}E}}{3}\\ 
        \frac{-l_{\mathbb{T}E}}{3}+\frac{2l_{\mathbb{T}E'}}{3}
    \end{pmatrix}  \thinspace
    \phi_{(2,0)}=\begin{pmatrix}
        2l_{\mathbb{T}E} \\ l_{\mathbb{T}E} \\ \frac{l_{\mathbb{T}E'}}{3}-\frac{2l_{\mathbb{T}E}}{3}
    \end{pmatrix}
\end{diagram}
Compatibility is verified exemplarily for $\phi_{(0,0)}$ and $\phi_{(1,0)}$: $r_1$ is given by the equation $2l'=l$. Substituting in (\ref{diagram_inexampleTL}) yields
\begin{align}
  \phi_{(0,0)}=\begin{pmatrix}
        2l' \\ 0\\ 
        l'
    \end{pmatrix}  \text{ and } \phi_{(1,0)}= \begin{pmatrix}
        2l' \\ l'\\ 
        0
    \end{pmatrix}.
\end{align}
We conclude with a visualization in Figure \ref{figure_imageofDelta1inM2}. Repeating the same procedure for $\Delta_2$ we get
\begin{align}
    \phi_{(0,0)}=\begin{pmatrix}
        l_{\mathbb{T}E'} \\ \frac{l_{\mathbb{T}E}}{3}-\frac{2l_{\mathbb{T}E'}}{3}\\ 2l_{\mathbb{T}E'}
    \end{pmatrix} \thinspace
    \phi_{(1,0)}=\begin{pmatrix}
    \frac{-l_{\mathbb{T}E'}}{3}+\frac{2l_{\mathbb{T}E}}{3}\\
        \frac{2l_{\mathbb{T}E}}{3}+\frac{2l_{\mathbb{T}E'}}{3} \\  
        \frac{-l_{\mathbb{T}E}}{3}+\frac{2l_{\mathbb{T}E'}}{3}
    \end{pmatrix} \thinspace
    \phi_{(1,1)}=\begin{pmatrix}
        l_{\mathbb{T}E} \\ 2l_{\mathbb{T}E} \\ \frac{l_{\mathbb{T}E'}}{3}-\frac{2l_{\mathbb{T}E}}{3}
    \end{pmatrix}
\end{align}
and see that the image of $\Delta_1$ and $\Delta_2$ in $M_2^{tr}$ agree.
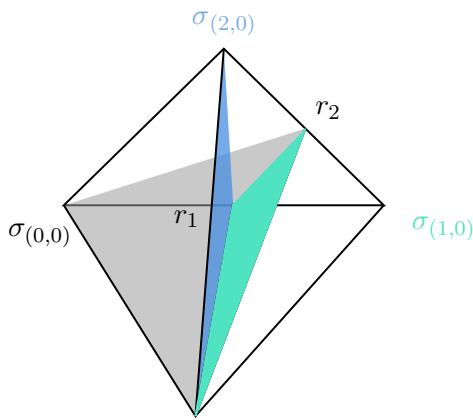
\begin{figure}[H]
    \centering
    \tikzset{every picture/.style={line width=0.75pt}} 

\begin{tikzpicture}[x=0.75pt,y=0.75pt,yscale=-0.76,xscale=0.76]

\draw   (276.17,37) -- (382.33,141) -- (170,141) -- cycle ;
\draw [fill={rgb, 255:red, 155; green, 155; blue, 155 }  ,fill opacity=0.54 ]   (170,141) -- (257.33,281) -- (330.33,90) ;
\draw [color={rgb, 255:red, 0; green, 0; blue, 0 }  ,draw opacity=1 ][fill={rgb, 255:red, 74; green, 144; blue, 226 }  ,fill opacity=0.76 ]   (276.17,37) -- (257.33,281) -- (282.33,140) ;
\draw [line width=0.75]    (382.33,141) -- (257.33,281) ;
\draw [color={rgb, 255:red, 80; green, 227; blue, 194 }  ,draw opacity=1 ][fill={rgb, 255:red, 80; green, 227; blue, 194 }  ,fill opacity=1 ]   (282.33,140) -- (257.33,281) -- (330.33,90) ;

\draw (253,10.4) node [anchor=north west][inner sep=0.75pt]  [color={rgb, 255:red, 74; green, 144; blue, 226 }  ,opacity=0.76 ]  {$\textcolor[rgb]{0.29,0.56,0.89}{\sigma }\textcolor[rgb]{0.29,0.56,0.89}{_{( 2,0)}}$};
\draw (399,146.4) node [anchor=north west][inner sep=0.75pt]  [color={rgb, 255:red, 80; green, 227; blue, 194 }  ,opacity=1 ]  {$\textcolor[rgb]{0.31,0.89,0.76}{\sigma _{( 1,0)}}$};
\draw (131,151.4) node [anchor=north west][inner sep=0.75pt]    {$\sigma _{( 0,0)}$};
\draw (241,142.4) node [anchor=north west][inner sep=0.75pt]    {$r_{1}$};
\draw (335,70.4) node [anchor=north west][inner sep=0.75pt]    {$r_{2}$};

\end{tikzpicture}

    \caption{The image of $\Delta_1$ in $M^{tr}_2$.}
    \label{figure_imageofDelta1inM2}
\end{figure}
\end{example}
\paragraph{Further Questions}
The combinatorial structure associated to $\mathbb{T}\mathcal{L}_d$ (see Theorem \ref{theorem_locusofcurveswithdsplitJacinintro}) organizes $\mathcal{Q}$ according to \emph{Selling's reduction algorithm}. Since we work with a set description of $\mathcal{Q}$ that is not a set of representatives, we expect a certain amount of redundancy. 
\begin{conjecture}
    The fans associated to $L_1$ and $L_{d-1}$ have the same image in $M^{tr}_2$.
\end{conjecture}
The case $d=3$ is verified in \ref{example_L1L2agree}, sporadic examples for higher degree point in the same direction. More generally:
\begin{q}\label{question_onintersectionbeh}
(Can we describe)/(What is) the intersection behavior in $M^{tr}_2$ of the fans associated to $L_k$ for all $k$?
\end{q}
An answer to Question \ref{question_onintersectionbeh} would not only give a more comprehensive understanding of the moduli space structure of $\mathbb{T}\mathcal{L}_d$, but also describe how far our "working description" of $\mathcal{Q}$ is from being a set of representatives. 
\section{Appendix}\label{section_appendix}
 This section follows Subsection \ref{subsection_step2} and the notation therein. We review Selling's reduction algorithm in greater detail (see \cite{thesisVallentin}, Section 2.3.3 and \cite{zbMATH02716590}) and provide an overview of some SINGULAR procedures, which may be used in the context Split Jacobians.

A 2-dimensional positive definite quadratic form $q$ with associated matrix $Q:=(q_{ij})$ is represented by its \emph{Selling parameters}: 
\begin{align}
    q_{12},\thinspace q_{13}:= -q_{11}-q_{12},\thinspace q_{23}:= -q_{11}-q_{12}.
\end{align}
Selling's reduction algorithm reduces $q$ to a form whose Selling parameters are all non-positive: 
Without loss of generality suppose $q_{12}<0$, else consider $\begin{pmatrix}
    1 & 0\\
    0 & -1
\end{pmatrix} \sbullet Q$
and do:
\begin{itemize}
    \item Compute the Selling Parameters of $Q$.
    \item If $q_{13}>0$, set $Q:=\begin{pmatrix}
    1 & 1\\
    0 & 1
\end{pmatrix} \sbullet Q$.
    \item If $q_{23}>0$, set $Q:=\begin{pmatrix}
    1 & 0\\
    1 & 1
\end{pmatrix} \sbullet Q$.
\end{itemize}
Repeat until $q_{13} \leqslant 0$ and $q_{23} \leqslant 0$.

Procedures available in the context of split Jacobians are: 
 \paragraph{\textbf{Procedure:}} \texttt{SELLING}.\\
 Input: A symmetric and positive definite $2 \times 2$ matrix $Q$ with $q_{12}<0$.\\
 Output: Selling parameters of a represantative $\tilde{Q}$ of $Q$ in $\sigma$ and a list of integers that records the number of transformations of type (1) $\begin{pmatrix} 1 & 0\\
    1 & 1
        
    \end{pmatrix}$ and (2) $\begin{pmatrix} 1 & 1\\
    0 & 1 
    \end{pmatrix}$ needed to obtain $\tilde{Q}$.
    \paragraph{\textbf{Procedure:}} \texttt{FD\_VERTETER}.\\
 Input: A symmetric and positive definite $2 \times 2$ matrix $Q$ with non-positive Selling parameters.\\
 Output: A represenative of $Q$ in the fundamental domain $F$ (see Subsection \ref{subsection_step2}).
 \paragraph{\textbf{Procedure:}} \texttt{LENGTH\_OUTPUT}\\
 Input: A matrix $Q\in F$.\\
 Output: A curve $\Gamma$ whose Jacobian is defined by $Q$.

 Now Steps 1 and 2 of Plan \ref{Masterplan} may be retraced computationally: Let $(\mathbb{T}E,\mathbb{T}E',G)$ be splitting data, where $\mathbb{T}E$ and $\mathbb{T}E'$ are elliptic curves of length $1$, respectively $3$, and $G$ is the graph of 
 \begin{align}
     f: \mathbb{R}/3\mathbb{Z}[18] \rightarrow \mathbb{R}/\mathbb{Z}[18], \qquad \frac{7\cdot 3}{18} \longmapsto \frac{1}{18}.
 \end{align}
 \textbf{\emph{Step 1 of Plan \ref{Masterplan}:}} The splitting data $(\mathbb{T}E,\mathbb{T}E',G)$  is recorded as 
 \begin{align}
     (d,k,l_{\mathbb{T}E'},l_{\mathbb{T}E}) = (18, 7, 3, 1).
 \end{align}
 We compute the pptav $J^{pp}$ by running \texttt{SETMATRIX}, which returns:
 \begin{verbatim}
> matrix Q=SETMATRIX(18,7,3,1);
> Q;
Q[1,1]=54
Q[1,2]=-21
Q[2,1]=-21
Q[2,2]=74/9
 \end{verbatim}
  i.e. $Q^{pp}=\begin{pmatrix} 54 & -21\\
    -21 & \frac{74}{9}
    \end{pmatrix}$ (see Subsection \ref{subsection_step2}).  

\textbf{\emph{Step 2 of Plan \ref{Masterplan}:}} Running \texttt{SELLING} and \texttt{RED\_VERTRETER} returns
\begin{verbatim}
> list L=SELLING(Q);L;
[1]:
   -1/3
[2]:
   -11/9
[3]:
   -5/3
[4]:
   1
[5]:
   2
\end{verbatim}
   and 
\begin{verbatim}
> matrix Qn=RED_VERTRETER(L);
> Qn;
Qn[1,1]=26/9
Qn[1,2]=-5/3
Qn[2,1]=-5/3
Qn[2,2]=2
\end{verbatim}
where $Qn$ is a representative of $Q^{pp}$ in $\sigma$ whose Selling parameters are $\frac{-1}{3},\frac{-11}{9},\frac{-5}{3} $ and  $\begin{pmatrix} 1 & 1\\
    0 & 1
    \end{pmatrix},\begin{pmatrix} 1 & 0\\
    2 & 1
    \end{pmatrix}$ is a sequence of transformations that relates $Q^{pp}$ to $Qn$. The unique representative inside $F$ is obtained by invoking 
    \begin{verbatim}
> matrix Qnn=FD_VERTRETER(Qn);Qnn;
Qnn[1,1]=14/9
Qnn[1,2]=-1/3
Qnn[2,1]=-1/3
Qnn[2,2]=2
    \end{verbatim}
    and 
    \begin{verbatim}
> list Ln=LENGTH_OUTPUT(Qnn);Ln;
[1]:
   Jacobian of tropical curve of type T:
[2]:
   [1]:
      11/9
   [2]:
      5/3
   [3]:
      1/3
 \end{verbatim}
    returns the length data $L=\frac{11}{9},\frac{5}{3},\frac{1}{3}$ (see Algorithm \ref{algorithm_PreimageofTormap}) together with the combinatorial type of the curve $\Gamma$ whose Jacobian is $J^{pp}$.
 
\bibliographystyle{plain} 
\bibliography{SplitJacobians.bib}

\begin{thebibliography}{10}

\bibitem{MR3375652}
Omid Amini, Matthew Baker, Erwan Brugall\'e, and Joseph Rabinoff.
\newblock Lifting harmonic morphisms {I}: metrized complexes and {B}erkovich skeleta.
\newblock {\em Res. Math. Sci.}, 2:Art. 7, 67, 2015.

\bibitem{MR2772537}
Matthew Baker and Xander Faber.
\newblock Metric properties of the tropical {A}bel-{J}acobi map.
\newblock {\em J. Algebraic Combin.}, 33(3):349--381, 2011.

\bibitem{MR2525845}
Matthew Baker and Serguei Norine.
\newblock Harmonic morphisms and hyperelliptic graphs.
\newblock {\em Int. Math. Res. Not. IMRN}, (15):2914--2955, 2009.

\bibitem{MR3717092}
Janko B\"ohm, Kathrin Bringmann, Arne Buchholz, and Hannah Markwig.
\newblock Tropical mirror symmetry for elliptic curves.
\newblock {\em J. Reine Angew. Math.}, 732:211--246, 2017.

\bibitem{MR3752493}
Barbara Bolognese, Madeline Brandt, and Lynn Chua.
\newblock From curves to tropical {J}acobians and back.
\newblock In {\em Combinatorial algebraic geometry}, volume~80 of {\em Fields Inst. Commun.}, pages 21--45. Fields Inst. Res. Math. Sci., Toronto, ON, 2017.

\bibitem{MR2739784}
Silvia Brannetti, Margarida Melo, and Filippo Viviani.
\newblock On the tropical {T}orelli map.
\newblock {\em Adv. Math.}, 226(3):2546--2586, 2011.

\bibitem{MR3278571}
Lucia Caporaso.
\newblock Gonality of algebraic curves and graphs.
\newblock In {\em Algebraic and complex geometry}, volume~71 of {\em Springer Proc. Math. Stat.}, pages 77--108. Springer, Cham, 2014.

\bibitem{MR2641941}
Lucia Caporaso and Filippo Viviani.
\newblock Torelli theorem for graphs and tropical curves.
\newblock {\em Duke Math. J.}, 153(1):129--171, 2010.

\bibitem{MR2661417}
Renzo Cavalieri, Paul Johnson, and Hannah Markwig.
\newblock Tropical {H}urwitz numbers.
\newblock {\em J. Algebraic Combin.}, 32(2):241--265, 2010.

\bibitem{MR2968636}
Melody Chan.
\newblock Combinatorics of the tropical {T}orelli map.
\newblock {\em Algebra Number Theory}, 6(6):1133--1169, 2012.

\bibitem{arXiv:2410.13459}
Lou-Jean~Leila Cobigo.
\newblock Tropical {Split} {Jacobians} of genus 2 and optimal covers.
\newblock Preprint, {arXiv}:2410.13459 (2024), 2024.

\bibitem{MR1085258}
Gerhard Frey and Ernst Kani.
\newblock Curves of genus {$2$} covering elliptic curves and an arithmetical application.
\newblock In {\em Arithmetic algebraic geometry ({T}exel, 1989)}, volume~89 of {\em Progr. Math.}, pages 153--176. Birkh\"auser Boston, Boston, MA, 1991.

\bibitem{zbMATH02692043}
E.~Goursat.
\newblock Sur les substitutions orthogonales et les divisions r{\'e}guli{\`e}res de l'espace.
\newblock {\em Ann. Sci. {\'E}c. Norm. Sup{\'e}r. (3)}, 6:9--102, 1889.

\bibitem{MR2457725}
Eric Katz, Hannah Markwig, and Thomas Markwig.
\newblock The {$j$}-invariant of a plane tropical cubic.
\newblock {\em J. Algebra}, 320(10):3832--3848, 2008.

\bibitem{MR4261102}
Yoav Len and Martin Ulirsch.
\newblock Skeletons of {P}rym varieties and {B}rill-{N}oether theory.
\newblock {\em Algebra Number Theory}, 15(3):785--820, 2021.

\bibitem{MR4382460}
Yoav Len and Dmitry Zakharov.
\newblock Kirchhoff's theorem for {P}rym varieties.
\newblock {\em Forum Math. Sigma}, 10:Paper No. e11, 54, 2022.

\bibitem{zbMATH05564780}
Kay Magaard, Tanush Shaska, and Helmut V{\"o}lklein.
\newblock Genus 2 curves that admit a degree 5 map to an elliptic curve.
\newblock {\em Forum Math.}, 21(3):547--566, 2009.

\bibitem{MR2275625}
Grigory Mikhalkin.
\newblock Tropical geometry and its applications.
\newblock In {\em International {C}ongress of {M}athematicians. {V}ol. {II}}, pages 827--852. Eur. Math. Soc., Z\"urich, 2006.

\bibitem{MR2457739}
Grigory Mikhalkin and Ilia Zharkov.
\newblock Tropical curves, their {J}acobians and theta functions.
\newblock In {\em Curves and abelian varieties}, volume 465 of {\em Contemp. Math.}, pages 203--230. Amer. Math. Soc., Providence, RI, 2008.

\bibitem{zbMATH03975095}
J.~S. Milne.
\newblock Abelian varieties.
\newblock Arithmetic geometry, {Pap}. {Conf}., {Storrs}/{Conn}. 1984, 103-150 (1986)., 1986.

\bibitem{MR0282985}
David Mumford.
\newblock {\em Abelian varieties}, volume~5 of {\em Tata Institute of Fundamental Research Studies in Mathematics}.
\newblock Tata Institute of Fundamental Research, Bombay; by Oxford University Press, London, 1970.

\bibitem{MR1341069}
Jo\~ao~F. Queir\'o and Eduardo~M. S\'a.
\newblock Singular values and invariant factors of matrix sums and products.
\newblock {\em Linear Algebra Appl.}, 225:43--56, 1995.

\bibitem{röhrle2024tropicalngonalconstruction}
Felix Röhrle and Dmitry Zakharov.
\newblock The tropical $n$-gonal construction, 2024.

\bibitem{zbMATH02716590}
E.~Selling.
\newblock On binary and ternary quadratic forms.
\newblock {\em J. Reine Angew. Math.}, 77:143--229, 1873.

\bibitem{shaska2024machinelearningmodulispace}
Elira Shaska and Tony Shaska.
\newblock Machine learning for moduli space of genus two curves and an application to isogeny based cryptography, 2024.

\bibitem{zbMATH02074360}
T.~Shaska.
\newblock Genus 2 fields with degree 3 elliptic subfields.
\newblock {\em Forum Math.}, 16(2):263--280, 2004.

\bibitem{thesisVallentin}
Frank Vallentin.
\newblock Sphere coverings, lattices, and tilings (in low dimensions), 2003.

\end{thebibliography}
\end{document}